\documentclass{imsart}

\RequirePackage{amsthm,amsmath,amsfonts,amssymb}
\RequirePackage[numbers,sort&compress]{natbib}
\RequirePackage[colorlinks,citecolor=blue,urlcolor=blue]{hyperref}
\RequirePackage{graphicx}
\RequirePackage{natbib}
\RequirePackage[OT1]{fontenc}
\usepackage{epsfig}
\usepackage{graphicx}
\usepackage{epstopdf,color}
\usepackage{booktabs}
\usepackage{multirow}
\usepackage{listings}
\usepackage{cases}
\usepackage{subcaption}
\usepackage{indentfirst}
\allowdisplaybreaks[4]
\usepackage{euscript}

\startlocaldefs
\theoremstyle{plain}

\newtheorem{theorem}{Theorem}[section]
\newtheorem{lemma}[theorem]{Lemma}
\newtheorem{corollary}{Corollary}[section]
\theoremstyle{remark}
\newtheorem{assumption}{Assumption}
\newtheorem{remark}{Remark}[section]

\allowdisplaybreaks[4]
\newcommand{\be}{{\mathbf e}} 
 
\newcommand{\X}{{\mathbf X}} 
\newcommand{\Y}{{\mathbf Y}} 
\newcommand{\Z}{{\mathbf Z}} 
 
\newcommand{\A}{{\mathbf A}} 
\newcommand{\bB}{{\mathbf B}} 
 
\newcommand{\bK}{{\mathbf K}} 

\newcommand{\bN}{{\mathbf N}} 
\newcommand{\bR}{{\mathbf R}} 
\newcommand{\bS}{{\mathbf S}}
\newcommand{\W}{{\mathbf W}}
\newcommand{\I}{{\mathbf I}}
\newcommand{\s}{{\mathbf s}} 
\newcommand{\e}{{\mathbf e}}
\newcommand{\bx}{{\mathbf x}} 
\newcommand{\y}{{\mathbf y}} 
\newcommand{\z}{{\mathbf z}}

\newcommand{\E}{{\mathbb{E}}}
 
\newcommand{\Cov}{{\rm Cov}} 
 
\newcommand{\bL}{{\mathbf L}}
\newcommand{\bU}{{\mathbf U}} 
 
\newcommand{\bI}{{\mathbf I}} 
\newcommand{\bV}{{\mathbf V}} 
\newcommand{\bM}{{\mathbf M}} 
\newcommand{\bbR}{{\mathbf R}} 
\newcommand{\bQ}{{\mathbf Q}} 
\newcommand{\bH}{{\mathbf H}} 
 
\newcommand{\tr}{\mathop{\text{\rm tr}}}

\newcommand{\bmu}{{\boldsymbol \mu}}

\newcommand{\bSigma}{{\boldsymbol \Sigma}}

\endlocaldefs

\begin{document}

\begin{frontmatter}
\title{On spiked eigenvalues of a renormalized sample covariance matrix from multi-population}
\runtitle{spiked eigenvalues from multi-population}

\begin{aug}
\author[A]{\fnms{Weiming}~\snm{Li}\ead[label=e1]{li.weiming@shufe.edu.cn}\orcid{0000-0002-5181-3437}},
\author[B]{\fnms{Zeng}~\snm{Li}\ead[label=e2]{liz9@sustech.edu.cn}\orcid{0000-0002-3406-5051}}
\and
\author[B]{\fnms{Junpeng}~\snm{Zhu}\ead[label=e3]{zhujp@sustech.edu.cn}\orcid{0000-0003-2528-0662}}
\address[A]{School of Statistics and Management, Shanghai University of Finance and
	Economics\printead[presep={,\ }]{e1}}

\address[B]{Department of Statistics and Data Science, Southern University of Science and Technology\printead[presep={,\ }]{e2,e3}}
\end{aug}

\begin{abstract}
Sample covariance matrices from multi-population typically exhibit several large spiked eigenvalues, which stem from differences between population means and are crucial for inference on the underlying data structure. This paper investigates the asymptotic properties of spiked eigenvalues of a renormalized sample covariance matrices from multi-population in the ultrahigh dimensional context where the dimension-to-sample size ratio \(p/n \to\infty\).  The first- and second-order convergence of these spikes are established based on asymptotic properties of three types of sesquilinear forms from multi-population.  These findings are further applied to two scenarios, including determination of total number of subgroups and a new criterion for evaluating clustering results in the absence of true labels. Additionally, we provide a unified framework  with \(p/n \to c\in (0,\infty]\) that integrates the  asymptotic results in both high and ultrahigh dimensional settings.

\end{abstract}



\end{frontmatter}

	\section{Introduction}

Consider \(n\) independent observations from \(\tau\) populations in \(\mathbb{R}^p\), represented as
\begin{align}\label{xij}
\mathbf{x}_{ij} = \boldsymbol{\mu}_i + \boldsymbol{\Sigma}_0^{\frac{1}{2}} \mathbf{z}_{ij} \triangleq \boldsymbol{\mu}_i + \mathbf{s}_{ij}, \quad i = 1, \ldots, \tau, \ j = 1, \ldots, n_i,
\end{align}
where the population labels \(i\) and the sizes \(n_i\) are  unknown parameters, satisfying \(n_1 + \cdots + n_\tau = n\).
For each \(i\) and \(j\), the random vector \(\mathbf{z}_{ij} \in \mathbb{R}^p\) consists of \(p\) i.i.d. components with zero mean and unit variance. Hence, the \(\tau\) populations share a common covariance matrix \(\boldsymbol{\Sigma}_0 \in \mathbb{R}^{p \times p}\) but different mean vectors \(\{\boldsymbol{\mu}_i\}\). Denote \(\mathbf{X} = (\mathbf{x}_{11}, \ldots, \mathbf{x}_{1n_1}, \ldots, \mathbf{x}_{\tau1}, \ldots, \mathbf{x}_{\tau n_\tau})_{p\times n}\) as the data matrix, then its sample covariance matrix can be written as
\begin{align*}
\mathbf{S}_n = \frac{1}{n} \mathbf{X} \mathbf{\Phi} \mathbf{X}^{\top}
\end{align*}
where \(\mathbf{\Phi} = \mathbf{I}_n - \mathbf{1}_n \mathbf{1}_n^{\top}/n\).
The eigenvalues of \(\mathbf{S}_n\) serve as important statistics and often play crucial roles in the inference on population parameters, see \cite{Anderson}.

For the single population case when \(\tau = 1\), although entrywisely $\mathbf{S}_n$ is a consistent estimator of its population counterpart $\boldsymbol{\Sigma}_0$, the eigenvalues of $\mathbf{S}_n$ exhibit significant deviations from those of $\boldsymbol{\Sigma}_0$ in the high dimensional settings. Consider the following regime,
\begin{align}\label{mp}
n \rightarrow \infty, \ p = p_n \rightarrow \infty, \ p / n \rightarrow c \in (0, \infty),
\end{align}
referred to as the \emph{Mar\v{c}enko-Pastur} (MP) regime.
\cite{MP67,Bai2010a} showed that the empirical spectral distribution of \(\mathbf{S}_n\) converges weakly to a limit known as the MP law. Additionally, \cite{BS04,zheng2015substitution,MR3189088,PZ08} established the \emph{central limit theorem} (CLT) for the \emph{linear spectral statistics} (LSS) of \(\mathbf{S}_n\). Applications of these theories have been extensively discussed, especially in the areas of hypothesis testing, principal component analysis, and signal processing, see \cite{yao2015sample}.

\begin{figure}[htbp]
	\begin{subfigure}{0.45\linewidth}
		\centering
		\includegraphics[width=1\linewidth]{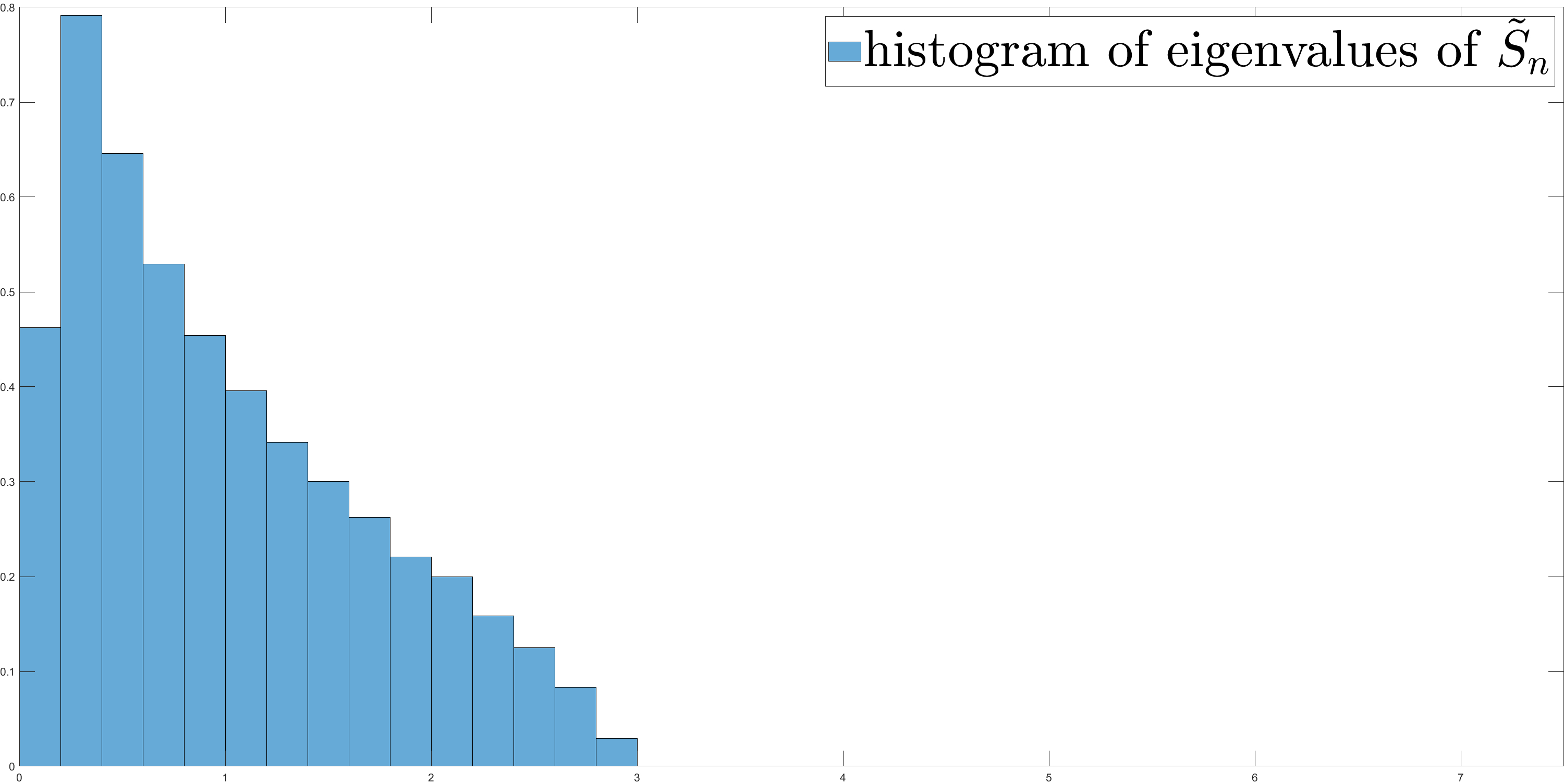}
		\caption{Histogram of eigenvalues of $\tilde\bS_n$.}
	\end{subfigure}
	\begin{subfigure}{0.45\linewidth}
		\centering
		\includegraphics[width=1\linewidth]{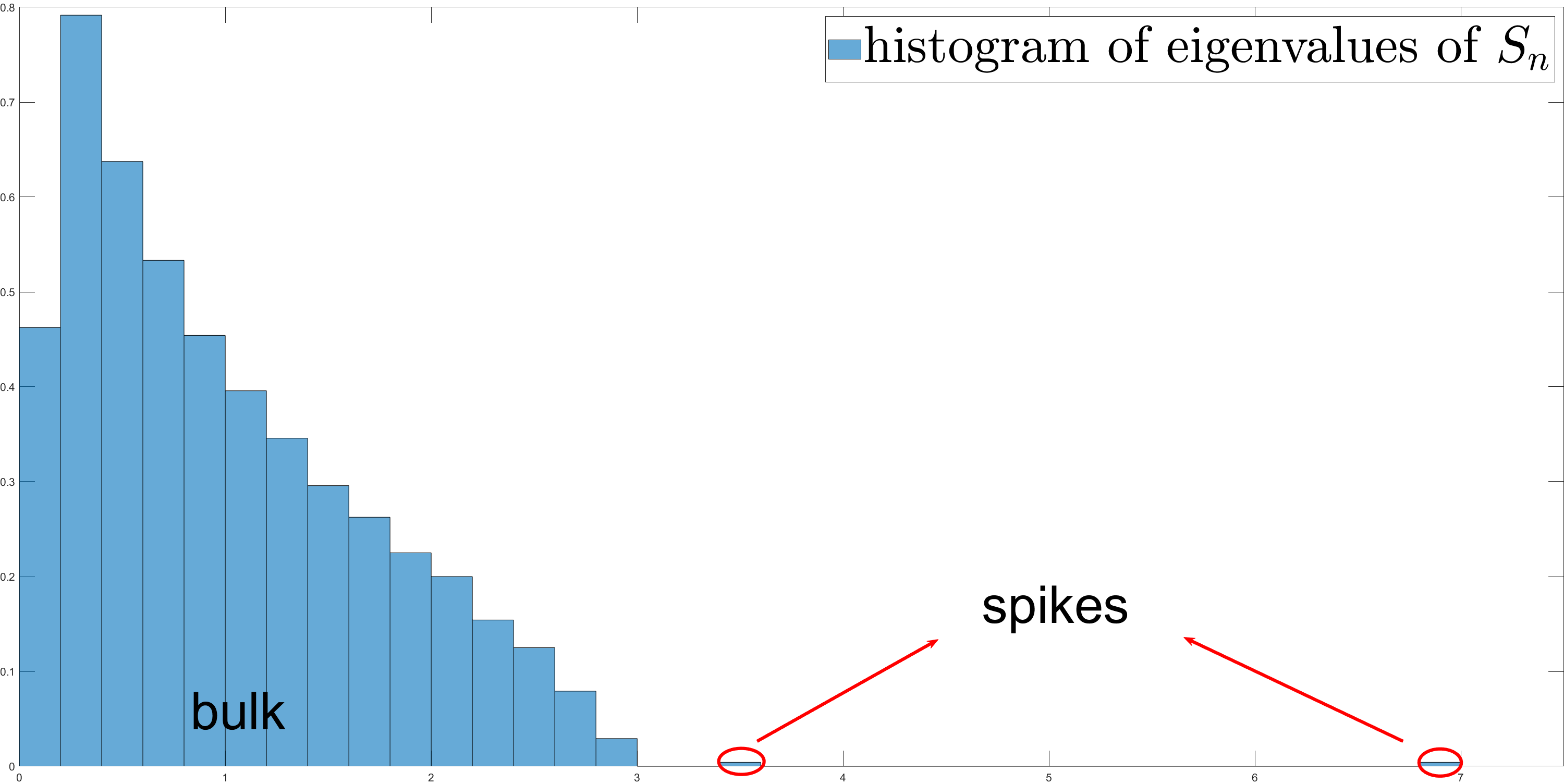}
		\caption{Histogram of eigenvalues of $\bS_n$.}
	\end{subfigure}
	\caption{Histograms of eigenvalues of $\tilde\bS_n$ and $\bS_n$, where $\tau=3$, $\bSigma_0=\bI_p$, $n_1=n_2=n_3=800$, $p=1200$, $\bmu_1=\mathbf{0}$, $\bmu_2=(4,0,\ldots,0)$ and $\bmu_3=(0,4,0,\ldots,0)$.
	}
	\label{intro-Sn}
\end{figure}

For the multi-population case when $\tau>1$,  \(\mathbf{S}_n\) additionally carries information regarding the differences among the subgroups, especially the distances between mean vectors. 
To illustrate, we decompose the matrix into two parts:
\[
\mathbf{S}_n = \tilde{\mathbf{S}}_n + \mathbf{P}_n,
\]
where
\[
\tilde{\mathbf{S}}_n = \frac{1}{n} \sum_{i=1}^{\tau} \sum_{j=1}^{n_i} (\mathbf{s}_{ij} - \bar{\mathbf{s}})(\mathbf{s}_{ij} - \bar{\mathbf{s}})^{\top}\quad\text{with}\quad 
\bar{\mathbf{s}} = \frac{1}{n} \sum_{i=1}^{\tau} \sum_{j=1}^{n_i} \mathbf{s}_{ij},~\E(\mathbf{s}_{ij})=\mathbf{0},
\]
and
\[
\mathbf{P}_n = \sum_{i < j}^\tau \frac{n_i n_j}{n^2} \left\{ (\boldsymbol{\mu}_i - \boldsymbol{\mu}_j)(\boldsymbol{\mu}_i - \boldsymbol{\mu}_j)^{\top} + (\bar{\mathbf{s}}_i - \bar{\mathbf{s}}_j)(\boldsymbol{\mu}_i - \boldsymbol{\mu}_j)^{\top} + (\boldsymbol{\mu}_i - \boldsymbol{\mu}_j)(\bar{\mathbf{s}}_i - \bar{\mathbf{s}}_j)^{\top} \right\}
\]
with \(\bar{\mathbf{s}}_i =  \sum_{j=1}^{n_i} \mathbf{s}_{ij}/n_i\) for \(i = 1, \ldots, \tau\). 
Some algebra can show that the expectation of the two parts are 
\begin{align}\label{sigma-mu}
\E(\tilde{\mathbf{S}}_n) = \frac{n - 1}{n}\boldsymbol{\Sigma}_0 \quad \text{and} \quad \E(\mathbf{P}_n) \triangleq \boldsymbol{\Sigma}_{\boldsymbol{\mu}} = \sum_{1 \leq i < j \leq \tau} \frac{n_i n_j}{n^2} (\boldsymbol{\mu}_i - \boldsymbol{\mu}_j)(\boldsymbol{\mu}_i - \boldsymbol{\mu}_j)^{\top}.
\end{align}
This implies that \(\tilde{\mathbf{S}}_n\) contains solely information pertaining to \(\boldsymbol{\Sigma}_0\), while \(\mathbf{P}_n\) captures the differences in mean vectors between subgroups.
In the context of mixture models,  \(\mathbf{P}_n\) is usually unobservable. However, it can convey information through several of the largest eigenvalues of \(\mathbf{S}_n\), referred to as \emph{spiked eigenvalues}. Since \(\operatorname{rank}(\mathbf{P}_n) \leq \tau - 1\), the number of spiked eigenvalues is at most \(\tau - 1\). An example is shown in Figure \ref{intro-Sn}. Therefore, studying these spiked eigenvalues is crucial for understanding \(\boldsymbol{\Sigma}_{\boldsymbol{\mu}}\) and the underlying data structure.
Recently, the first-order convergence of these spiked eigenvalues has been investigated in \cite{yiminghigh,liu2023asymptotic}. For the second-order convergence, \cite{GMM} derived a CLT for the spikes under Gaussian assumptions. \cite{lin2024asymptotic} investigated the asymptotic distribution of spikes in the  signal-plus-noise model. All these studies are conducted under the MP regime \eqref{mp}, i.e., \(p / n \to c \in (0, \infty)\).

However, in the ultrahigh dimensional case where \(p\gg n\), the eigenvalues of \(\mathbf{S}_n\) exhibit behaviors markedly different from those in the MP regime. Properties of the spiked eigenvalues of \(\mathbf{S}_n\) induced by  \(\mathbf{P}_n\) when $p\gg n$ remain largely unknown in current literature. To fill this gap, we consider a new regime where  \(p/n \to \infty\) as \(n \to \infty\). In this scenario, unlike the MP regime, most eigenvalues of the matrix \(\mathbf{S}_n\) are zero, and all non-zero eigenvalues diverge to infinity. To address this, we renormalize the sample covariance matrix as follows:
\[
\mathbf{A}_n = \sqrt{\frac{p}{n b_p}} \left[\frac{1}{p} \mathbf{\Phi} \mathbf{X}^{\top} \mathbf{X} \mathbf{\Phi} - a_p \mathbf{\Phi}\right], \quad a_p = \frac{1}{p} \operatorname{tr}(\boldsymbol{\Sigma}_0), \quad b_p = \frac{1}{p} \operatorname{tr}(\boldsymbol{\Sigma}_0^2).
\]
$\mathbf{A}_n$ is \(n \times n\) and has \(n-1\) non-zero eigenvalues, which connect to the non-zero eigenvalues of \(\mathbf{S}_n\) through the following identity:
\[
\lambda^{\mathbf{A}_n} = \sqrt{\frac{n}{p b_p}} \lambda^{\mathbf{S}_n} - \sqrt{\frac{p}{n b_p}} a_p.
\]
Most existing studies on \(\mathbf{A}_n\) have been conducted in the context of single population case (\(\tau = 1\)). \cite{baiyin,wanglili} demonstrated that the empirical spectral distribution of \(\mathbf{A}_n\) converges weakly to the standard semicircle law, which differs from the MP law for \(\mathbf{S}_n\) in the MP regime. Similar discrepancies were reported for the largest eigenvalues of \(\mathbf{A}_n\) in \cite{bbchen}  and CLT for LSS in \cite{cbbclt,MR3416046,qiujiaxin}. While in this paper, we primarily focus on the multi-population case ($\tau>1$) when $p\gg n$.

Specifically, this paper investigates the spiked eigenvalues of the renormalized sample covariance matrix $\mathbf{A}_n$ from multi-population scenarios under the ultrahigh dimensional contexts, where 
\[
n \rightarrow \infty, \ p = p_n \rightarrow \infty, \ p / n \rightarrow \infty.
\]
Firstly, we demonstrate a phase transition phenomenon for the spiked eigenvalues of \(\mathbf{A}_n\). A critical condition for this transition is derived, showing that if the smallest eigenvalue of the matrix \(\boldsymbol{\Sigma}_{\boldsymbol{\mu}}\) exceeds a certain threshold, the first \(\tau-1\) eigenvalues of \(\mathbf{A}_n\) will fall outside the support of its limiting spectral distribution. These outliers  are referred to as \emph{distant spiked eigenvalues}. Secondly, we establish a CLT for these distant spiked eigenvalues.  Our theoretical findings are further applied to two scenarios: one is to determine the total number of subgroups, where our estimator exhibits superior numerical performance compared to existing estimators. The other is to establish a new criterion for assessing clustering outcomes when true labels are unknown.  Our new criterion integrates traditional metrics such as accuracy, recall and precision, providing highly informative guidance for tasks in unsupervised learning. Last but not least, in order to encompass most existing results derived under the MP regime, we propose a unified framework wherein
$$
n \rightarrow \infty, p=p_n \rightarrow \infty, c_n=p / n \rightarrow c \in(0, \infty].
$$
This framework incorporates all the asymptotic results across both high and ultrahigh dimensional settings, thereby broadening the applicative landscape of existing findings.

From a technical point of view, our approach to constructing the CLT of spiked eigenvalues relies on three types of random sesquilinear forms in the multi-population setting, i.e.,
\begin{align}\label{ses-infty}
\left\{\bar{\mathbf{s}}_i^{\top}\left(\mathbf{B}_n-\tilde{z}_n \mathbf{I}\right)^{-1} \bar{\mathbf{s}}_j, ~\boldsymbol{\mu}_i^{\top}\left(\mathbf{B}_n-\tilde{z}_n \mathbf{I}\right)^{-1} \boldsymbol{\mu}_j, ~\bar{\mathbf{s}}_i^{\top}\left(\mathbf{B}_n-\tilde{z}_n \mathbf{I}\right)^{-1} \boldsymbol{\mu}_j,
~ 1\leq i, j\leq\tau\right\},
\end{align}
where \(\mathbf{B}_n\) denotes the group-wise centered sample covariance matrix,
$$
\bB_{n}=\frac{1}{n} \sum_{i=1}^\tau \sum_{j=1}^{n_i}\left(\mathbf{s}_{i j}-\bar{\mathbf{s}}_i\right)\left(\mathbf{s}_{i j}-\bar{\mathbf{s}}_i\right)^{\top},\quad \bar{\mathbf{s}}_i=\frac1{n_i}\sum_{j=1}^{n_i}\mathbf{s}_{i j},
$$
and \(\tilde{z}_n\) denotes a complex number lying away from the eigenvalues of \(\mathbf{B}_n\).
Existing theories regarding \eqref{ses-infty} are primarily derived under the MP regime with $\tau\leq 2$. For $\tau=1$, the CLT for \(\bar{\mathbf{s}}_1^\top (\mathbf{B}_n - \tilde{z}_n \mathbf{I})^{-1} \bar{\mathbf{s}}_1\) was established in \cite{pan-T, 2011A}. The asymptotic properties of \(\boldsymbol{\mu}_1^\top (\mathbf{B}_n - \tilde{z}_n \mathbf{I})^{-1} \boldsymbol{\mu}_1\) were studied in \cite{bai2007asymptotics, PZ08, MR3189088}. For $\tau=2$, \cite{lihaoran} established a joint CLT for \(\{\bar{\mathbf{s}}_i^\top (\mathbf{B}_n - \tilde{z}_n \mathbf{I})^{-1} \bar{\mathbf{s}}_j : i, j = 1, 2\}\) under sub-Gaussian assumptions.  We have extended these results to the ultrahigh dimensional context and established a unified joint CLT for these quantities in \eqref{ses-infty} when $\tau\geq 1$ and $p/n\rightarrow c\in(0,\infty]$. This is by far the most general result for sesquilinear forms and is valuable in its own right. Potential applications of this joint CLT  include discriminant analysis, multivariate analysis of variance, and canonical correlation analysis, among others.

The rest of the paper is organized as follows. Section~\ref{main-cinf} details our main results, including phase transition,  CLT for distant spiked eigenvalues and the random sesquilinear forms. Section~\ref{sec:4} discusses the two applications of our findings. Section~\ref{sec:2} provides a unified framework with $p / n \rightarrow c \in(0, \infty]$.
 Section~\ref{sec:3} presents examples and simulations. 
Technical proofs are outlined in Section~\ref{sec:5} and detailed in the Supplementary Material.

\section{Main results}\label{main-cinf}
 \subsection{Preliminary}
\label{sec:concept}

For a \( p \times p \) real symmetric matrix \(\mathbf{M}_p\) with eigenvalues \(\lambda_1\geq\lambda_2\geq\cdots\geq\lambda_p\), its \emph{empirical spectral distribution} (ESD) is defined as the following probability measure:
\[ 
F^{\mathbf{M}_p} = \frac{1}{p} \sum_{j=1}^p \delta_{\lambda_j}, 
\] 
where \(\delta_{\lambda_j}\) denotes the Dirac measure at \(\lambda_j\). The limit of the  ESD sequence \(\{F^{\mathbf{M}_p}\}\)  as \(p \to \infty\), if exists, is called  \emph{limiting spectral distribution} (LSD).
For any measure \(G\) supported on the real line, the Stieltjes transform of \(G\) is defined as
\[ 
s_G(z) = \int \frac{1}{x - z} \, dG(x), \quad z \in \mathbb{C}^+,
\] 
where \(\mathbb{C}^+ = \{z \in \mathbb{C} : \Im(z) > 0\}\) denotes the upper complex plane. 
 \subsection{First-order convergence of the eigenvalues of 
	$\A_n$}
%

In this section, we present results on the first-order convergence of eigenvalues of \(\mathbf{A}_n\) in the multi-population setting when \(\tau \geq 2\) and $p\gg n$. We begin by positing some assumptions to characterize the LSD of \(\mathbf{A}_n\).

\begin{assumption}\label{as-moment-cinf}
$\mathbf{x}_{ij} = \boldsymbol{\mu}_i + \boldsymbol{\Sigma}_0^{\frac{1}{2}} \mathbf{z}_{ij}$, $\Z=(\z_{11},\ldots, \z_{1n_1},\ldots, \z_{\tau 1},\ldots,\z_{\tau n_\tau})= \left(z_{ijq}\right)_{p\times n}$, where  ${\{z_{ijq}\}}$ i.i.d. satisfy
	$$
	\E\left(z_{ijq}\right)=0, \quad 
	\E\left(z_{ijq}^{2}\right)=1, \quad 
	\E\left(z_{ijq}^{3}\right)=v_3, \quad 
	\E\left(z_{ijq}^{4}\right)=v_4<\infty.
	$$
\end{assumption}

\begin{assumption}\label{as-mp-cinf}
	The dimension $p$ and the subgroup sample sizes $\{n_1,\ldots, n_\tau\}$ are functions of the total sample size $n$ and all tend to infinity such that 
	$$
	c_n=\frac{p}{n}\to 
	\infty,\quad
	p\asymp n^t, t> 1,\quad
	k_{ni}=\frac{n_i}{n}\to k_i\in(0,1),\ i\in[1:\tau].
	$$
		Here $[a:b]$ represents the set of all integers between $a$ and $b$.
\end{assumption}

\begin{assumption}\label{as-sigma0-cinf}  
	The ESD $ H_p$ of $\bSigma_{0}$  weakly converges to a probability distribution $ H $, supported on a compact set 	$\Gamma_H\subset\mathbb R^+$.
\end{assumption}


\begin{lemma}\label{th-limit1-cinf}
	Suppose that Assumptions \ref{as-moment-cinf}-\ref{as-sigma0-cinf} hold, then almost surely, the ESD \(F^{\mathbf{A}_n}\) of \(\mathbf{A}_n\) converges weakly to the standard semicircle law with density function
	\[
f(x) = \frac{1}{2\pi} \sqrt{4 - x^2}, \quad |x| \leq 2.
	\]
\end{lemma}


As mentioned in the introduction \(\mathbf{S}_n = \tilde{\mathbf{S}}_n + \mathbf{P}_n\), where \(\mathbf{S}_n\) can be viewed as a finite rank perturbation of \(\tilde{\mathbf{S}}_n\). If the eigenvalues of \(\mathbf{P}_n\) are all significantly larger than the spectral norm of \(\tilde{\mathbf{S}}_n\), then the \(\tau - 1\) largest eigenvalues of \(\mathbf{S}_n\) will be clearly separated from its remaining eigenvalues, named as spikes.  At the population level, \(\E (\mathbf{S}_n) = (1 - 1/n)\boldsymbol{\Sigma}_0 + \boldsymbol{\Sigma}_{\boldsymbol{\mu}}\), with \(\boldsymbol{\Sigma}_{\boldsymbol{\mu}}=\E (\mathbf{P}_n)\) as defined in \eqref{sigma-mu}. The relative distance between the largest eigenvalues of \(\E (\mathbf{S}_n)\) induced by \(\boldsymbol{\Sigma}_{\boldsymbol{\mu}}\) and the spectrum of \(\boldsymbol{\Sigma}_0\) determines the total number of spikes in the spectrum of \(\mathbf{S}_n\). Similarly, for the renormalized sample covariance matrix \(\mathbf{A}_n\), we can consider a renormalized version of \(\E (\mathbf{S}_n)\), given by
\[
\boldsymbol{\Sigma}_{\mathbf{x}} = \frac{1}{\sqrt{c_n b_p}}\left(\boldsymbol{\Sigma}_0 + \boldsymbol{\Sigma}_{\boldsymbol{\mu}}\right).
\]
Some assumptions related to the spectrum of \(\boldsymbol{\Sigma}_0\) and the largest \(\tau - 1\) eigenvalues of \(\boldsymbol{\Sigma}_{\mathbf{x}}\) are listed below.


\begin{assumption}\label{as-eig0-cinf} 
	The largest eigenvalue \(\lambda_1^{\boldsymbol{\Sigma}_0}\) of \(\boldsymbol{\Sigma}_0\) satisfies
	\(
	dist\left(\lambda_1^{\boldsymbol{\Sigma}_0}, \Gamma_H\right) \rightarrow 0.
	\)
\end{assumption}

\begin{assumption}\label{as-mu-cinf} 
	The $\tau$ population means \(\{\boldsymbol{\mu}_i\}\) satisfy
	\(
	(\boldsymbol{\mu}_i - \boldsymbol{\mu}_j)^{\top} (\boldsymbol{\mu}_i - \boldsymbol{\mu}_j) \asymp \sqrt{c_n}$, ~$\forall~ i \neq j \in [1:\tau].
	\)
Moreover, the eigenvalues of $\boldsymbol{\Sigma}_{\boldsymbol{\mu}}$ in \eqref{sigma-mu} satisfy $\lambda_{j}^{\boldsymbol{\Sigma}_{\boldsymbol{\mu}}}\asymp \sqrt{c_n},~\forall j\in[1:\tau-1]$.
\end{assumption}

\begin{assumption}\label{as-sigman-cinf}
	The \(\tau - 1\) largest eigenvalues of \(\boldsymbol{\Sigma}_{\mathbf{x}}\) form \(M\) clusters, denoted as
	\[
	\{\alpha_{n,k\ell} : ~1\leq \ell \leq m_k\}, \quad k\in[1:M],
	\]
	where \(\{m_k\}\) are constants with \(\sum_{k=1}^{M} m_k = \tau - 1\). For each cluster, the eigenvalues have a common limit, i.e.,
	\[
	\alpha_{n,k\ell} = \alpha_{nk} + o\left(n^{-\frac{1}{2}}\right) \to \alpha_k, \mbox{ as } n \to \infty, ~\forall~ 1\leq \ell \leq m_k.
	\]
	And these limits are pairwise distinct, i.e., \(\alpha_i \neq \alpha_j\) for \(i \neq j \in [1:M]\).
\end{assumption}

\begin{remark}
 Assumption \ref{as-eig0-cinf}-\ref{as-mu-cinf} guarantee that all the $\tau-1$ largest eigenvalues of \(\boldsymbol{\Sigma}_{\mathbf{x}}\) are well separated from the spectrum of \(\boldsymbol{\Sigma}_0\) and all the $\tau-1$ spiked eigenvalues of \(\mathbf{A}_n\) originate from the perturbation matrix \(\boldsymbol{\Sigma}_{\boldsymbol{\mu}}\).
	
\end{remark}

\begin{remark}
	Assumption \ref{as-sigman-cinf} is general, allowing the $\tau-1$ largest eigenvalues of \(\boldsymbol{\Sigma}_{\mathbf{x}}\) to form \(M\) clusters, with the eigenvalues in each cluster sharing a common limit. Consequently, the largest \(\tau - 1\) eigenvalues of \(\mathbf{A}_n\) can be grouped into \(M\) clusters in the same manner and order as those of \(\boldsymbol{\Sigma}_{\mathbf{x}}\), denoted as
	\begin{align}\label{group}
	\left\{ \lambda_{k\ell}^{\A_{n}}: ~1\leq \ell \leq m_k\right\},\quad k\in[1:M].
	\end{align}
\end{remark}

\renewcommand{\theassumption}{1*}
\begin{assumption}\label{as-2star}
$\mathbf{x}_{ij} = \boldsymbol{\mu}_i + \boldsymbol{\Sigma}_0^{\frac{1}{2}} \mathbf{z}_{ij}$, $\Z=(\z_{11},\ldots, \z_{1n_1},\ldots, \z_{\tau 1},\ldots,\z_{\tau n_\tau})= \left(z_{ijq}\right)_{p\times n}$, where ${\{z_{ijq}\}}$ i.i.d. satisfy $E\left|z_{i j q}\right|^\kappa<\mathfrak{M}_\kappa<\infty$ for any integer $\kappa \in \mathbb{N}$. 
\end{assumption}

\setcounter{assumption}{6}
\renewcommand{\theassumption}{\arabic{assumption}}

\begin{theorem}[\emph{Phase transition}]
	\label{th-limit-cinf}
	Suppose that Assumptions \ref{as-2star},\ref{as-mp-cinf}-\ref{as-sigman-cinf} hold. The largest $\tau-1$ eigenvalues of $\A_n$ converge almost surely.
	For $k \in[1:M]$,
	\begin{itemize}
		\item [(i)] if $\alpha_k>1$, then $\lambda_{k\ell}^{\A_{n}}\xrightarrow{a.s.} \alpha_k+1/\alpha_k>2,~ \forall~1\leq \ell \leq m_k$ (distant spike);
		\item [(ii)] if $\alpha_k\leq1$, then
		$\lambda_{k\ell}^{\A_{n}}\xrightarrow{a.s.} 2, ~ \forall~1\leq \ell \leq m_k$  (close spike).
	\end{itemize}
\end{theorem}

\begin{figure}[htbp]
	\begin{subfigure}{0.495\linewidth}
		\centering
		\includegraphics[width=1\linewidth]{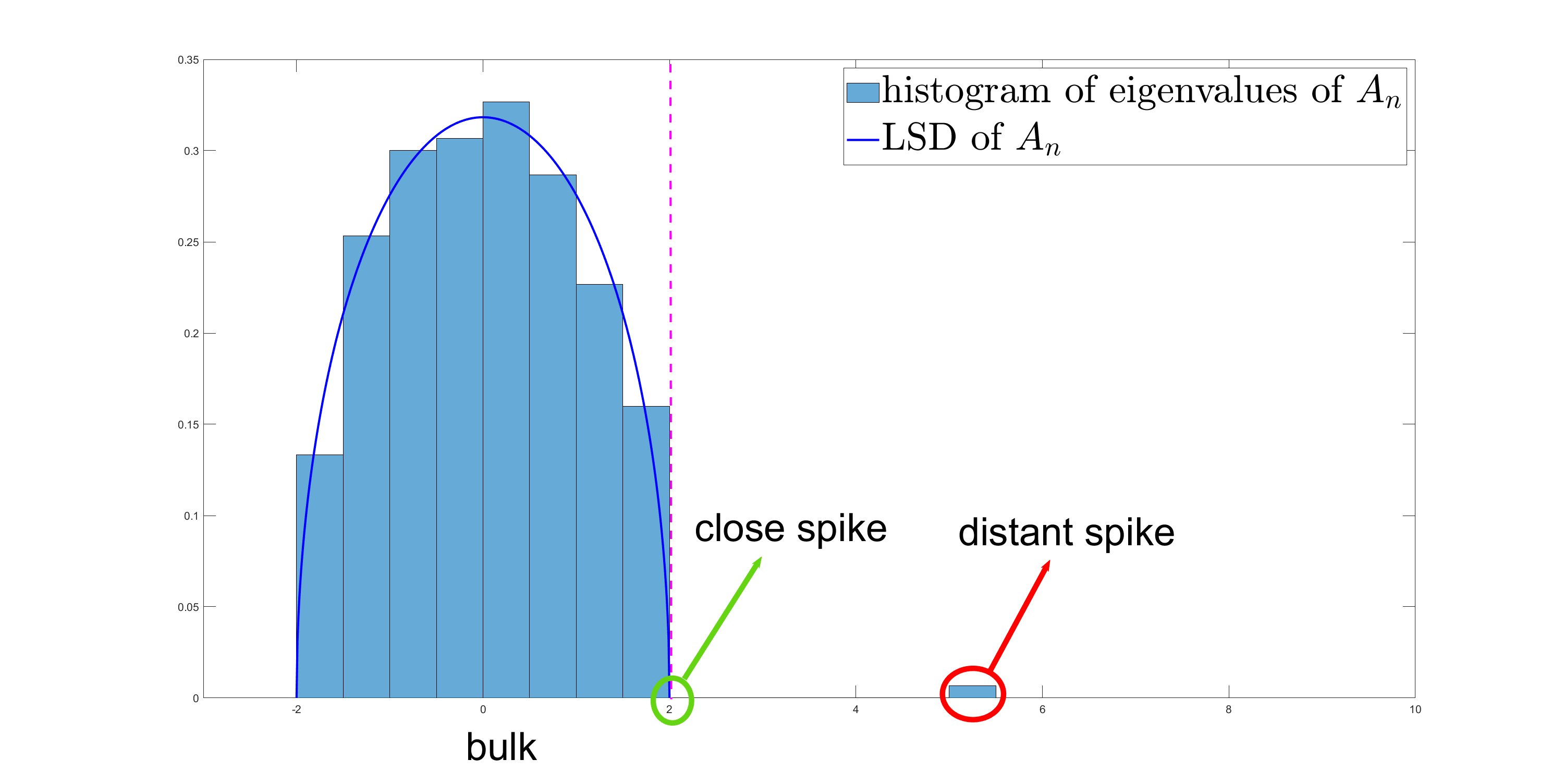}
		\caption{One distant spike.}
	\end{subfigure}
	\begin{subfigure}{0.495\linewidth}
		\centering
		\includegraphics[width=1\linewidth]{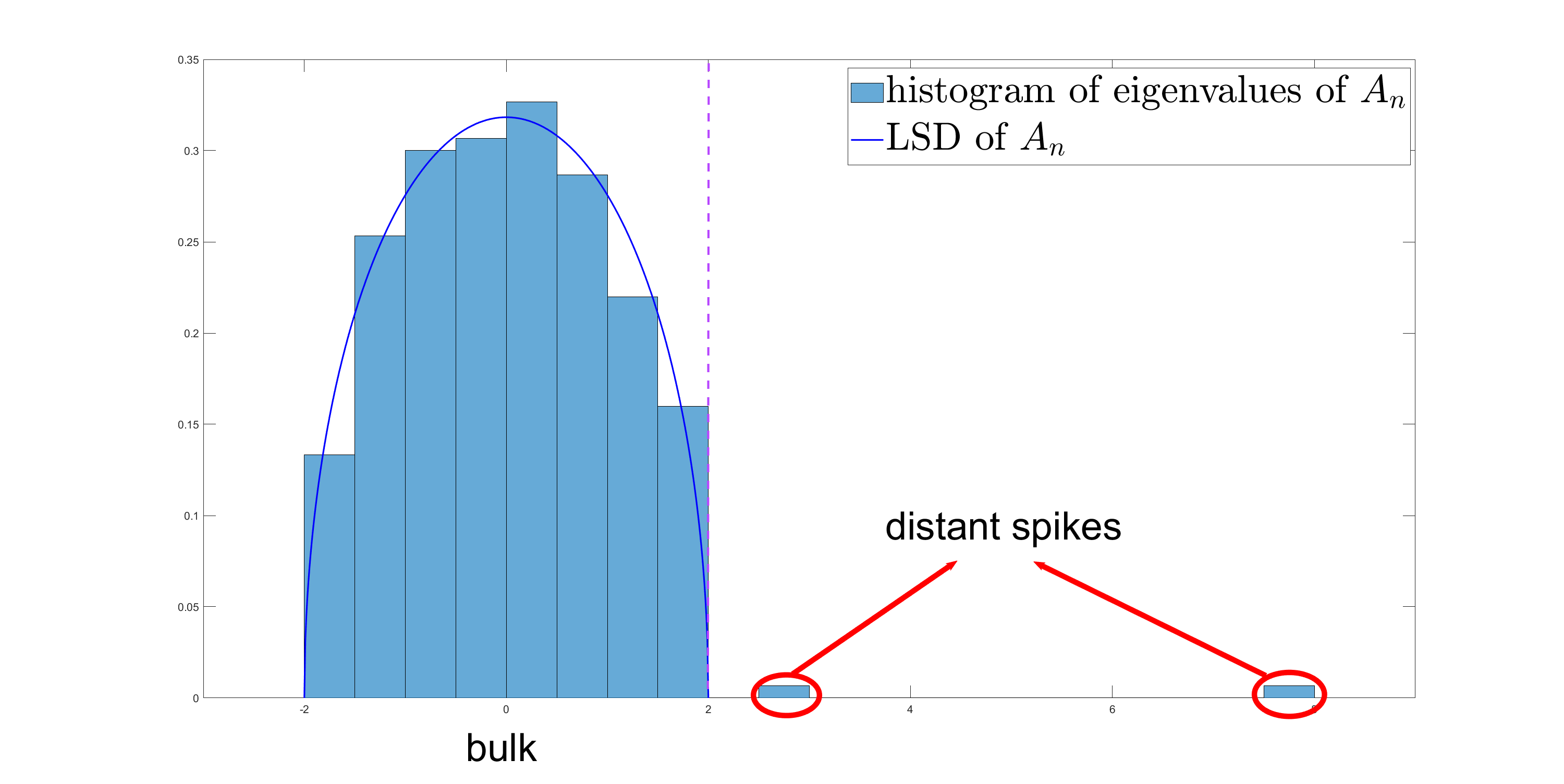}
		\caption{Two distant spikes.}
	\end{subfigure}
	\caption{The histogram and LSD of eigenvalues of  $\A_n$, where $\tau=3$, $\bSigma_0=\bI_p$, $n_1=n_2=n_3=100$, $p=300^2$, $\bmu_1=\mathbf{0}$, $\bmu_2=(-20,0,\ldots,0)$; $\bmu_3=(0,4,0,\ldots,0)$ in case (a) and $\bmu_3=(0,20,0,\ldots,0)$ in case (b). The purple dashed vertical line represents the right boundary of the standard semicircle law.
	}
	\label{intro-An}
\end{figure}

Theorem \ref{th-limit-cinf} illustrates a phase transition phenomenon of the spiked eigenvalues \(\{\lambda_{k\ell}^{\mathbf{A}_n},~1\leq \ell \leq m_k\}\) of \(\mathbf{A}_n\). The critical condition for this transition is provided, 
indicating that if \(\alpha_k > 1\), the eigenvalues \(\{\lambda_{k\ell}^{\mathbf{A}_n},~1\leq \ell \leq m_k\}\) will converge to a limit larger than 2. These spikes are referred to as \emph{distant spikes}. Otherwise, \(\{\lambda_{k\ell}^{\mathbf{A}_n},~1\leq \ell \leq m_k\}\) will converge to 2, the right edge point of the support of the semicircle law in Lemma~\ref{th-limit1-cinf}, and are called \emph{close spikes}. An example is shown in Figure \ref{intro-An}.

\subsection{CLT for distant spiked eigenvalues of $\A_n$}\label{sec:CLT-spike}
This section investigates the fluctuation of distant spiked eigenvalues. Let 
$$
\phi_{c_n, H_p}(x)=x+\frac{1}{b_p} \int \frac{t^2}{x-t / \sqrt{c_n b_p}} d H_p(t).
$$
For the \(k\)th cluster of distant spikes, i.e., \(\{\lambda_{k\ell}^{\mathbf{A}_n} : ~1\leq \ell\leq m_k\}\), we normalize it as follows:
\[
\mathbf{\Lambda}_{nk} \triangleq \sqrt{n}\left\{\lambda_{k\ell}^{\mathbf{A}_n} - \lambda_{nk} : ~1\leq \ell\leq m_k\right\}, \quad \text{with} \quad \lambda_{nk} = \phi_{c_n,H_p}(\alpha_{nk}).
\]
Its limiting distribution involves two auxiliary matrices denoted as \(\mathbf{U}_n\) and \(\mathbf{N}_n\), coming from the decomposition 
\[
\boldsymbol{\Sigma}_{\boldsymbol{\mu}} = \mathbf{U}_n \mathbf{N}_n \mathbf{U}_n^{\top},
\]
where
\begin{align}\label{con:D}
\bU_n =\left(\sqrt{k_{n1}}\boldsymbol{\mu}_{1}, \ldots, \sqrt{k_{n\tau}}\boldsymbol{\mu}_{\tau}\right)_{p\times \tau},~
\bN_n=\left(\begin{array}{cccc}
1-k_{n1} & -\sqrt{k_{n1} k_{n2}} & \ldots & -\sqrt{k_{n1} k_{n\tau}} \\
-\sqrt{k_{n1} k_{n2}} & 1-k_{n2} & \ldots & -\sqrt{k_{n2} k_{n\tau}} \\
\vdots & \vdots & \vdots & \vdots \\
-\sqrt{k_{n1} k_{n\tau}} & -\sqrt{k_{n2} k_{n\tau}} & \cdots & 1-k_{n\tau}
\end{array}\right)_{\tau\times \tau},
\end{align} 
 $k_{ni}= n_i/n$ for $i\in[1:\tau]$ and   $\bN_n$ is a projection matrix of rank $\tau-1$.

\begin{theorem}\label{th-c-cinf}
	Under Assumptions \ref{as-moment-cinf}-\ref{as-sigman-cinf},  the $m_{k}$-dimensional random vector
	${\bf \Lambda}_{nk}$
	converges in distribution to the joint distribution of the $m_{k}$ eigenvalues of the following Gaussian random matrix 
	$$
	-\sqrt{\alpha_k^{2}-1} \mathbf{Q}_k^{\top} \mathbf{N} \mathbf{W} \mathbf{N} \mathbf{Q}_k,
	$$
	where $\bN$ is the limit of $\bN_n$ defined in \eqref{con:D}, $\bQ_k$ is a $\tau\times m_k$ matrix such that 
	\begin{align}\label{mq-cinf}
	\bQ_k^{\top}\bQ_k={\bf I}_{m_k}\quad \text{and}\quad \bQ_k^{\top}\bN\bV(\alpha_k) \bN\bQ_k=-\I_{m_k},
	\end{align}
	with $\mathbf{V}(\alpha_k)=-\alpha_k^{-1} \lim _{n \rightarrow \infty} \mathbf{U}_n^{\top} \mathbf{U}_n/\sqrt{c_n b_p}$.
	And $\W=(W_{i j})$ is a symmetric $ \tau\times \tau $ random matrix with independent zero-mean Gaussian entries satisfying
		\begin{align*}
		W_{ii}\sim N(0,2\alpha_k^{-2}),\quad
		W_{ij}\sim N(0,\alpha_k^{-2}), \quad
		1\leq i< j\leq \tau.
		\end{align*}
\end{theorem}

\begin{remark}
	Theorem \ref{th-c-cinf} establishes the asymptotic distribution of the sample distant spikes \(\{\lambda_{k\ell}^{\mathbf{A}_n} : 1\leq \ell\leq  m_k\}\). This limiting distribution is jointly determined by two deterministic matrices \(\mathbf{N}\), \(\mathbf{Q}_k\)\, and the Gaussian random matrix \(\mathbf{W}\). The matrix \(\mathbf{Q}_k\), defined in \eqref{mq-cinf}, consists of \(m_k\) eigenvectors of \(\mathbf{N} \mathbf{V}(\alpha_k) \mathbf{N}\) corresponding to its eigenvalue \(-1\) with multiplicity \(m_k\). The existence of \(\mathbf{Q}_k\)  is guaranteed by Proposition 1 in \cite{GMM}. 
\end{remark}

	Note that \(\mathbf{A}_n\) involves two parameters, \(a_p\) and \(b_p\), which are typically unknown in practice. We can replace them with consistent estimators to obtain
	\begin{align}\label{hatA}
	\hat{\mathbf{A}}_n=\sqrt{\frac{p}{n \hat{b}_p}}\left[\frac{1}{p} \boldsymbol{\Phi} \mathbf{X}^{\top} \mathbf{X} \boldsymbol{\Phi}-\hat{a}_p \boldsymbol{\Phi}\right],
	\end{align}
	where
	\begin{align}\label{hata}
	\hat{a}_p=\frac{1}{p} \operatorname{tr} \hat\bS_n,\quad
	\hat{b}_p=\frac{1}{p} \operatorname{tr} \hat\bS_n^2-\frac{1}{(n-1) p}\left(\operatorname{tr} \hat\bS_n\right)^2,\quad\text{with}\quad
	\hat\bS_n=\frac{1}{n-1} \boldsymbol{\Phi} \mathbf{X}^{\top} \mathbf{X} \boldsymbol{\Phi}.
	\end{align}


\begin{theorem}\label{le-hatA-cinf}
	 Theorems \ref{th-limit-cinf}-\ref{th-c-cinf} also hold for \(\hat{\mathbf{A}}_n\) under the same assumptions.
\end{theorem}

\subsection{Asymptotic properties for random sesquilinear forms}\label{sec-ses-cinf}

The proofs of Theorem \ref{th-limit-cinf}-\ref{th-c-cinf} rely on three types of random sesquilinear forms in the multi-population setting, i.e.,
\begin{align*}
	\left\{\bar{\mathbf{s}}_i^{\top}\left(\mathbf{B}_n-\tilde{z}_n \mathbf{I}\right)^{-1} \bar{\mathbf{s}}_j, ~\boldsymbol{\mu}_i^{\top}\left(\mathbf{B}_n-\tilde{z}_n \mathbf{I}\right)^{-1} \boldsymbol{\mu}_j, ~\bar{\mathbf{s}}_i^{\top}\left(\mathbf{B}_n-\tilde{z}_n \mathbf{I}\right)^{-1} \boldsymbol{\mu}_j,
	~ 1\leq i, j\leq\tau\right\},
\end{align*}
where $\bB_{n}=\frac{1}{n} \sum_{i=1}^\tau \sum_{j=1}^{n_i}\left(\mathbf{s}_{i j}-\bar{\mathbf{s}}_i\right)\left(\mathbf{s}_{i j}-\bar{\mathbf{s}}_i\right)^{\top},~\bar{\mathbf{s}}_i=\frac1{n_i}\sum_{j=1}^{n_i}\mathbf{s}_{i j}$.
This section establishes the first- and second-order convergence of these random sesquilinear forms where  \(\tilde{z}_n\) is specified as
\[
\tilde{z}_n = c_n a_p + \sqrt{c_n b_p} z, \quad z \in \mathbb{C}\setminus[-2,2].
\]
These results are not only fundamental for proving Theorems \ref{th-limit-cinf} and \ref{th-c-cinf}, but also hold independent  significance for other inference procedures. We begin by establishing the first-order convergence.

\begin{theorem}\label{th-as-cinf}
	Suppose that Assumptions \ref{as-moment-cinf}-\ref{as-eig0-cinf} hold and $\bmu_{i}^{\top}\bmu_i\asymp\sqrt{c_n}$ for all $1\leq i\leq \tau$, then for any $1\leq i,j\leq \tau$,  we have,
	\begin{align*}
	&\frac{\tilde{z}_n}{\sqrt{c_n b_p}} \bar{\mathbf{s}}_i^{\top}\left(\mathbf{B}_n-\tilde{z}_n \mathbf{I}\right)^{-1} \bar{\mathbf{s}}_j
	+\frac{1}{k_{i}}\left[\sqrt{c_n/b_p}a_p+z+1/s(z)\right]I(i=j)
	\xrightarrow{a . s .} 0,
	\\
	&\frac{\tilde{z}_n}{\sqrt{c_n b_p}} \bar\s_i^{\top}\left(\mathbf{B}_n-\tilde{z}_n \mathbf{I}\right)^{-1} \boldsymbol{\mu}_j	\xrightarrow{a . s .} 0,\\
	&\frac{\tilde{z}_n}{\sqrt{c_n b_p}} \bmu_i^{\top}\left(\mathbf{B}_n-\tilde{z}_n \mathbf{I}\right)^{-1} \boldsymbol{\mu}_j +1/\sqrt{c_nb_p}\bmu_i^{\top}\bmu_j\xrightarrow{a . s .} 0.
	\end{align*}
	where $I(\cdot)$ denotes the indicator function.
\end{theorem}


To illustrate the second-order convergence, we denote 
 $ s_0=s_0(z) $  as the solution to the equation 
\begin{equation*}
z=-\frac{1}{s_0}-\frac{s_0}{b_p} \int \frac{t^2}{1+s_0 t / \sqrt{c_n b_p}} d H_p(t),\quad z\in \mathbb C^+,
\end{equation*} 
which is a finite-sample proxy for $ s(z) $. Denote
\begin{align*} 
\bM_{\bar\s}=\left(\bar\s_1,\ldots,\bar\s_\tau\right)_{p\times\tau},\quad
\bM_{\bmu} =\left(\bmu_1,\ldots,\bmu_\tau\right)_{p\times\tau},\quad
\bK_n={\rm diag}\left(
k_{n1},\ldots, k_{n\tau}\right),
\end{align*}
then we rearrange these sesquilinear forms into a 
$ 2\tau\times 2\tau $ matrix and normalize it as follows:
\begin{align}\label{eq:Ln}\nonumber
\mathbf{L}_n \triangleq &
\sqrt{n}\left\{ 
\frac{\tilde{z}_n}{\sqrt{c_n b_p}}
\left(\mathbf{M}_{\bar{\mathbf{s}}}, \mathbf{M}_{\boldsymbol{\mu}}\right)^{\top}\left(\mathbf{B}_n-	\tilde{z}_n \mathbf{I}\right)^{-1}\left(\mathbf{M}_{\bar{\mathbf{s}}}, \mathbf{M}_{\boldsymbol{\mu}}\right) \right.\\
& \left.-\left(\begin{array}{cc}
-\left(\sqrt{\frac{c_n}{b_p}}a_p+z+\frac{1}{ s_0(z)}\right) \mathbf{K}_n^{-1} & \mathbf{0}_{\tau \times \tau} \\
\mathbf{0}_{\tau \times \tau} & -\mathbf{M}_{\boldsymbol{\mu}}^{\top}\left(\sqrt{c_nb_p}\mathbf{I}+s_0(z) \boldsymbol{\Sigma}_0\right)^{-1} \mathbf{M}_{\boldsymbol{\mu}}
\end{array}\right)\right\}_{2\tau\times 2\tau}.
\end{align}
The asymptotic distribution of $\mathbf{L}_n$ is established in the following theorem.	
\begin{theorem}\label{th-clt-cinf}
Suppose that Assumptions \ref{as-moment-cinf}-\ref{as-eig0-cinf} hold and $\bmu_{i}^{\top}\bmu_i\asymp\sqrt{c_n}$ for all $1\leq i\leq \tau$, then \(\mathbf{L}_n\) converges in distribution to a symmetric \(2\tau \times 2\tau\) random matrix \(\mathbf{L} = (L_{ij})\) of Gaussian entries with zero-mean and covariance
	\begin{align*}
	\Cov(L_{ij},L_{lt})=\begin{cases}
	\frac{2}{1-s^2(z)}\frac{1}{k_i^2}&\mbox{if $ i=j=l=t\in[1:\tau] $},\\
	\frac{1}{1-s^2(z)}\frac{1}{k_ik_j}&\mbox{if $ i=l\neq j=t\in[1:\tau] $},\\
	0&\mbox{o.w..}
	\end{cases}
	\end{align*}

\end{theorem}

	\section{Applications}\label{sec:4}
	
This section explores two scenarios in which our theoretical findings can be applied.
\begin{itemize}
    \item [1.] {\bf Determination of number of subgroups.} 
Estimating the number of subgroups is crucial for revealing the underlying data structure. In the MP regime (\(p/n\to c \in(0,\infty)\)), \cite{liu2023asymptotic} proposed two methods, referred to as EDA and EDB, based on distances of adjacent eigenvalues. However both EDA and EDB cannot handle the case when $p\gg n$. For the ultrahigh dimensional case where \(p/n\to \infty\), \cite{MR4606323} introduced a method based on eigenvalue ratios, referred to as the ER method.
Yet the authors did not provide a theoretical analysis, and its performance remains unknown in general cases. In this section, we have developed a new method, based on the theoretical results obtained in previous sections, to accurately determine the number of subgroups in ultrahigh dimensional data. Our method leverages the properties of spiked eigenvalues generated from the pairwise distances of the population means and exhibits superior numerical performance compared to others.

\item [2.] {\bf Assessment of clustering results.} 
Accuracy, recall and precision are usually used to evaluate the clustering results. However, these criteria can only be obtained when the true labels are known. In this section, we propose a novel criterion for evaluating clustering results when the true labels are unknown.  Our criterion is designed based on the asymptotic properties of the spiked eigenvalues from multi-population and provides a robust measure for assessing the quality of clustering in unsupervised learning.
\end{itemize}

	\subsection{Determination of the number of subgroups}
	Consider the data matrix $\mathbf{X}=\left(\mathbf{x}_{ij}\right)_{p\times n}$ from multi-population
	\begin{align*}
		\mathbf{x}_{ij} = \boldsymbol{\mu}_i + \boldsymbol{\Sigma}_0^{\frac{1}{2}} \mathbf{z}_{ij} ,~ i = 1, \ldots, \tau, \ j = 1, \ldots, n_i,
	\end{align*}
	and $\sum_{i=1}^\tau n_i=n,~k_i=\frac{n_i}{n}$, $c_n=\frac{p}{n}$.  
	Utilizing the eigenvalues of $\hat\A_n$ in \eqref{hatA},
	we estimate the total number  of subgroups \(\tau\) as 
	\begin{align}\label{hat-tau}
	\hat{\tau}=\max \left\{i: \lambda_{i}^{\hat\A_n} \geq 2+d_n\right\}+1,
	\end{align}
	where $d_n$ is a positive vanishing constant.

Suppose that the $\tau$ subgroups are well separated such that
the first \(\tau - 1\) largest eigenvalues of \(\boldsymbol{\Sigma}_{\mathbf{x}}\) are all distant spikes. Consequently, there will be \(\tau - 1\) spiked eigenvalues of \(\hat{\mathbf{A}}_n\), while the rest are bulk eigenvalues, bounded by the right edge of the support of the LSD \(F\). Thus, our method is to find a critical value to distinguish the spiked eigenvalues from the bulk ones, i.e., \(\lambda_{\tau-1}^{\hat{\mathbf{A}}_n}\) and \(\lambda_{\tau}^{\hat{\mathbf{A}}_n}\). 
The consistency of our estimator naturally holds as follows.




	\begin{theorem}\label{th-num}
		Suppose that  Assumptions~\ref{as-2star}, \ref{as-mp-cinf}-\ref{as-sigman-cinf}  hold and the first $\tau-1$ large eigenvalues of $\bSigma_{\bx}$ 
		are all distant spikes, i.e., $\alpha_k>1$ for $k \in[1:M]$.
		Let $d_n\to 0$ and 
		$n^{ \kappa}d_n\to \infty$  for any $\kappa>0$ as $n\to \infty$.
		Then  $\hat\tau \xrightarrow{p}\tau$ as $n\to\infty$.
	\end{theorem}

To examine the performance of our estimate,  we compare with the EDA and EDB methods in \cite{liu2023asymptotic} and the ER method in \cite{MR4606323}. The total number of subgroups  \(\tau = 4\) and the sample sizes $n_i$ are balanced with \(k_{n1} = k_{n2} = k_{n3} = k_{n4} = 1/4\). $\boldsymbol{\Sigma}_0^{\frac{1}{2}}$ is tridiagonal, with all main diagonal elements equal to $1$ and all subdiagonal and superdiagonal elements equal to $0.5$. For the mean vectors, we consider two cases:
	\begin{itemize}
	\item []Case 1 (weak signals). 
	\begin{align*}
	\boldsymbol{\mu}_1 & =c_n^{\frac{1}{4}}(0,  0, 0,5,0\ldots, 0)^{\top}, \quad \boldsymbol{\mu}_2  =c_n^{\frac{1}{4}}(4\sqrt{1.8}, 0,0,5,0\ldots, 0)^{\top}, \\ \boldsymbol{\mu}_3 & =c_n^{\frac{1}{4}}(2\sqrt{1.8}, 2\sqrt{5.4} ,  0, 5,0\ldots, 0)^{\top}, \quad \boldsymbol{\mu}_4  =c_n^{\frac{1}{4}}(2\sqrt{1.8}, 2\sqrt{0.6}, -4 \sqrt{1.2},  5,0, \ldots, 0)^{\top}.
	\end{align*}
	\item []Case 2 (strong signals). 
	\begin{align*} \boldsymbol{\mu}_1 & =c_n^{\frac{1}{4}}(0,  0, 0,5,0\ldots, 0)^{\top}, \quad \boldsymbol{\mu}_2  =c_n^{\frac{1}{4}}(4\sqrt{2.8}, 0,0,5,0\ldots, 0)^{\top}, \\ \boldsymbol{\mu}_3 & =c_n^{\frac{1}{4}}(2\sqrt{2.8},  2\sqrt{8.4},  0, 5,0\ldots, 0)^{\top}, \quad \boldsymbol{\mu}_4  =c_n^{\frac{1}{4}}(2\sqrt{2.8}, \frac{2}{3}\sqrt{8.4}, -\frac{8}{3}\sqrt{6},  5,0, \ldots, 0)^{\top}.\end{align*}
\end{itemize}

The variables \(\z_{ij}=\{z_{ijq}~,q\in[1:p]\}\) are generated from
		\begin{itemize}
			\item [(1)] Gaussian distribution \(N(0,1)\);
			\item [(2)] Symmetric Bernoulli distribution with outcomes 1 and -1. 
		\end{itemize}
The tuning parameter \(d_n\) is set to be \(1/(\log n)^2\). The dimensional settings are $n=300$ and \(c_n =p/n=\{ 250, 500, 1000,2000,5000,10^4\}\). Table \ref{app-num} reports estimation accuracy based on 5000 replicates. It's clear that our estimator \(\hat{\tau}\) has superior performance across all settings.

			\begin{table}[t!]
			\centering
			\caption{Empirical accuracy of \(\hat{\tau}\), EDA, EDB, and ER from 5000 replications for Cases I and II.}
			\label{app-num}\par
			\vspace{1em}
			\resizebox{\linewidth}{!}{
				\begin{tabular}{lcccccccccccc}
					\hline
					&\multicolumn{6}{c} {Case I}&\multicolumn{6}{c}{Case II}
					\\
					\cmidrule(r){2-7}\cmidrule(r){8-13}
					$p/n$&250&500&1000
					&2000&5000&$10^4$
					&250&500&1000
					&2000&5000&$10^4$\\
					\hline
					\multicolumn{2}{c}{Gaussian data}
					\\
					$\hat\tau$
					&\textbf{0.9352}&\textbf{0.9866}&\textbf{0.9956}&\textbf{0.9972}&\textbf{0.9936}&\textbf{0.9948}&\textbf{0.9992}
					&\textbf{0.9998}&\textbf{1}&\textbf{1}&\textbf{1}&\textbf{1}
					\\
					EDA&0.0538&0.0630&0.0276&0.0026&0&0
					&0.8434&0.8722&0.8848&0.8368&0.6864&0.4388
					\\
					EDB&0&0&0&0&0&0
					&0.6544&0.6280&0.4932&0.2258&0&0\\
					ER&0.8370&0.8870&0.9134&0.9346&0.9452&0.9544
					&0.9958&0.9970&0.9960&0.9982&0.9996
					&0.9990\\
					\hline
					\multicolumn{2}{c}{Non-Gaussian data}
					\\
					$\hat\tau$
					&\textbf{0.9386}&\textbf{0.9894}&\textbf{0.9950}&\textbf{0.9968}&\textbf{0.9924}&\textbf{0.9950}&\textbf{0.9998}
					&\textbf{1}&\textbf{1}&\textbf{1}&\textbf{1}&\textbf{1}
					\\
					EDA&0.0580&0.0598&0.0304&0.0034&0&0
					&0.8500&0.8836&0.8900&0.8456&0.6968&0.4650
					\\
					EDB&0&0&0&0&0&0
					&0.6578&0.6408&0.4956&0.2226&0.0012&0\\
					ER&0.8540&0.8978&0.9208&0.9350&0.9402&0.9512
					&0.9962&0.9970&0.9992&0.9980&0.9988
					&0.9978\\
					\hline
			\end{tabular}}
		\end{table}

	\subsection{Assessment of clustering results}
	Consider data matrix $\X=(\bx_{ij})_{p\times n}$  from two populations, $$
	\bx_{ij}=\bmu_i+\bSigma_0^{\frac{1}{2}}\z_{ij},~ i=1,2,\ j\in[1:n_i],~n_1+n_2=n,~ n_1/n\leq 1/2.
	$$
 Denote the true labels as $\y=(y_i)_{1\times n}$ and estimated ones as $\check{\mathbf{y}}=\left(\check{y}_i\right)_{1\times n}$. Suppose 
	\begin{align*}
	\mathbf{y} = (y_i) = (\underbrace{1, \ldots, 1}_{n_1}, \underbrace{2, \ldots, 2}_{n_2}), 
	\end{align*}
then the accuracy, recall and precision of this clustering result $\check{\mathbf{y}}$ can be written as
	\[
	\text{ACC} = \frac{\sum_{i=1}^n I(y_i = \check{y}_i)}{n}, \quad \text{REC} = \frac{\sum_{i=1}^{n_1} I(\check{y}_i = 1)}{n_1},\quad 
	\text{PRE} = \frac{\sum_{i=1}^{n_1} I(\check{y}_i = 1)}{\sum_{i=1}^{n} I(\check{y}_i = 1)}.
	\]
These metrics can only be calculated when the true labels $\y$ are known. 

In this section, we propose a new criterion for evaluating clustering results $\check{\mathbf{y}}$ when $\y$ is unknown.
Specifically, denote \[
\check{n}_i=\sum_{k=1}^n I\left(\check{y}_k=i\right),~\check{\mathbf{x}}_{i}=\sum_{k=1}^n\bx_kI\left(\check{y}_k=i\right)/\check{n}_i,~i=1,2,
\] the metric we use to evaluate the discrepancy between $\check{\mathbf{y}}$ and $\y$ is given by
\[
T = \frac{\check{\alpha}_{n}}{\hat{\alpha}_{n}},
\]
where 
\begin{align*}
	\hat\alpha_{n}&=\frac{1}{2}\left(\lambda_{\max}^{\hat\A_n}+\sqrt{(\lambda_{\max}^{\hat\A_n})^2-4}\right),~\check{\alpha}_{n} = \sqrt{\frac{n}{p \hat{b}_p}} \frac{\check{n}_1 \check{n}_2}{n^2}\left(\check{\mathbf{x}}_1-\check{\mathbf{x}}_2\right)^{\top}\left(\check{\mathbf{x}}_1-\check{\mathbf{x}}_2\right) - \sqrt{\frac{p}{n \hat{b}_p}} \hat{a}_p.
\end{align*} 
Here \(\{\hat{a}_p, \hat{b}_p\}\) and \(\hat\A_n\) in \eqref{hatA} are directly obtained from $
\X$.

Our method is motivated by the following decomposition:
\[
\E (\bS_n) = \frac{n-1}{n} \boldsymbol{\Sigma}_0 + \boldsymbol{\Sigma}_{\boldsymbol{\mu}},~ \text{with}~ \bS_n=\frac{1}{n} \mathbf{X} \boldsymbol{\Phi} \mathbf{X}^{\top},\ \boldsymbol{\Sigma}_{\boldsymbol{\mu}}=\frac{n_1 n_2}{n^2}\left(\boldsymbol{\mu}_1-\boldsymbol{\mu}_2\right)\left(\boldsymbol{\mu}_1-\boldsymbol{\mu}_2\right)^{\top},
\]
which indicates that the information about the clusters encoded in \(\bS_n\) arises from the perturbation \(\boldsymbol{\Sigma}_{\boldsymbol{\mu}}\). Note that
\(\boldsymbol{\Sigma}_{\boldsymbol{\mu}}/\sqrt{c_nb_p}\) is of rank one, and its  eigenvalue is given by
\[\alpha_n=\sqrt{\frac{n}{p b_p}} \frac{n_1 n_2}{n^2}\left(\boldsymbol{\mu}_1-\boldsymbol{\mu}_2\right)^{\top}\left(\boldsymbol{\mu}_1-\boldsymbol{\mu}_2\right).\] 
Our strategy is to estimate $\alpha_n$ from two different perspectives.   
On one hand, we can directly use $\hat\alpha_{n}$, because from Theorems \ref{th-limit-cinf} and \ref{le-hatA-cinf}, we know that \(\hat\alpha_n\) is a consistent estimator of \(\alpha_n\) as long as $\alpha_n > 1$.  On the other hand, if the clustering result $\check{\mathbf{y}}$ is accurate, then \(\check\alpha_n\) can also consistently estimate \(\alpha_n\). In general, the more accurate $\check{\mathbf{y}}$ is, the closer our criterion \(T\) is to 1. 

In fact, as the theorem below suggests, \(T\) can be regarded as an approximation of a composite measure of accuracy and recall, represented as
\begin{align*}
T_0 = \frac{k_{n1}\left(1 - \text{ACC} - k_{n1} - \text{REC} + 2k_{n1} \text{REC} \right)^2}{(1 - k_{n1})(1 - \text{ACC} - k_{n1} + 2k_{n1} \text{REC})(\text{ACC} + k_{n1} - 2k_{n1} \text{REC})},
\end{align*}
where $k_{n1}=n_1/n$. 	
\begin{theorem}\label{th-T}
	For $\tau=2$,	suppose Assumptions \ref{as-moment-cinf}-\ref{as-mu-cinf} hold  and  $\lim_{n\rightarrow\infty}\alpha_n>1$,  then
		\[
		T-T_0\xrightarrow{i.p.}0.
		\]
  	\end{theorem}

The measure \(T_0\) takes values in \([0,1]\) and, in particular, \(T_0 = 1\) if and only if \(\text{ACC} = 1\). When \(\text{ACC} < 1\),  \(T_0\) will vary according to \(\text{ACC}\) and \(\text{REC}\). For a fixed \(\text{ACC}\), the maximum value of \(T_0\) is given by
\begin{align}\label{max-T}
\max_{REC} T_{0} = \frac{k_{n1}}{1-k_{n1}} \cdot \frac{\text{ACC} - k_{n1}}{1 - (\text{ACC} - k_{n1})}, \mbox{ if $\text{REC} = 1$},
\end{align}
and the minimum is given by
\begin{align}\label{min-T}
\min_{REC} T_{0} =
\begin{cases}
\displaystyle 0, & \text{if } k_{n1} < \frac{1}{2} \text{ and } k_{n1} + \text{ACC} \leq 1,\\
\displaystyle \frac{k_{n1} - \text{ACC}}{k_{n1}} \cdot \frac{1 - k_{n1} - \text{ACC}}{1 - k_{n1}}, & \text{otherwise}.
\end{cases}
\end{align}
Here the minimum of $T_0$ is achieved when
\(\text{REC} = (1 - \text{ACC} - k_{n1})/(1 - 2k_{n1})\) in the first case and \(\text{REC} = \text{ACC}(\text{ACC} + k_{n1} - 1)/\{(2\text{ACC} - 1)k_{n1}\}\) in the second case. Clearly, both the minimum and maximum values of \(T_0\) monotonically increase with \(\text{ACC}\). As a result, our metric \(T\) generally captures the trend of \(\text{ACC}\).

\begin{figure}[h!]
		\centering
		\begin{tabular}{cc}
			\begin{minipage}[t]{2.5in}
				\includegraphics[width=2.5in]{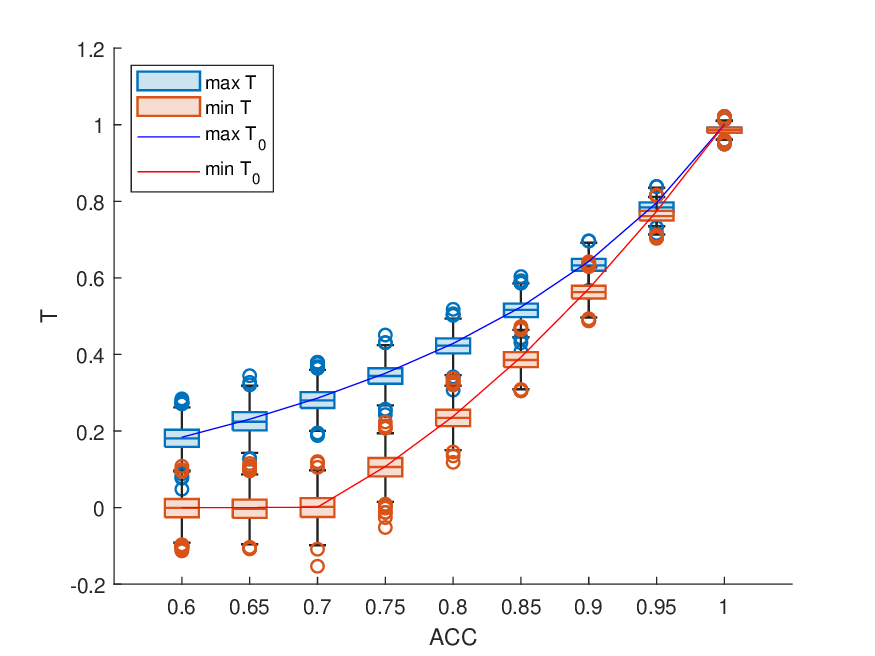}
			\end{minipage}
			\begin{minipage}[t]{2.5in}
				\includegraphics[width=2.5in]{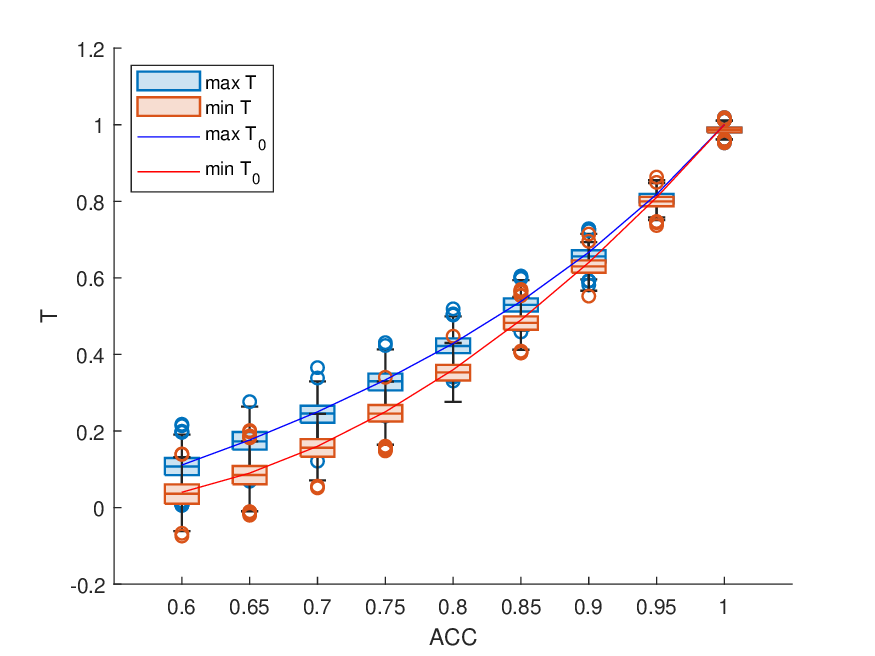}
			\end{minipage}\\
		\end{tabular}
		\caption{
Box plots of $T$ when $T_0$ achieves maximum (denoted by $\max T$) and minimum (denoted by $\min T$).  \(k_{n1}=0.3\) (left) and \(k_{n1}=0.5\) (right).  \(\text{ACC}\) ranges from $[0.6, 1]$. 
		}\label{fig-change}
	\end{figure}

A numerical illustration is presented in Figure \ref{fig-change}, where the population model is
\begin{align*}
\bSigma_0 = \operatorname{diag}(\underbrace{1, \ldots, 1}_{p / 2}, \underbrace{2, \ldots, 2}_{p / 2}),\quad
\bmu_1 = \mathbf{0},\quad
\bmu_2 = w \mathbf{1}_p,\quad (z_{ijq}) \text{ iid from } N(0,1).
\end{align*}
The dimensional setting is \(n=400\) and \(p=n^2\). The parameter \(w\) is 0.0531 for \(k_{n1}=0.3\) and \(w=0.0487\) for \(k_{n1}=0.5\) satisfying $\alpha_n=3$ in both cases. Box plots of \(T\) when $T_0$ achieves maximum and minimum,  where \(\text{REC}\) is set as \eqref{max-T} and \eqref{min-T} respectively,  are shown in Figure~\ref{fig-change} from 1000 independent replications.

\begin{figure}[h!]
    	\centering
    	\begin{tabular}{ccc}
    		\begin{minipage}[t]{1.8in}
    			\includegraphics[width=1.8in]{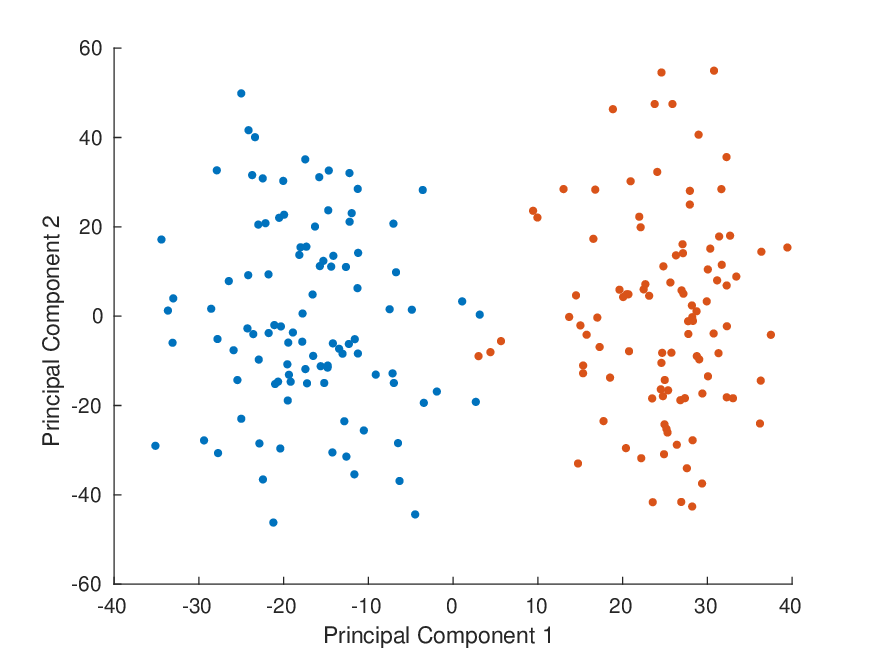}
    		\end{minipage}
    		\begin{minipage}[t]{1.8in}
    			\includegraphics[width=1.8in]{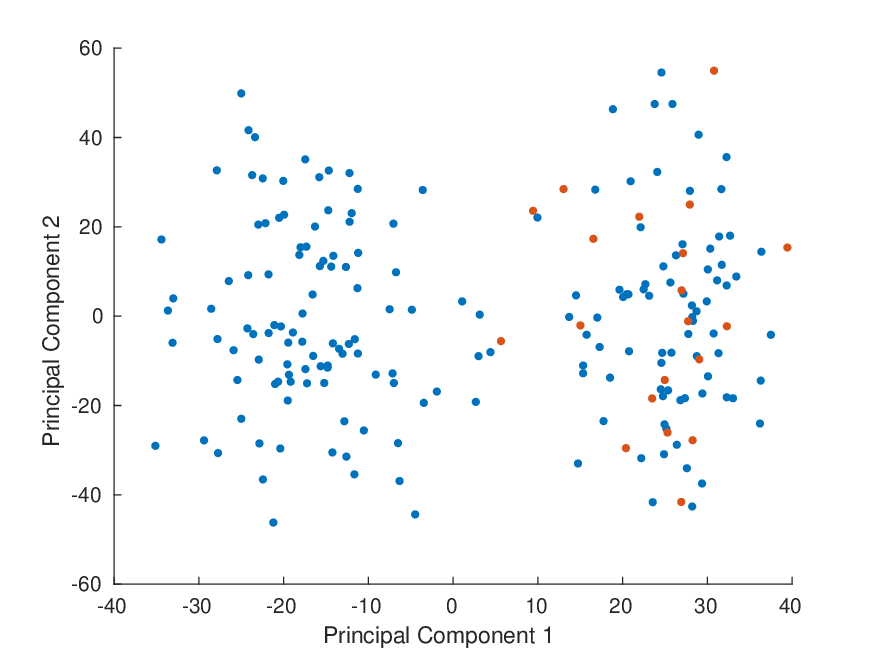}
    		\end{minipage}
    		\begin{minipage}[t]{1.8in}
    			\includegraphics[width=1.8in]{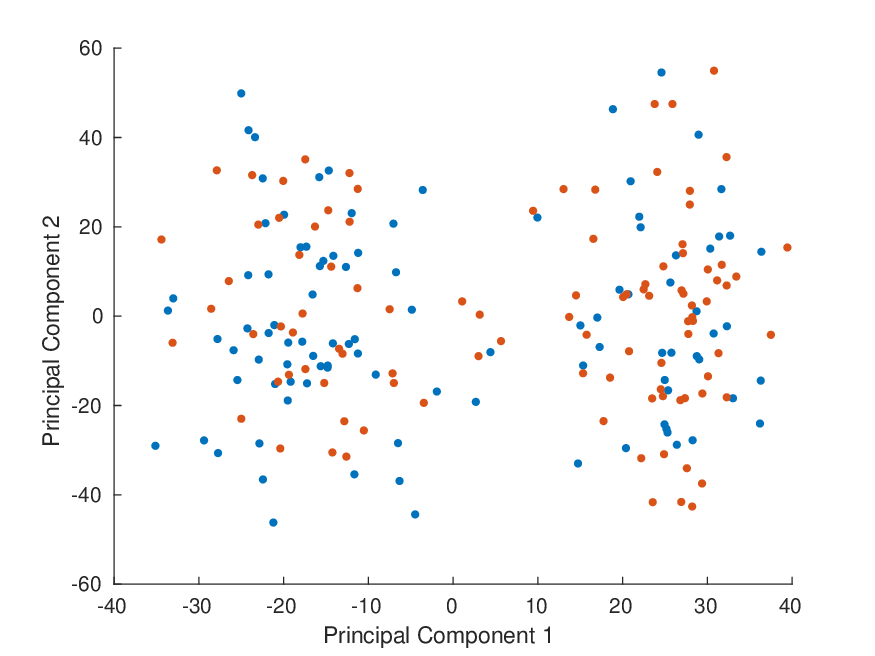}
    		\end{minipage}
    	\end{tabular}
    	\caption{
    \(n=200\) sample observations through the first two principal components, \(p=n^2\). The left panel shows the true labels. The middle and right panels display two clustering results with \(\{\text{ACC}=0.6, \text{REC}=1, \max T_0=1/9\}\) and \(\{\text{ACC}=0.6, \text{REC}=0.6, \min T_0=1/25\}\).	}\label{fig-pca}
    \end{figure}

For a fixed ACC,  \(T_0\) varies with \(\text{REC}\), reflecting the dispersion of incorrect labels. \(T_0\) achieves  maximum when \(\text{REC} = 1\), indicating that all wrong labels occur within the second cluster. Conversely, if the incorrect labels appear evenly across two clusters, then \(T_0\) reaches  minimum . This occurs, for instance, when \(\text{REC} = \text{ACC}\) with \(k_{n1} = 1/2\). We illustrate this phenomenon in Figure \ref{fig-pca}. The left panel shows the true labels of \(n = 200\) sample observations with \(n_1=100\) (blue) and \(n_2=100\) (red). The middle and right panels  exhibit two clustering results with the same accuracy level (\(\text{ACC} = 0.6\)). The middle panel represents the case of \(\max T_0\) while the right  representing the case of \(\min T_0\).

Lastly, we generalize our method to the case of $\tau>2$. Suppose that the \(\tau-1\) spiked eigenvalues of \(\hat{\mathbf{A}}_n\) are all larger than 2. Accordingly, the sample observations are clustered into \(\tau\) groups, with sample sizes \(\{\check n_i,~i\in[1:\tau]\}\) and sample means 
\(\{\check \bx_i~i\in[1:\tau]\}\). Then, we evaluate this clustering result using 
\[
T_{\tau} \triangleq \frac{\text{tr}(\hat{\boldsymbol{\Sigma}}_{\tau})}{\sum_{k=1}^{\tau-1} \hat\alpha_{nk}},
\]
where $\hat\alpha_{nk} = \frac{1}{2}\left\{\lambda_{k}^{\hat{\mathbf{A}}_n} + \sqrt{\left(\lambda_{k}^{\hat{\mathbf{A}}_n}\right)^2 - 4}\right\}$ and 
\begin{align*}
\hat{\boldsymbol{\Sigma}}_{\tau} = \sqrt{\frac{n}{p \hat b_p}} \check{\mathbf{N}}_n 
\begin{pmatrix}
\sqrt{\check k_{n 1}} \check{\mathbf{x}}_1^{\top} \\
\vdots \\
\sqrt{\check k_{n \tau}} \check{\mathbf{x}}_{\tau}^{\top}
\end{pmatrix}
\begin{pmatrix}
\sqrt{\check k_{n 1}} \check{\mathbf{x}}_1 & \cdots & \sqrt{\check k_{n \tau}} \check{\mathbf{x}}_{\tau}
\end{pmatrix}
\check{\mathbf{N}}_n - \sqrt{\frac{p}{n \hat b_p}} \hat a_p \check{\mathbf{N}}_n.
\end{align*}
Here, \(\check k_{n i} = \check n_i/n\) and the matrix \(\check {\mathbf{N}}_n\) is an analogue of \(\mathbf{N}_n\) in \eqref{con:D} with \(k_{ni}\) replaced by \(\check k_{ni}\). Similar to the two-sample case, \(T_{\tau}\) generally aligns with the \(\text{ACC}\) measure, and will tend to 1 when \(\text{ACC} = 1\) under certain mild conditions.

 	\section{A unified framework} 
 \label{sec:2} 
 

In this section, we have integrated all the theoretical results derived under the celebrated MP regime and the ultrahigh dimensional setting, thereby broadening the applicability of our theory. We focus on the asymptotic properties of spiked eigenvalues of the renormalized random matrix $\mathbf{A}_n$ from multi-population scenarios, under the general asymptotic regime where
\[
n \rightarrow \infty, \ p = p_n \rightarrow \infty, \ c_n=p / n \rightarrow c \in (0, \infty].
\]
Specifically, the unified LSD, phase transition, CLT for distant spikes and the asymptotic properties of the random sesquilinear \eqref{ses-infty} are established here.  In particular, we have eliminated the Gaussian constraints for the CLT of distant spikes in \cite{GMM}. More importantly, we are the first to establish joint CLT for several sesquilinear forms from multi-population scenarios, which holds its own value in high dimensional inference problems.

 \subsection{Unified LSD and phase transition of $\A_n$}

 \renewcommand{\theassumption}{2*}
 	\begin{assumption}\label{as-mp}
 		The dimension $p$ and the sample sizes $\{n_1,\ldots, n_\tau\}$ are functions of the total sample size $n$ and all tend to infinity such that 
 		$$
 		c_n=\frac{p}{n}\to 
 		c\in(0,\infty],\quad
 		p\asymp n^t, t\geq 1,\quad
 		k_{ni}=\frac{n_i}{n}\to k_i\in(0,1),\ i\in[1:\tau].
 		$$
 	\end{assumption}

 \setcounter{assumption}{6}
 \renewcommand{\theassumption}{\arabic{assumption}}

 \begin{lemma}\label{th-limit1}
 	Suppose that Assumptions \ref{as-moment-cinf},\ref{as-mp},\ref{as-sigma0-cinf} hold. Almost surely, the ESD \(F^{\mathbf{A}_n}\) of \(\mathbf{A}_n\) converges weakly to a (non-random) probability measure \(F\), the Stieltjes transform \(s = s(z)\) of which is the unique solution to the equation
 	\begin{align}\label{mp1} 
 	z = -\frac{1}{s} - \frac{s}{b} \int \frac{t^2}{1 + st/\sqrt{cb}} \, dH(t), \quad z \in \mathbb{C}^+,
 	\end{align}
 	in the set \(\left\{s: z \in \mathbb{C}^+,~ s \in \mathbb{C}^{+}\right\}\), where \(b = \int t^2 \, dH(t)\). 
 \end{lemma}
 
 \begin{remark}
 	Lemma \ref{th-limit1} provides a unified  LSD of \(\mathbf{A}_n\) when \(p/n\rightarrow c \in (0, \infty]\). This result is consistent with the generalized MP law of \(\mathbf{S}_n\)  when \( p/n\rightarrow c \in (0, \infty)\) in \citep{MP67}. To determine the boundary of the support of the LSD $F$, we introduce the function \(\phi_{c,H}(x)\):
 	\[
 	\phi_{c,H}(x) = x + \frac{1}{b} \int \frac{t^2}{x - t/\sqrt{cb}} \, dH(t),
 	\]
 	and define 
 	\[
 	\mathfrak{a} = \max_{x \in \mathbb{R}^+} \{x : \phi_{c,H}'(x) = 0\}.
 	\]
 	Then the right edge point \(\mathfrak{b}\) of the support is given by \(\mathfrak{b} = \phi_{c,H}(\mathfrak{a})\). When \(p/n \to \infty\), we have
 	\[
 	\phi_{c,H}(x) = x + \frac{1}{x}, \quad \mathfrak{a} = 1, \quad \text{and} \quad \mathfrak{b} = 2,
 	\]
 	corresponding to the case of the semicircle law in Lemma~\ref{th-limit1-cinf}.
 \end{remark}

 \begin{theorem}[Phase transition]
 	\label{th-limit}
 	Suppose that Assumptions \ref{as-2star}, \ref{as-mp},\ref{as-sigma0-cinf}-\ref{as-sigman-cinf} hold. The largest $\tau-1$ eigenvalues of $\A_n$ converge almost surely.
 	For $k \in[1:M]$,
 	\begin{itemize}
 		\item [(i)] if $\alpha_k>\mathfrak{a}$, $\lambda_{k\ell}^{\A_{n}}\xrightarrow{a.s.} \phi_{c,H}(\alpha_k)>\mathfrak{b},~\forall~ 1\leq \ell \leq m_k$ (distant spike);
 		\item [(ii)] if $\alpha_k\leq\mathfrak{a}$,
 		$\lambda_{k\ell}^{\A_{n}}\xrightarrow{a.s.} \mathfrak{b}, ~\forall~ 1\leq \ell \leq m_k$  (close spike).
 	\end{itemize}
 \end{theorem}

 \subsection{Unified CLT for distant spiked eigenvalues of $\A_n$}
 Under the unified framework $p/n\rightarrow c\in(0,\infty]$, 
the limiting distribution of $\mathbf{\Lambda}_{nk}$ defined in Section \ref{sec:CLT-spike}, involves one more auxiliary matrix, denoted as  \(\mathcal{Q}(\alpha)\),
 \[
 \mathcal{Q}(\alpha) = \boldsymbol{\Sigma}_0 - \sqrt{c_n b_p} \alpha \mathbf{I}_p.
 \]
Here \(\alpha\) is arbitrary real number satisfying \(\lim_{n \to \infty} \sqrt{c_n b_p} \alpha \notin \Gamma_H\) (the support of $H$), which guarantees that $\mathcal{Q}(\alpha)$ is invertible for all large \(p\) and \(n\).

 \begin{assumption}\label{as-limits}
 	As $n\to \infty$,
 	\begin{align*}
 	\bV_n( \alpha )&\triangleq \bU_n^{\top}\mathcal{Q}^{-1}(\alpha) \bU_n\to \bV( \alpha ),\quad
 	\bV_n^{\prime}( \alpha )\triangleq \sqrt{c_nb_p}\bU_n^{\top}\mathcal{Q}^{-2}(\alpha) \bU_n\to \bV^{\prime}( \alpha ).
 	\end{align*}
 	In addition,  
 	if $v_3\neq 0$, then
 	\begin{align}\label{hijl}
 	h_{n,ijl}(\alpha)\triangleq
 	&\sqrt{k_{n i} k_{n j} k_{n l}}
 	\sum_{q=1}^{p}\be_{q}^{\top}\bSigma_{0}^{\frac{1}{2}}\mathcal{Q}^{-1}(\alpha)\bmu_i
 	\be_{q}^{\top}\bSigma_{0}^{\frac{1}{2}}\mathcal{Q}^{-1}(\alpha)\bmu_j\be_{q}^{\top}\bSigma_{0}^{\frac{1}{2}}\mathcal{Q}^{-1}(\alpha)\bmu_l\nonumber \\
 	\to &h_{ijl}(\alpha),\quad \forall~ i, j, l \in[1:\tau];
 	\end{align}
 	if $v_4\neq 3$, then
 	\begin{align}\label{rho}
 	\rho_{n,ijlt}(\alpha)\triangleq
 	&\sqrt{k_{n i} k_{n j} k_{n l}k_{n t}}
 	\sum_{q=1}^{p}
 	\be_{q}^{\top}\bSigma_{0}^{\frac{1}{2}}\mathcal{Q}^{-1}(\alpha)\bmu_i\be_{q}^{\top}\bSigma_{0}^{\frac{1}{2}}\mathcal{Q}^{-1}(\alpha)\bmu_j\nonumber\\
 	&\times
 	\be_{q}^{\top}\bSigma_{0}^{\frac{1}{2}}\mathcal{Q}^{-1}(\alpha)\bmu_l\be_{q}^{\top}\bSigma_{0}^{\frac{1}{2}}\mathcal{Q}^{-1}(\alpha)\bmu_t\to \rho_{ijlt}(\alpha),\quad \forall~ i, j, l, t \in [1:\tau],
 	\end{align}
 	where $ \e_q \in \mathbb R^p$  denotes the unit vector with the $q$th element being 1 and all others being 0.
 \end{assumption}
 
 \begin{remark}
 	Assumption \ref{as-limits} states the existence of four limits, which will be used to define the limiting distribution of the spikes. In particular, the quantities in \eqref{hijl} and \eqref{rho} arise when dealing with non-Gaussian distributions of \(\mathbf{z}_{11}\), originating from the expectations of quadratic forms like
 	$\E\left(\mathbf{z}_{11}^{\top} \mathbf{M}_1 \mathbf{z}_{11}-\operatorname{tr} \mathbf{M}_1\right)\left(\mathbf{z}_{11}^{\top} \mathbf{M}_2 \mathbf{z}_{11}-\operatorname{tr} \mathbf{M}_2\right)$ 
 	and $\E\left\{\left(\mathbf{z}_{11}^{\top} \mathbf{M}_1 \mathbf{z}_{11}-\operatorname{tr} \mathbf{M}_1\right)\mathbf{z}_{11}^{\top} \mathbf{M}_1 \bmu_1\right\}$. When $p/n\rightarrow\infty$, these four limits degenerates to
 	\begin{align*}
 	\mathbf{V}(\alpha)&=
 	-\frac{1}{\alpha} \lim _{n \rightarrow \infty} 
 	\frac{\mathbf{U}_n^{\top} \mathbf{U}_n}{\sqrt{c_n b_p}},\quad
 	\mathbf{V}^{\prime}(\alpha)=
 	\frac{1}{\alpha^2} \lim _{n \rightarrow \infty} \frac{\mathbf{U}_n^{\top} \mathbf{U}_n}{\sqrt{c_n b_p}},\quad
 	h_{ijl}(\alpha)= \rho_{ijlt}(\alpha)\equiv 0,
 	\end{align*}
 	for $	i, j, l, t \in[1:\tau]$, which is consistent with Theorem~\ref{th-c-cinf}.
 	%
 	%
 	%
 \end{remark}
 

 \begin{theorem}\label{th-c}
 	Under Assumptions \ref{as-moment-cinf}, \ref{as-mp},\ref{as-sigma0-cinf}-\ref{as-limits},  the $m_{k}$-dimensional random vector
 	${\bf \Lambda}_{nk}$
 	converges in distribution to the joint distribution of the $m_{k}$ eigenvalues of the following Gaussian random matrix 
 	\begin{align*}
 	-\sqrt{ \phi_{c,H}^{\prime}\left(\alpha_{k}\right)}
 	\bQ_k^{\top} \bN \W \bN \bQ_k \mathbf{G},
 	\end{align*}
 	where $\bN$ is the limit of $\bN_n$ defined in \eqref{con:D}, $\bQ_k$ is a $\tau\times m_k$ matrix such that 
 	\begin{align*}
 	\bQ_k^{\top}\bQ_k={\bf I}_{m_k}\quad \text{and}\quad \bQ_k^{\top}\bN\bV(\alpha_k) \bN\bQ_k=-\I_{m_k},
 	\end{align*}
 	$ \mathbf{G} $ is the inverse matrix of $ \bQ_k^{\top} \bN\bV^{\prime}(\alpha_k) \bN\bQ_k$,	and
 	$\W=(W_{i j})$ is a symmetric $ \tau\times \tau $ random matrix with zero-mean Gaussian entries. Their covariances are, for $1 \leq i \neq j\neq l \neq t \leq \tau$,
 	\begin{align*}
 	\operatorname{Var}(W_{ii}) &=2\left[2\theta_{ii}+\theta_{ii}^2+1-\phi_{c,H}^{\prime}\left(\alpha_{k}\right)\right]+\phi_{c,H}^{\prime}\left(\alpha_{k}\right)(v_4-3)\rho_{iiii}, \\
 	\operatorname{Var}(W_{ij}) &=\theta_{ii}+\theta_{jj}+\theta_{ij}^2+\theta_{jj}\theta_{ii}+1-\phi_{c,H}^{\prime}\left(\alpha_{k}\right)+\phi_{c,H}^{\prime}\left(\alpha_{k}\right)(v_4-3)\rho_{ijij},\\
 	\operatorname{Cov}(W_{ii},W_{it})&=2\left[\theta_{it}+\theta_{it}\theta_{ii}\right]+\phi_{c,H}^{\prime}\left(\alpha_{k}\right)(v_4-3)\rho_{iiit}, \\
 	\operatorname{Cov}(W_{ii},W_{ll})&=2\theta_{il}^2+\phi_{c,H}^{\prime}\left(\alpha_{k}\right)(v_4-3)\rho_{iill},\\
 	\operatorname{Cov}(W_{ii},W_{lt})&=2\theta_{ti}\theta_{il}+\phi_{c,H}^{\prime}\left(\alpha_{k}\right)(v_4-3)\rho_{iilt},\\
 	\operatorname{Cov}(W_{ij},W_{lj})&=\theta_{il}+\theta_{jl}\theta_{ji}+\theta_{jj}\theta_{il}+\phi_{c,H}^{\prime}\left(\alpha_{k}\right)(v_4-3)\rho_{ijlj},\\
 	\operatorname{Cov}(W_{ij},W_{lt})&=\theta_{jl}\theta_{ti}+\theta_{jt}\theta_{li}+\phi_{c,H}^{\prime}\left(\alpha_{k}\right)(v_4-3)\rho_{ijlt},
 	\end{align*}
 	where 
 	$$\theta_{ij}=\theta_{ij}\left(\alpha_k\right)=\lim_{n \to\infty}\sqrt{k_{ni} k_{nj}} \boldsymbol{\mu}_{i}^{\top}\mathcal{Q}^{-1}(\alpha_k) \bSigma_0\mathcal{Q}^{-1}(\alpha_k)\boldsymbol{\mu}_{j}$$ 
 	and $\rho_{ijlt}=\rho_{ijlt}\left(\alpha_k\right)$  defined in  \eqref{rho}. 
 \end{theorem}

 \begin{remark}
 	Theorem \ref{th-c} establishes the asymptotic distribution of the sample distant spikes \(\{\lambda_{k\ell}^{\mathbf{A}_n} : 1\leq \ell\leq  m_k\}\) when \(p/n\rightarrow c \in (0, \infty]\).  The covariances between entries of \(\mathbf{W}\) depend on  \(\theta_{ij}\) and \(\rho_{ijlt}\), whose existence is guaranteed by Assumption \ref{as-limits}. These quantities are functions of the mixing weights \(\{k_{ni}\}\) and the inner products of the means \(\{\boldsymbol{\mu}_i\}\) and the eigenvectors of \(\boldsymbol{\Sigma}_0\). When $p/n\rightarrow\infty$, Theorem \ref{th-c} reduces to Theorem \ref{th-c-cinf} with $\theta_{i j}=\rho_{i j l t}=0$ and  $\mathbf{G}=\alpha_k \mathbf{I}_\tau$. 
 \end{remark}

 \begin{theorem}\label{le-hatA}
 Theorem \ref{th-limit} and Theorem \ref{th-c} still hold for the approximated matrix \(\hat{\mathbf{A}}_n\) under the same assumptions.
 \end{theorem}

 \subsection{Unified asymptotic properties for random sesquilinear forms}\label{sec-ses}

 This section establishes the first- and joint second-order convergence of several random sesquilinear forms of \(\left(\mathbf{B}_n - \tilde{z}_n \mathbf{I}\right)^{-1}\) in \eqref{ses-infty} under the unified framework $p/n\rightarrow c\in(0,\infty]$. Here \(\tilde{z}_n\) is specified as
 \[
 \tilde{z}_n = c_n a_p + \sqrt{c_n b_p} z, \quad z \notin \Gamma_F,
 \]
 where \(\Gamma_F\) represents the support of \(F\) in Lemma \ref{th-limit1}. 
 These results are fundamental for proving Theorems \ref{th-limit}-\ref{th-c}, and possess independent interest for other inference procedures.
 

 \begin{theorem}\label{th-as}
 	Suppose that Assumptions \ref{as-moment-cinf},\ref{as-mp},\ref{as-sigma0-cinf}-\ref{as-eig0-cinf} hold and $\bmu_{i}^{\top}\bmu_i\asymp\sqrt{c_n}$ for all $ 1\leq i\leq \tau$, then for any $1\leq i,j\leq \tau$, we have
 	\begin{align*}
 	&\frac{\tilde{z}_n}{\sqrt{c_n b_p}} \bar{\mathbf{s}}_i^{\top}\left(\mathbf{B}_n-\tilde{z}_n \mathbf{I}\right)^{-1} \bar{\mathbf{s}}_j
 	+\frac{1}{k_{i}}\left[\sqrt{c_n/b_p}a_p+z+1/s(z)\right]I(i=j)
 	\xrightarrow{a . s .} 0,
 	\\
 	&\frac{\tilde{z}_n}{\sqrt{c_n b_p}} \bar\s_i^{\top}\left(\mathbf{B}_n-\tilde{z}_n \mathbf{I}\right)^{-1} \boldsymbol{\mu}_j	\xrightarrow{a . s .} 0,\\
 	&\frac{\tilde{z}_n}{\sqrt{c_n b_p}} \bmu_i^{\top}\left(\mathbf{B}_n-\tilde{z}_n \mathbf{I}\right)^{-1} \boldsymbol{\mu}_j+\bmu_i^{\top}\left(\sqrt{c_nb_p}\bI+s(z)\bSigma_0\right)^{-1}\bmu_j\xrightarrow{a . s .} 0,
 	\end{align*}
 	where $I(\cdot)$ denotes the indicator function.
 \end{theorem}

 \begin{theorem}\label{th-clt}
 	Suppose that Assumptions \ref{as-moment-cinf},\ref{as-mp},\ref{as-sigma0-cinf}-\ref{as-eig0-cinf}  and \ref{as-limits} hold with \(\boldsymbol{\mu}_i^{\top} \boldsymbol{\mu}_i \asymp \sqrt{c_n}\). Then,  \(\mathbf{L}_n\) in \eqref{eq:Ln} converges in distribution to a symmetric \(2\tau \times 2\tau\) random matrix \(\mathbf{L} = (L_{ij})\) of Gaussian entries with zero mean and covariance
 	\begin{itemize}
 		\item[(1)] for $ i,j,l,t\in[1:\tau] $,
 		\begin{align*}
 		\Cov(L_{ij},L_{lt})=\begin{cases}
 		2\frac{s^{\prime}\left(z\right)-s^{2}\left(z\right)}{s^4(z)}\frac{1}{k_i^2}&\mbox{if $ i=j=l=t $,}\\
 		\frac{s^{\prime}\left(z\right)-s^{2}\left(z\right)}{s^4(z)}\frac{1}{k_ik_j}&\mbox{if $ i=l\neq j=t $,}\\
 		0&\mbox{o.w.;}
 		\end{cases}
 		\end{align*}
 		
 		\item[(2)] for $ i,j,l,t\in[\tau+1:2\tau]$,
 		\begin{align*}
 		\Cov(L_{ij},L_{lt})&=\frac{s^{\prime}(z)}{s^{4}(z)}\left(\zeta_{(j-\tau)(l-\tau)} \zeta_{(t-\tau)( i-\tau)}+\zeta_{(j-\tau)(t-\tau)} \zeta_{(l-\tau)(i-\tau)}\right)\\
 		&\quad+\frac{v_4-3}{s^2(z)}g_{(i-\tau)(j-\tau)(l-\tau)(t-\tau)};
 		\end{align*}
 		
 		\item[(3)] for $ i\in[1:\tau] $ and $ t\in [\tau+1:2\tau]$,
 		\begin{align*}
 		\Cov(L_{ij},L_{lt})=\begin{cases}
 		\frac{s^{\prime}(z)}{s^4(z)}\frac{\zeta_{(j-\tau)(t-\tau)}}{k_i}&\mbox{if $i=l$, $j\in[\tau+1:2\tau]$,}\\
 		-\frac{v_3}{s^2(z)}f_{(j-\tau)(l-\tau)(t-\tau)}&\mbox{if $j, l\in[\tau+1: 2 \tau]$,}\\
 		0&\mbox{o.w.}
 		\end{cases}
 		\end{align*}
 		%
 		
 		
 		
 	\end{itemize}
 	where $$ \zeta_{ij}=\frac{\theta_{ij}(-1/s(z))}{\sqrt{k_ik_j}}, f_{ijl}=\frac{h_{ijl}(-1/s(z))}{\sqrt{k_ik_jk_l}}, g_{ijlt}=\frac{\rho_{ijlt}(-1/s(z))}{\sqrt{k_ik_jk_lk_t}}.$$
 \end{theorem}
 \begin{remark}
 	Theorem \ref{th-clt} establishes the asymptotic distribution of 
 	the random matrix \(\mathbf{L}_n\) when  \(p/n\rightarrow c \in (0, \infty]\). The limiting distribution is jointly determined by \(v_3\), \(v_4\), \(s(z)\) and \(\{\zeta_{ij}, f_{ijl}, g_{ijlt}\}\). Notably, it can be seen that \(\bar{\mathbf{s}}_i^{\top}\left(\mathbf{B}_n - \tilde{z}_n \mathbf{I}\right)^{-1} \bar{\mathbf{s}}_j\), \(\bar{\mathbf{s}}_i^{\top}\left(\mathbf{B}_n - \tilde{z}_n \mathbf{I}\right)^{-1} \boldsymbol{\mu}_j\), and \(\bar{\mathbf{s}}_j^{\top}\left(\mathbf{B}_n - \tilde{z}_n \mathbf{I}\right)^{-1} \boldsymbol{\mu}_i\) are asymptotically independent. \(\bar{\mathbf{s}}_i^{\top}\left(\mathbf{B}_n - \tilde{z}_n \mathbf{I}\right)^{-1} \bar{\mathbf{s}}_j\) is also asymptotically independent of \(\boldsymbol{\mu}_i^{\top}\left(\mathbf{B}_n - \tilde{z}_n \mathbf{I}\right)^{-1} \boldsymbol{\mu}_j\).
 \end{remark}

    \section{Examples and simulations}
    \label{sec:3}
    In this section, we provide examples and simulations for Theorem \ref{th-c} under the unified framework where $p/n\rightarrow c\in(0,\infty]$.
    
    \subsection{Example I} 
    Consider the two-sample case when \(\tau = 2\). Note that \(\mathbf{P}_n\) is of rank 1 in \(\mathbf{S}_n = \tilde{\mathbf{S}}_n + \mathbf{P}_n\), thus \(\mathbf{A}_n\) has at most one distant spiked eigenvalue in this case.

    
    \begin{corollary}\label{co-tau=2}
    	Suppose that Assumptions \ref{as-moment-cinf}, \ref{as-mp},\ref{as-sigma0-cinf}-\ref{as-limits} hold  for $\tau=2$.	If the largest eigenvalue of $\bSigma_{\bx}$ is a distant spike, i.e., $\alpha_1>\mathfrak{a}$, then 
    	$$
    	\sqrt{n}\left(\lambda_1^{\A_{n}}-\lambda_{n1}\right)\xrightarrow{D} N\left(0, \sigma^{2}\right),
    	$$
    	where the limiting variance is
    	\begin{align*}
    	\sigma^2=	2\phi_{c,H}^{\prime}\left(\alpha_1\right)\alpha_1^2\left[1-\frac{\phi_{c,H}^{\prime}\left(\alpha_1\right)}{(1+\omega)^2}
    	+\frac{\phi_{c,H}^{\prime}\left(\alpha_1\right)}{2(1+\omega)^2}(v_4-3)\eta
    	\right],
    	\end{align*}
    	with	
    	\begin{align*}
    	\omega&=\lim_{n \to \infty} k_1k_2\left(\bmu_1-\bmu_2\right)^{\top}\left(\bSigma_0-\sqrt{c_nb_p}\alpha_1\bI\right)^{-1}\bSigma_0\left(\bSigma_0-\sqrt{c_nb_p}\alpha_1\bI\right)^{-1}\left(\bmu_1-\bmu_2\right)\\
    	&=k_2\theta_{11}+k_1\theta_{22}-2\sqrt{k_1k_2}\theta_{12},\\
    	\eta&=\lim_{n \to \infty} k_1^2k_2^2\sum_{q=1}^{p}\left[\e_q^{\top}\bSigma_0^{\frac{1}{2}}\left(\bSigma_0-\sqrt{c_nb_p}\alpha_1\bI\right)^{-1}\left(\bmu_1-\bmu_2\right)\right]^4\\
    	&=k_2^2\rho_{1111}+k_1^2\rho_{2222}+6k_1k_2\rho_{1212}-4k_2\sqrt{k_1k_2}\rho_{1112}-4k_1\sqrt{k_1k_2}\rho_{1222}.
    	\end{align*}
    	Especially, when $p/n\rightarrow\infty$, $\sigma^2=2\left(1-\alpha_1^{-2}\right)$.
    \end{corollary}
    
    \subsection{Example II} 	
    Consider the case where the base covariance matrix \(\boldsymbol{\Sigma}_0 = \mathbf{I}_p\), then  \(\boldsymbol{\Sigma}_{\mathbf{x}}\) simplifies to
    \[
    \boldsymbol{\Sigma}_{\mathbf{x}} = \sqrt{\frac{n}{p}}\left(\mathbf{I}_p + \boldsymbol{\Sigma}_{\tau}\right).
    \]
    The \(k\)th cluster of eigenvalues of \(\boldsymbol{\Sigma}_{\mathbf{x}}\) are distant spikes if and only if \(\alpha_k > 1 + 1/\sqrt{c}\). 
    Accordingly, the \(k\)th cluster of the sample eigenvalues of \(\mathbf{A}_n\) converges jointly to the spectrum of a Gaussian matrix.

    \begin{corollary}\label{co-danwei}
    	Suppose that Assumptions \ref{as-moment-cinf}, \ref{as-mp},\ref{as-sigma0-cinf}-\ref{as-limits} hold with  $\bSigma_0=\I_p$ and \(\alpha_k > 1 + 1/\sqrt{c}\).
    	Then,
    	the $m_{k}$-dimensional random vector ${\bf \Lambda}_{nk}$
    	converges in distribution to the joint distribution of the $m_{k}$ eigenvalues of the Gaussian random matrix
    	$$
    	-\left[\left(\alpha_k-\frac{1}{\sqrt{c}}\right)^2-1\right]^{\frac{1}{2}}\bQ_k^{\top} \bN \W \bN \bQ_k 
    	$$
    	where $ \W $ is defined in Theorem \ref{th-c} with $ \phi_{c,H}^{\prime}\left(\alpha_{k}\right)=1- c/(\sqrt{c}\alpha_{k}-1)^2$ and 
    	\begin{align*}
    	\theta_{i j}&=\lim_{n \to \infty} \boldsymbol{\mu}_{i}^{\top} \boldsymbol{\mu}_{j}\sqrt{k_{i} k_{j}}(\sqrt{c_n}\alpha_k-1)^{-2}, \\
    	\rho_{i jlt}&=\lim_{n \to \infty}\sum_{q=1}^{p}\e_q^{\top}\bmu_i\e_q^{\top}\bmu_j\e_q^{\top}\bmu_l\e_q^{\top}\bmu_t\sqrt{k_{i} k_{j}k_lk_t}(\sqrt{c_n}\alpha_k-1)^{-4},
    	\quad i, j,l,t\in[1:\tau]. 
    	\end{align*}
    	Especially when $p/n\rightarrow\infty$, $ \alpha_k\W $ is a $\tau\times \tau$ Wigner matrix.
    \end{corollary}
    
    In addition, for the simplest case where $\bSigma_0=\mathbf{I}$ and $\tau=2$, the only spiked eigenvalue of $\bSigma_{\bx}$ is 
    $$\alpha_{n1}=\sqrt{\frac{n}{p}}\left[1+k_{n1}k_{n2}\left(\boldsymbol{\mu}_{1}-\boldsymbol{\mu}_{2}\right)^{\top}\left(\boldsymbol{\mu}_{1}-\boldsymbol{\mu}_{2}\right)\right].
    $$
    If $\alpha_{n1}\to\alpha_1>1+1/\sqrt{c}$,
    then
    $$
    \sqrt{n}\left[\lambda_{1}^{\A_{n}}-\lambda_{n1}\right] \xrightarrow{D} N\left(0, \sigma^{2}\right),
    $$
    where 
    $\lambda_{n1}=\alpha_{n1}+1/\left(\alpha_{n1}-\sqrt{n/p}\right)$ and
    \begin{align*}
    \sigma^2=2\left(1+\frac{2\alpha_1}{\sqrt{c}}-\frac{1}{c}\right)
    \left[1-\frac{1}{\left({\alpha_1-1/\sqrt{c}}\right)^2}\right]+
    (v_4-3)\left[\alpha_1-\frac{1}{\sqrt{c}}-\frac{1}{\alpha_1-1/\sqrt{c}}\right]^2\eta
    \end{align*}
    with $ \eta=k_1^2k_2^2\lim_{n \to \infty}\sum_{q=1}^{p}\left[ \e_q^{\top}\left(\bmu_1-\bmu_2\right) \right]^4\left(\sqrt{c_n} \alpha_k-1\right)^{-4}$.
    Especially, when $p/n\rightarrow\infty, \sigma^2=2\left(1-\alpha_1^{-2}\right)$.

    \subsection{Numerical results}
    In this section, we examine the numerical performance of CLT for distant spikes of \(\mathbf{A}_n\). The following multi-sample setting with \(\tau = 4\) is employed:
    $$
    \mathbf{x}_{i j}=\boldsymbol{\mu}_i+\mathbf{\Sigma}_0^{\frac{1}{2}} \mathbf{z}_{i j}, \quad 
    \bSigma_0=\operatorname{diag}(\underbrace{1, \ldots, 1}_{p / 2}, \underbrace{2, \ldots, 2}_{p / 2}),
    $$
    for $i\in[1:4]$ and $j\in[1:n_i], ~n_i=n/4$.     The dimensional settings are \((p, n)=(5000, 10^4), (4 \times 10^4, 4000),(10^6, 2000)\) with \(c_n = 0.5, 10, 500\).    
    The mean vectors are set to be 
    \begin{align*}
    \boldsymbol{\mu}_{1}&={\bf 0},\quad
    \boldsymbol{\mu}_{2}=\left(c_nb_p\right)^{\frac{1}{4}}\big(4,0, \ldots,0\big)^{\top},\\
    \boldsymbol{\mu}_{3}&=\left(c_nb_p\right)^{\frac{1}{4}}\big(2,3\sqrt{2},0, \ldots,0\big)^{\top},\quad
    \boldsymbol{\mu}_{4}=\left(c_nb_p\right)^{\frac{1}{4}}\big(2,\sqrt{2},-4, 0, \ldots,0\big)^{\top},
    \end{align*}
    where $b_p=\tr(\bSigma_0^2)/p=2.5$.
    The random variables \((z_{ijq})\) are generated from
    \begin{itemize}
    	\item[] {\bf Case I:} ~~\((z_{ijq})\) i.i.d.  \( \text{Exp}(1) - 1 \) with \( v_3 = 2 \) and \( v_4 = 9 \).
    	\item[] {\bf Case II:} \((z_{ijq})\)  i.i.d.  \(\{ \text{Bernoulli}(t) - t \} / \sqrt{t(1 - t)}\) with \( t = \frac{\sqrt{3} + 3}{6} \), \( v_3 = -\sqrt{2} \),  \( v_4 = 3 \).
    \end{itemize}

    Under this model, 
    $$
    \boldsymbol{\Sigma}_\bx=\sqrt{\frac{2n}{5p}}\left(\bSigma_0+\frac{1}{16} \sum_{1 \leq i<j \leq 4}\left(\boldsymbol{\mu}_i-\boldsymbol{\mu}_j\right)\left(\boldsymbol{\mu}_i-\boldsymbol{\mu}_j\right)^{\top}\right),
    $$
    with spikes $\alpha_{n, 11}=\alpha_{n, 12}=3+\sqrt{2/(5c_n)}$ and $\alpha_{n, 21}=2+\sqrt{2/(5c_n)}$. These eigenvalues are all distant spikes, and
    accordingly, the largest $\tau-1$ eigenvalues of $\A_n$ are also spikes, denoted as $\lambda_{j}^{\A_{n}}$ for $j\in[1:3]$.
    Let
    \begin{align}\label{delta}
    \delta_{nj}
    \triangleq 
    \begin{cases}
    \sqrt{n}\left(\lambda_{j}^{\A_{n}}-\lambda_{n1}\right )& j=1,2,\\
    \sqrt{n}\left(\lambda_{j}^{\A_{n}}-\lambda_{n2}\right )& j=3,
    \end{cases}\quad\text{and}\quad
    \hat\delta_{nj}\triangleq
    \begin{cases}
    \sqrt{n}\left(\lambda_{j}^{\hat\A_{n}}-\sqrt{b_p/\hat b_p}\lambda_{n1}\right )& j=1,2,\\
    \sqrt{n}\left(\lambda_{j}^{\hat\A_{n}}-\sqrt{b_p/\hat b_p}\lambda_{n2}\right )& j=3,
    \end{cases}  
    \end{align}
    where
    $$
    \lambda_{n1}=\frac{46}{15}+\sqrt{\frac{2}{5c_n}}+\frac{4}{5}\frac{1}{3-\sqrt{2/(5c_n)}},\quad \lambda_{n2}=2.1+\sqrt{\frac{2}{5c_n}}+\frac{4}{5}\frac{1}{2-\sqrt{2/(5c_n)}}.
    $$
     From Theorems \ref{th-c}, the vector \((\delta_{n1}, \delta_{n2})^{\top}\) converges in distribution to the spectrum of a Gaussian random matrix, and \(\delta_{n3}\) converges in distribution to a Gaussian variable. In particular, we use $\sqrt{b_p / \hat{b}_p}$ in $\hat\delta_{nj}$ as a  finite sample correction. Actually, $\hat\delta_{nj}$  share the same asymptotic distribution with $\delta_{nj}$ because as $p/n\rightarrow c\in(0,\infty]$, $n\to\infty$,  $\sqrt{n}(\sqrt{b_p / \hat{b}_p}-1)=o_p(1)$ (see proof in Section S6 in the Supplementary Material).

\begin{figure}
    	\rotatebox{0}{\tiny{Case I for $\A_n$}}
    	\begin{subfigure}{.1\textwidth}
    		\centering
    		\includegraphics[width=\linewidth]{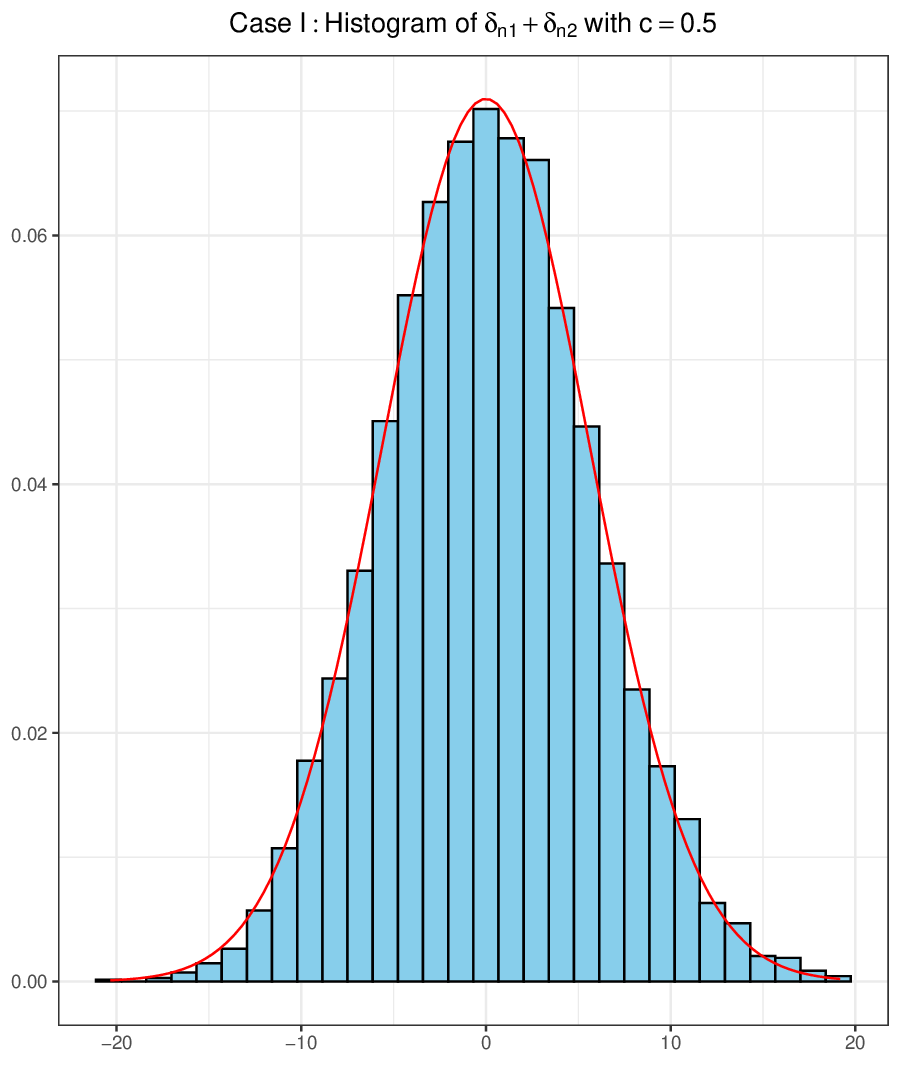}
    	\end{subfigure}
    	\begin{subfigure}{.1\textwidth}
    		\centering
    		\includegraphics[width=\linewidth]{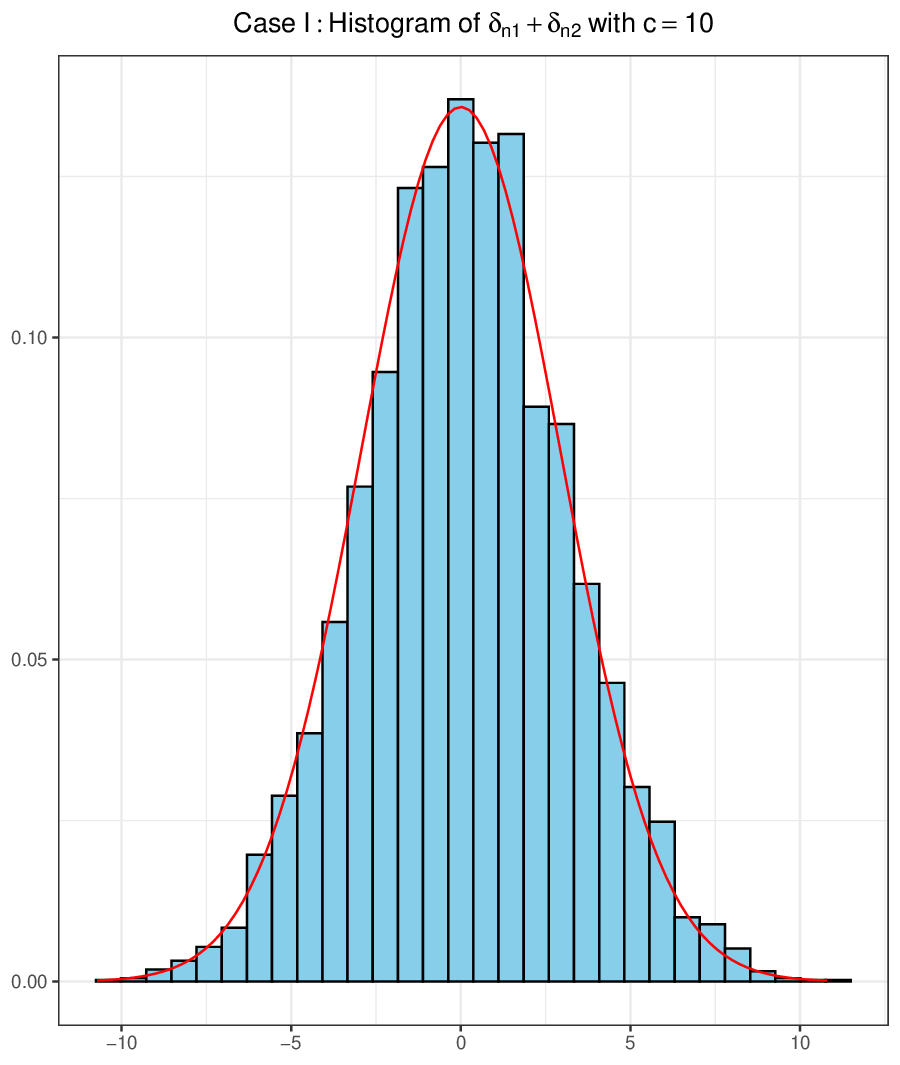}
    	\end{subfigure}
    	\begin{subfigure}{.1\textwidth}
    		\centering
    		\includegraphics[width=\linewidth]{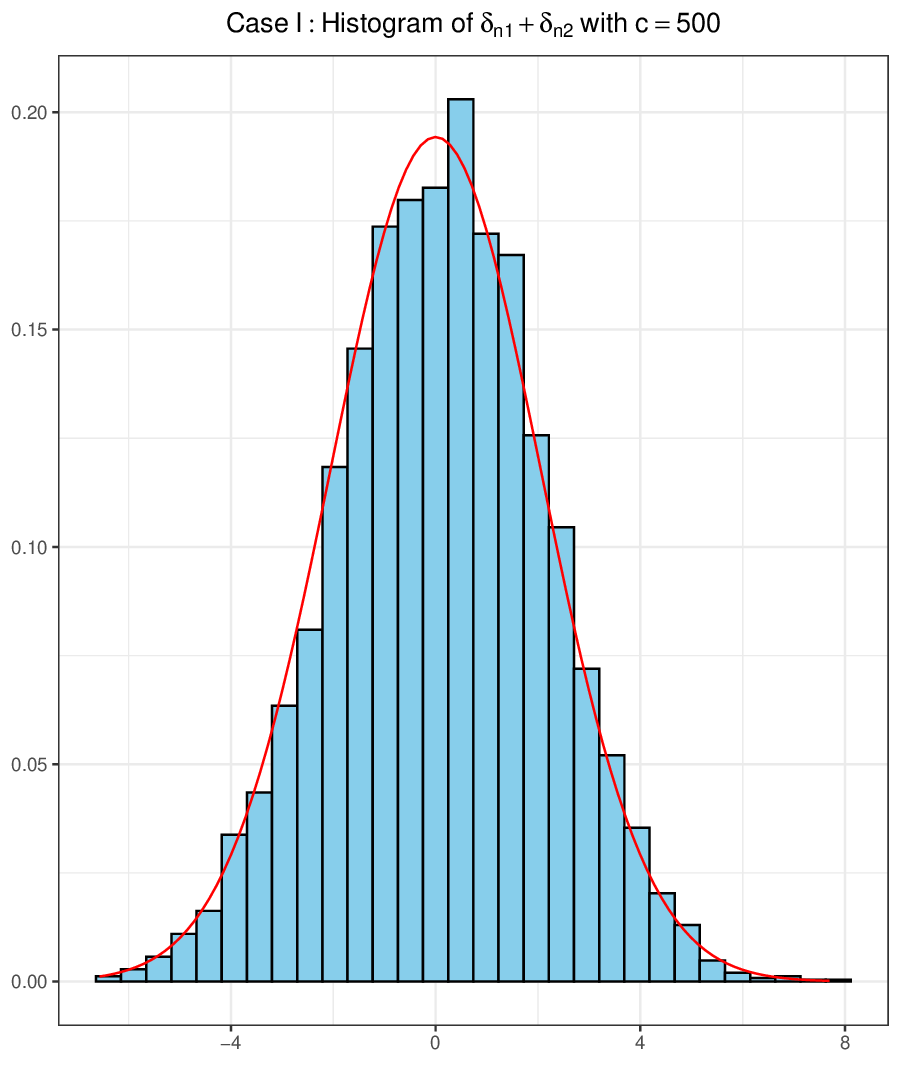}
    	\end{subfigure}
    	\qquad
    	\begin{subfigure}{.1\textwidth}
    		\centering
    		\includegraphics[width=\linewidth]{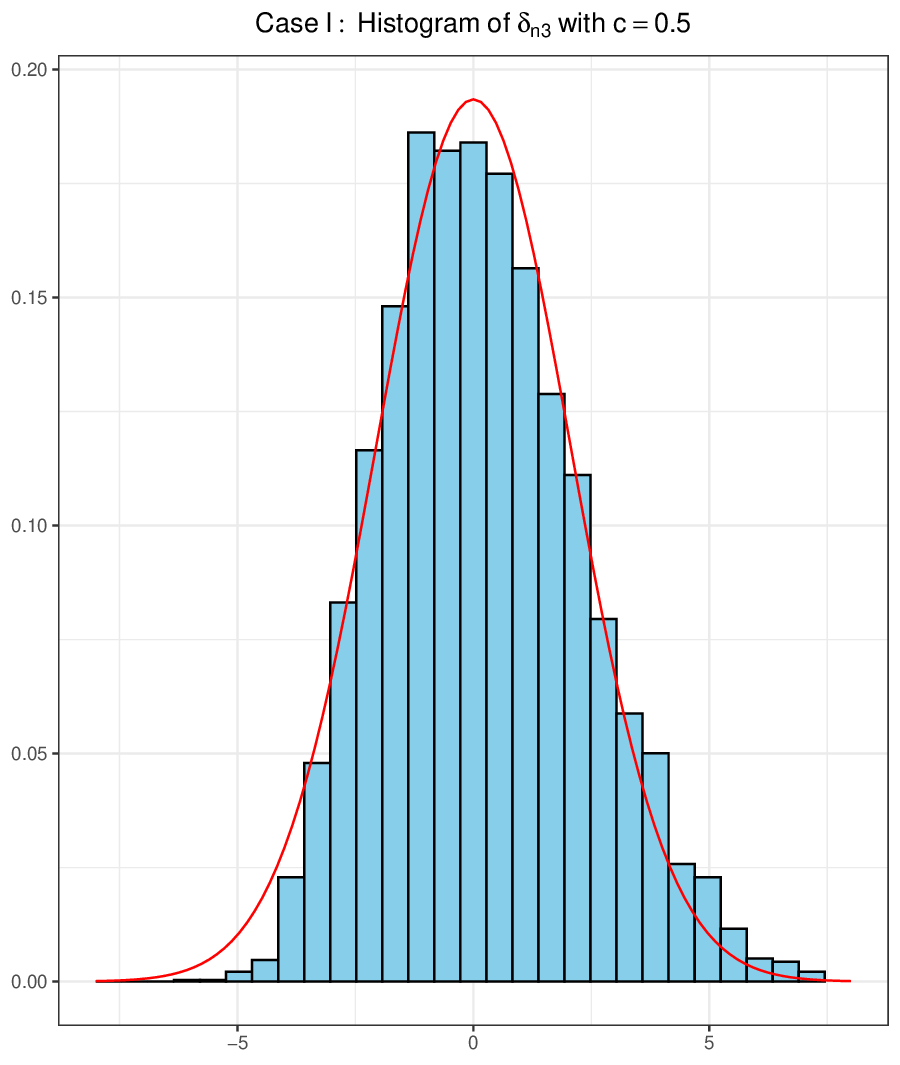}
    	\end{subfigure}
    	\begin{subfigure}{.1\textwidth}
    		\centering
    		\includegraphics[width=\linewidth]{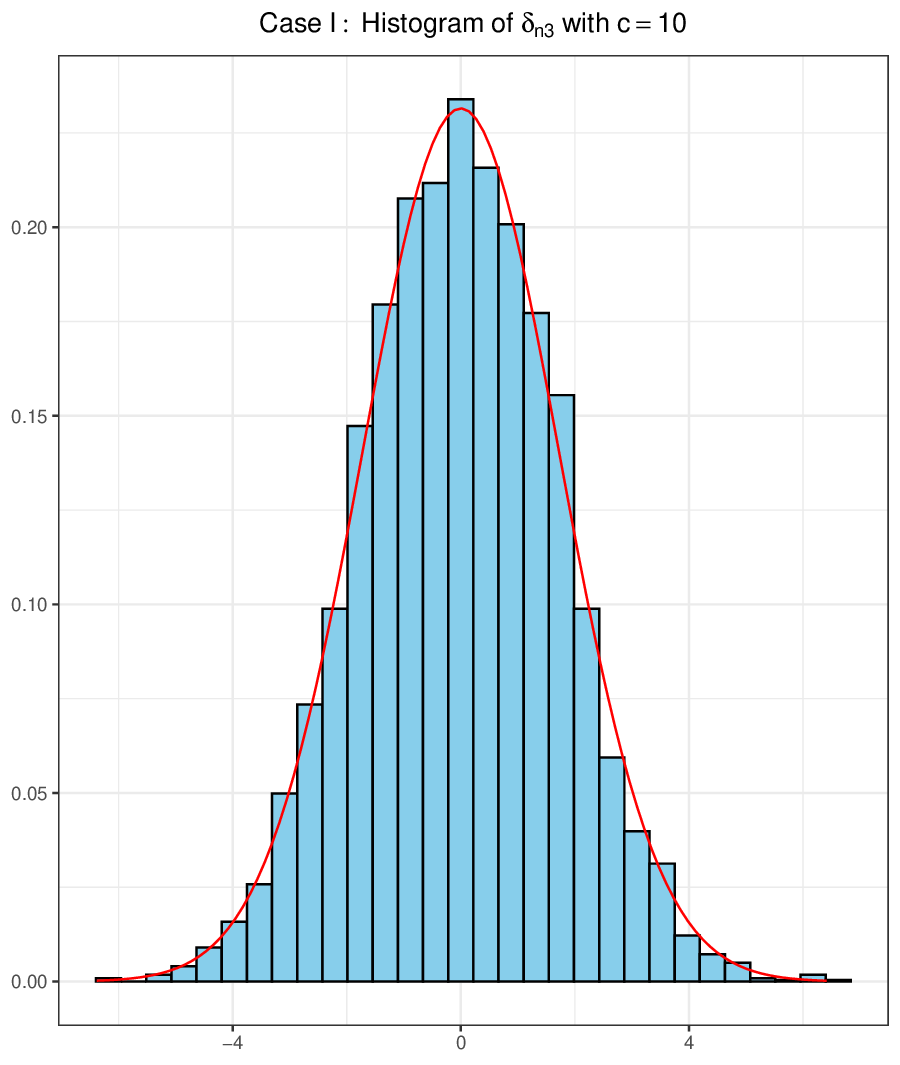}
    	\end{subfigure}
    	\begin{subfigure}{.1\textwidth}
    		\centering
    		\includegraphics[width=\linewidth]{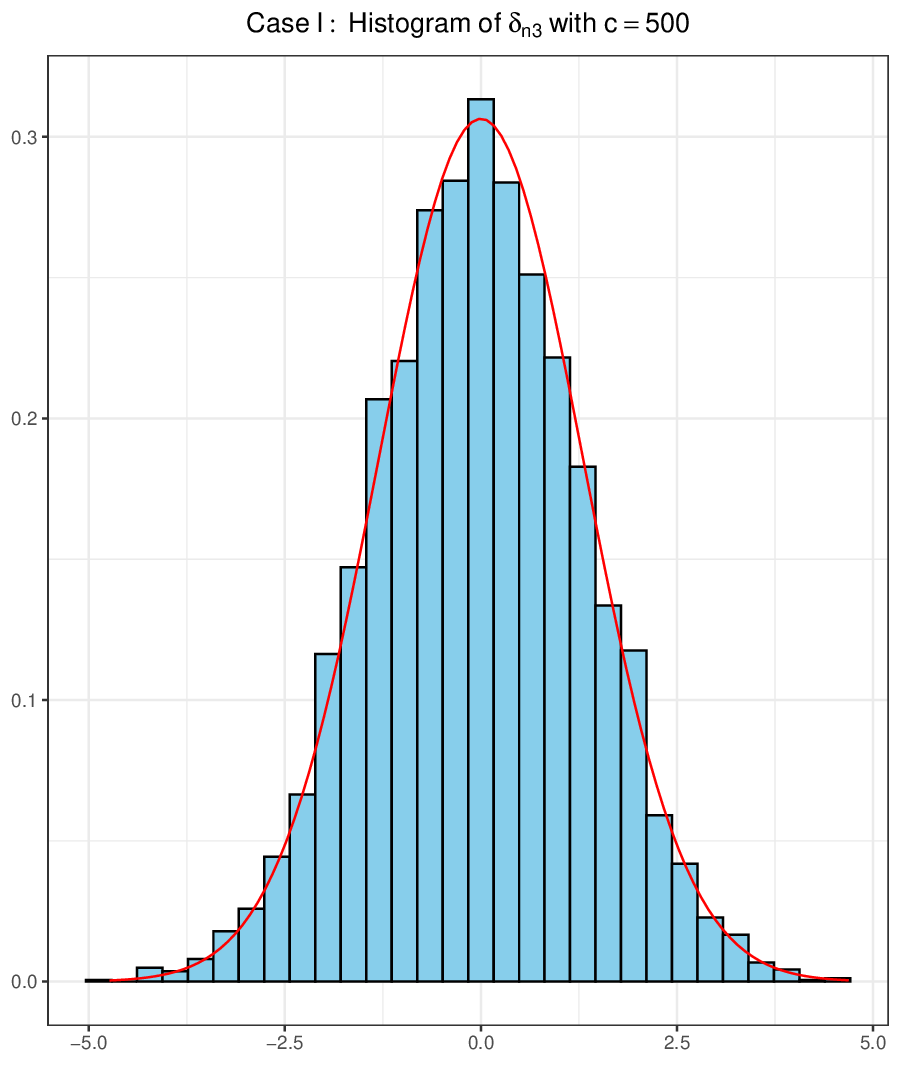}
    	\end{subfigure}
    	\\
    	\rotatebox{0}{\tiny{Case I for $\hat\A_n$}}
    	\begin{subfigure}{.1\textwidth}
    		\centering
    		\includegraphics[width=\linewidth]{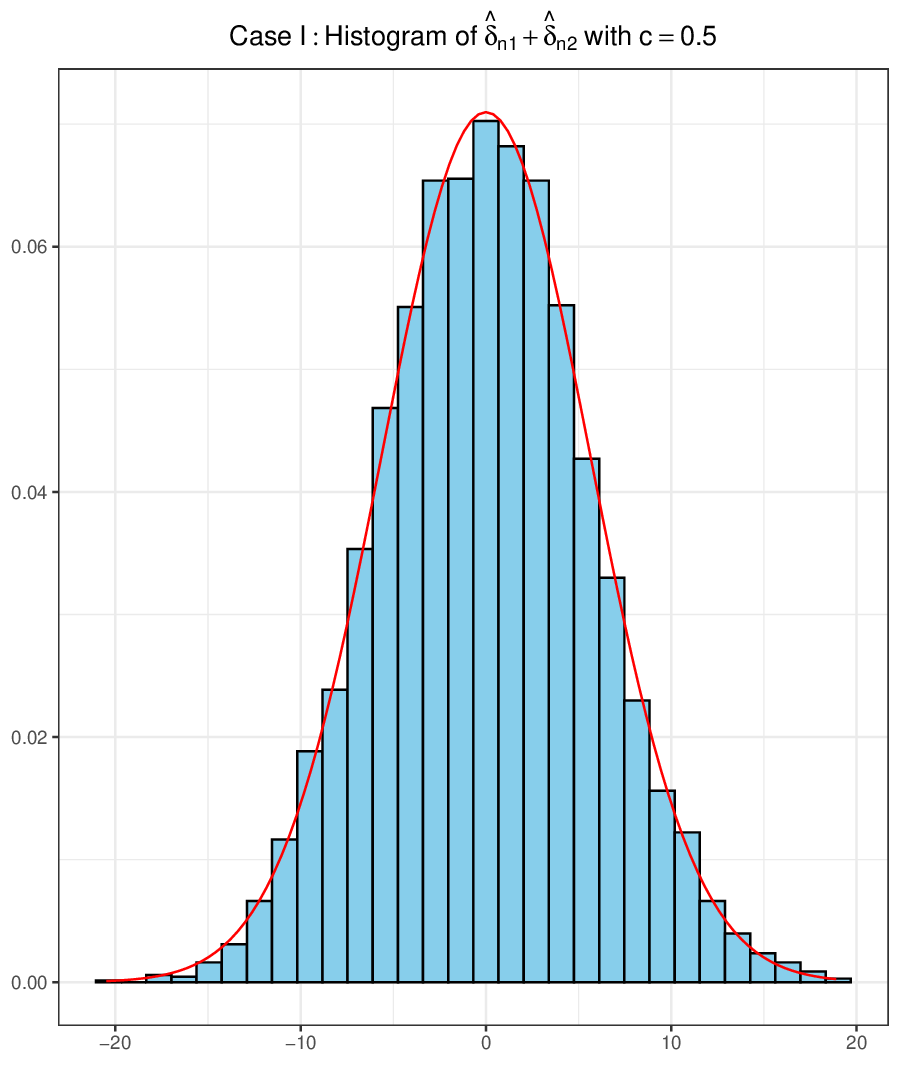}
    	\end{subfigure}
    	\begin{subfigure}{.1\textwidth}
    		\centering
    		\includegraphics[width=\linewidth]{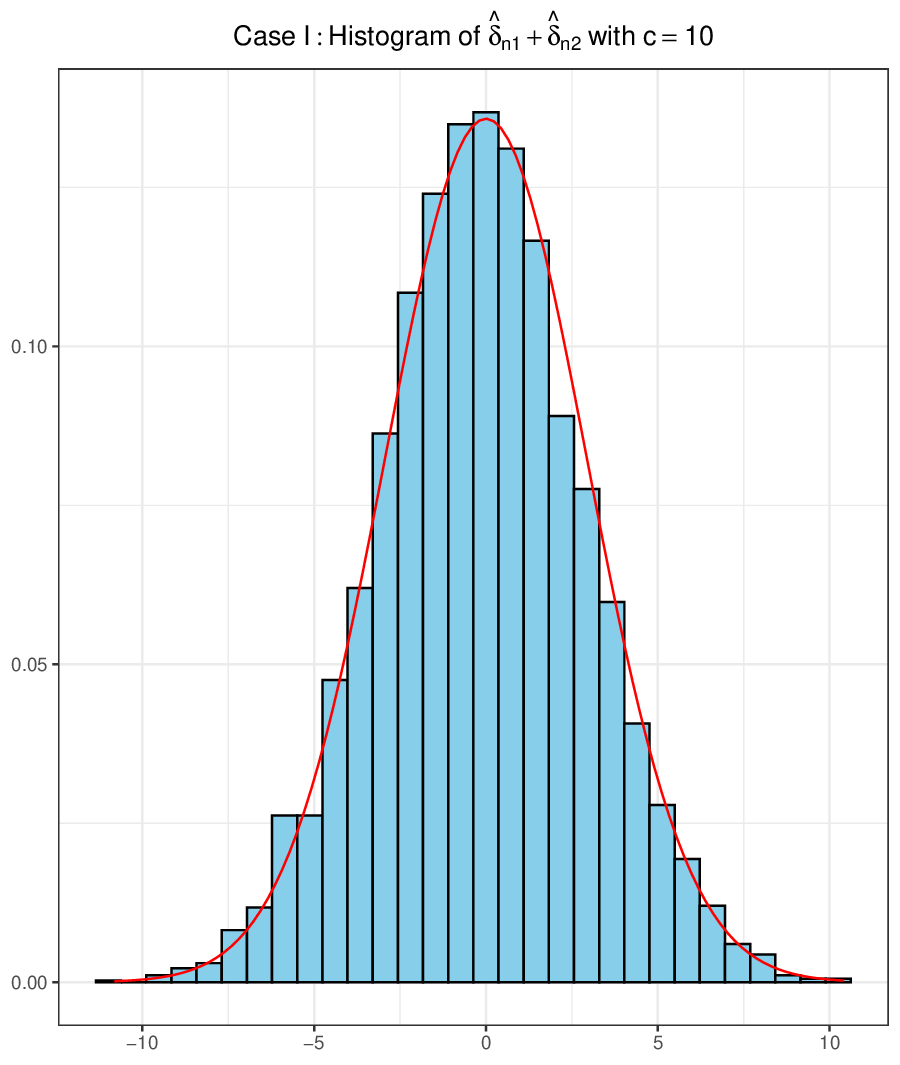}
    	\end{subfigure}
    	\begin{subfigure}{.1\textwidth}
    		\centering
    		\includegraphics[width=\linewidth]{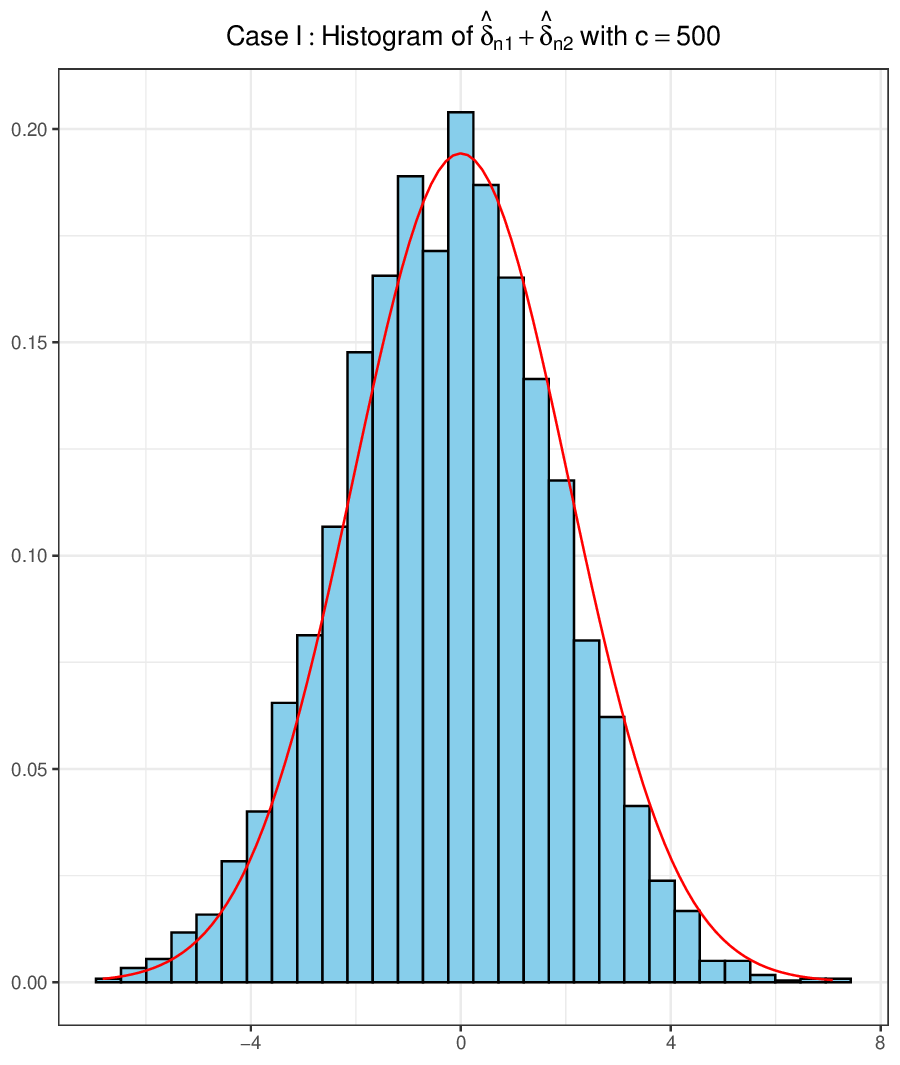}
    	\end{subfigure}
    	\qquad
    	\begin{subfigure}{.1\textwidth}
    		\centering
    		\includegraphics[width=\linewidth]{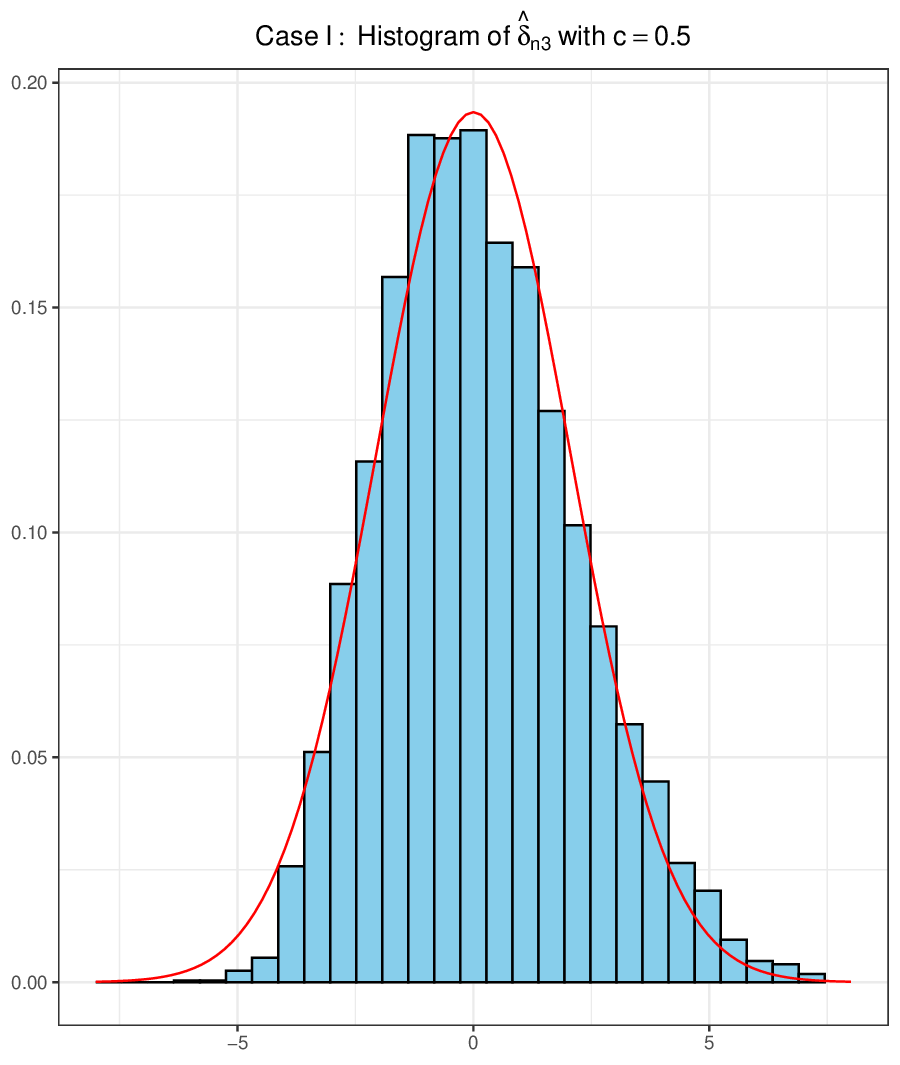}
    	\end{subfigure}
    	\begin{subfigure}{.1\textwidth}
    		\centering
    		\includegraphics[width=\linewidth]{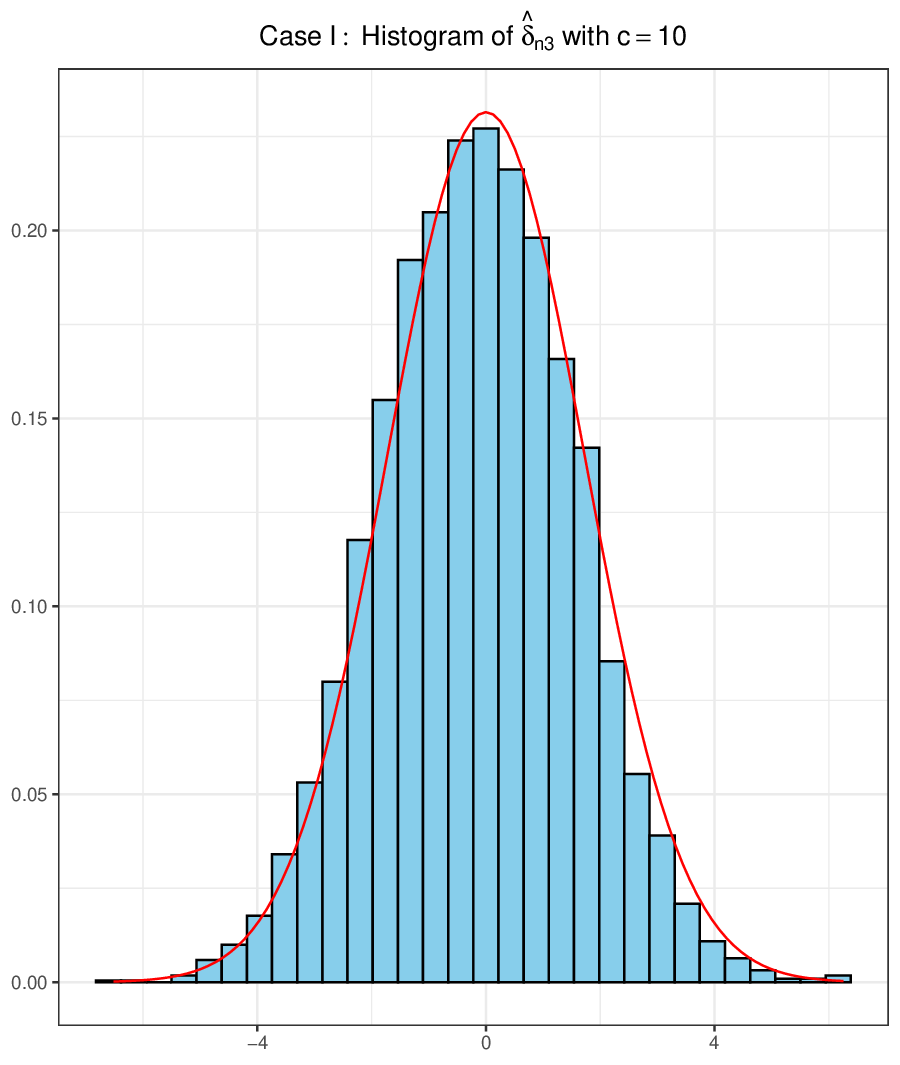}
    	\end{subfigure}
    	\begin{subfigure}{.1\textwidth}
    		\centering
    		\includegraphics[width=\linewidth]{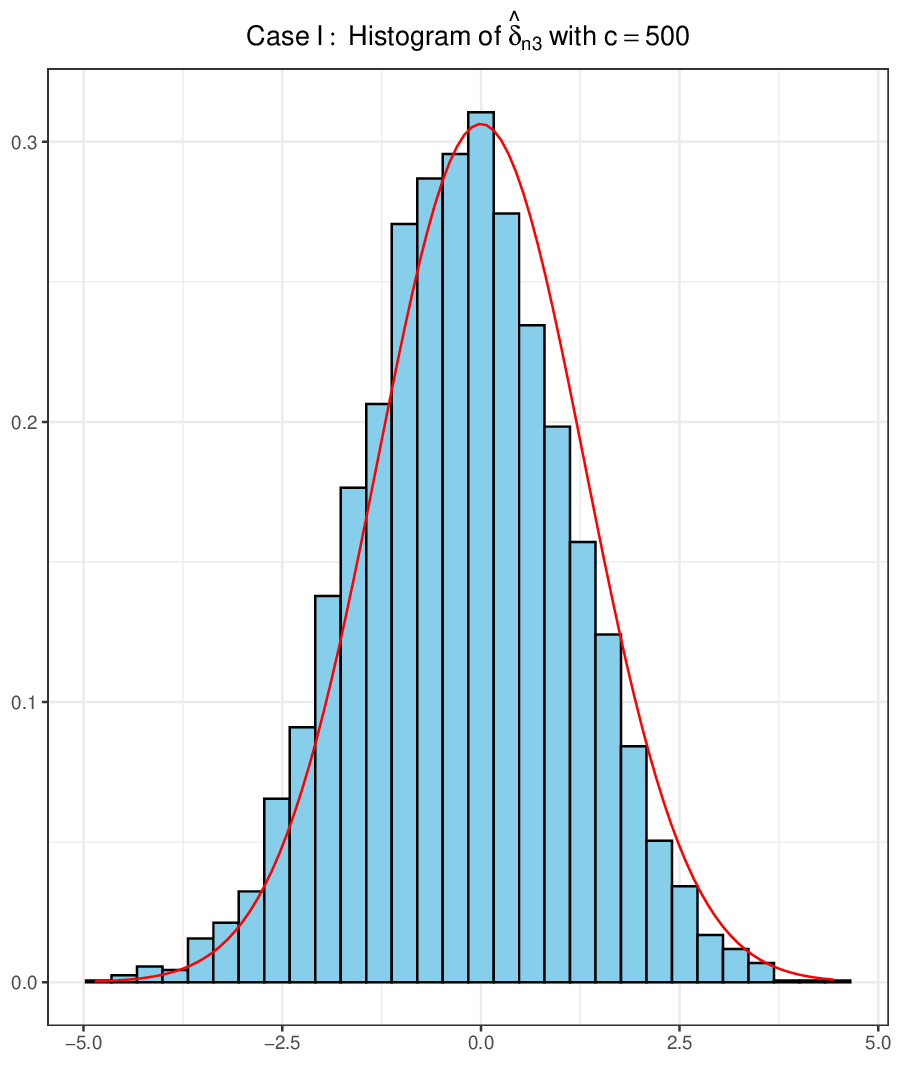}
    	\end{subfigure}\\
    	\rotatebox{0}{\tiny{Case II for $\A_n$}}
    	\begin{subfigure}{.1\textwidth}
    		\centering
    		\includegraphics[width=\linewidth]{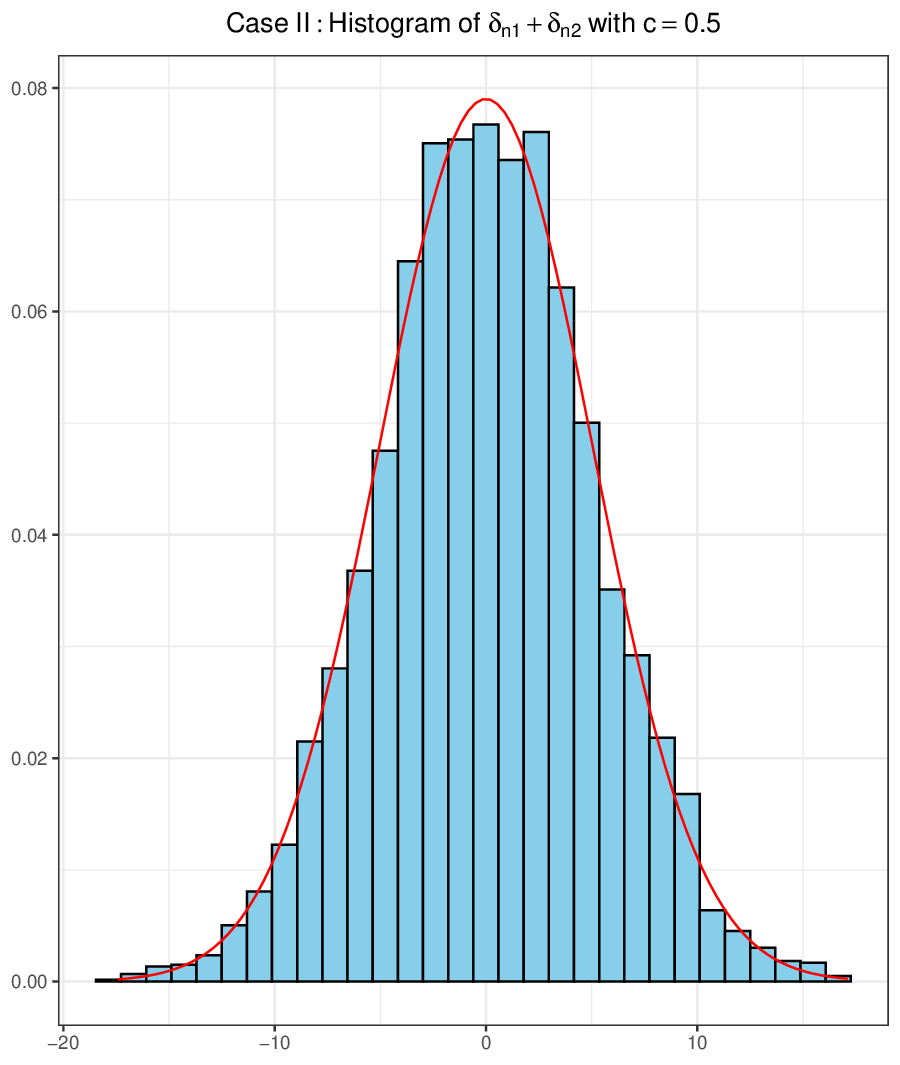}
    	\end{subfigure}
    	\begin{subfigure}{.1\textwidth}
    		\centering
    		\includegraphics[width=\linewidth]{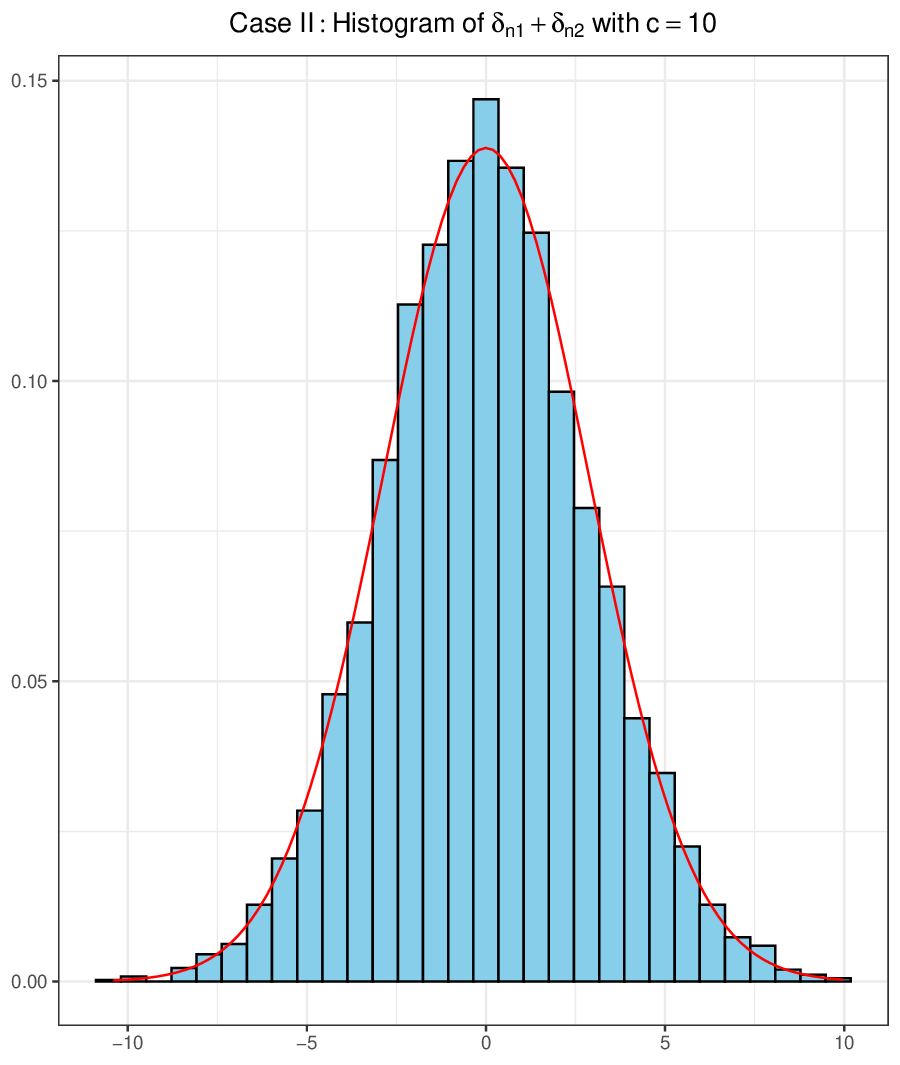}
    	\end{subfigure}
    	\begin{subfigure}{.1\textwidth}
    		\centering
    		\includegraphics[width=\linewidth]{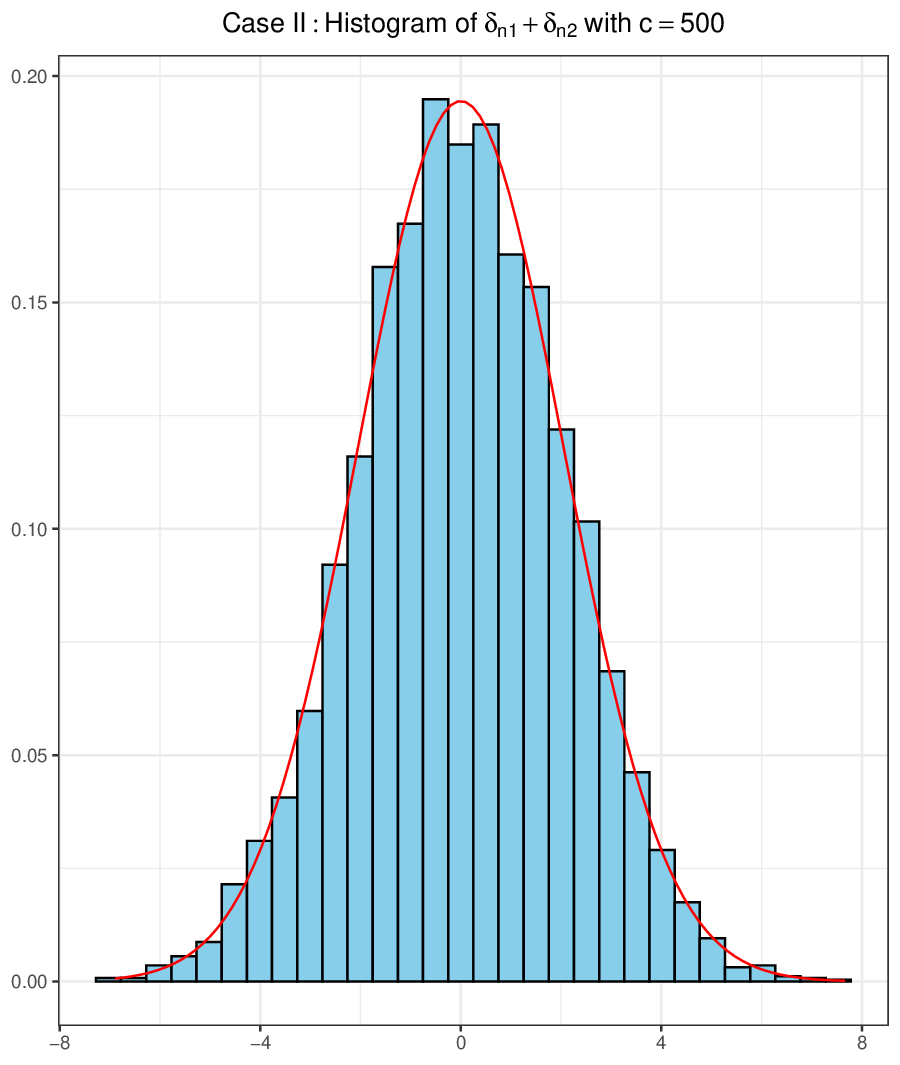}
    	\end{subfigure}\qquad
    	\begin{subfigure}{.1\textwidth}
    		\centering
    		\includegraphics[width=\linewidth]{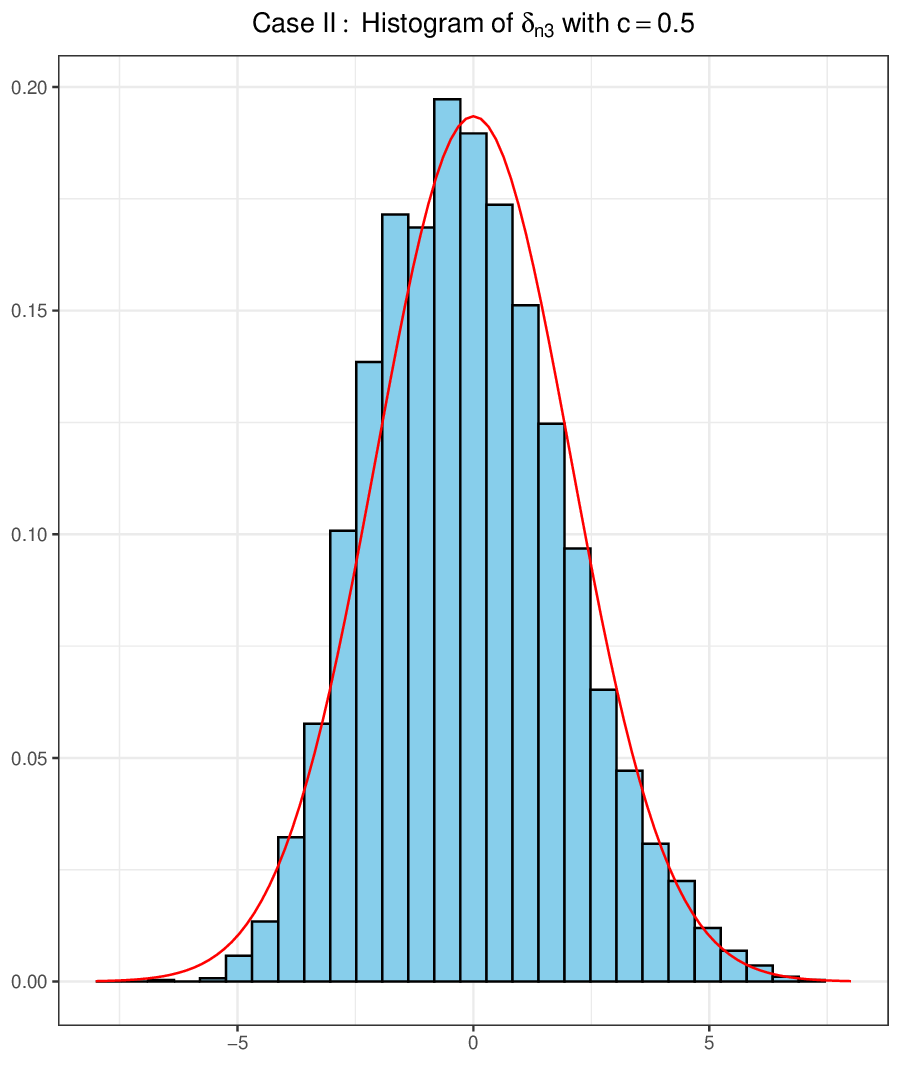}
    	\end{subfigure}
    	\begin{subfigure}{.1\textwidth}
    		\centering
    		\includegraphics[width=\linewidth]{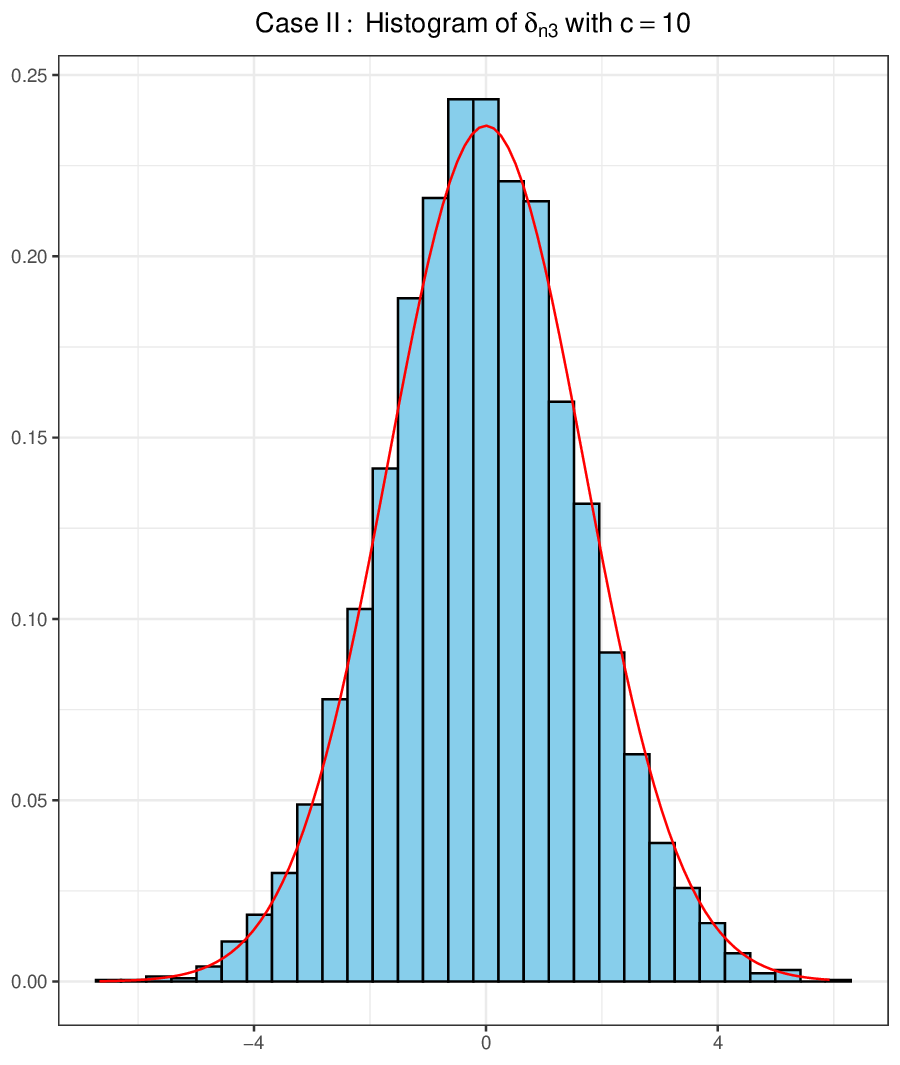}
    	\end{subfigure}
    	\begin{subfigure}{.1\textwidth}
    		\centering
    		\includegraphics[width=\linewidth]{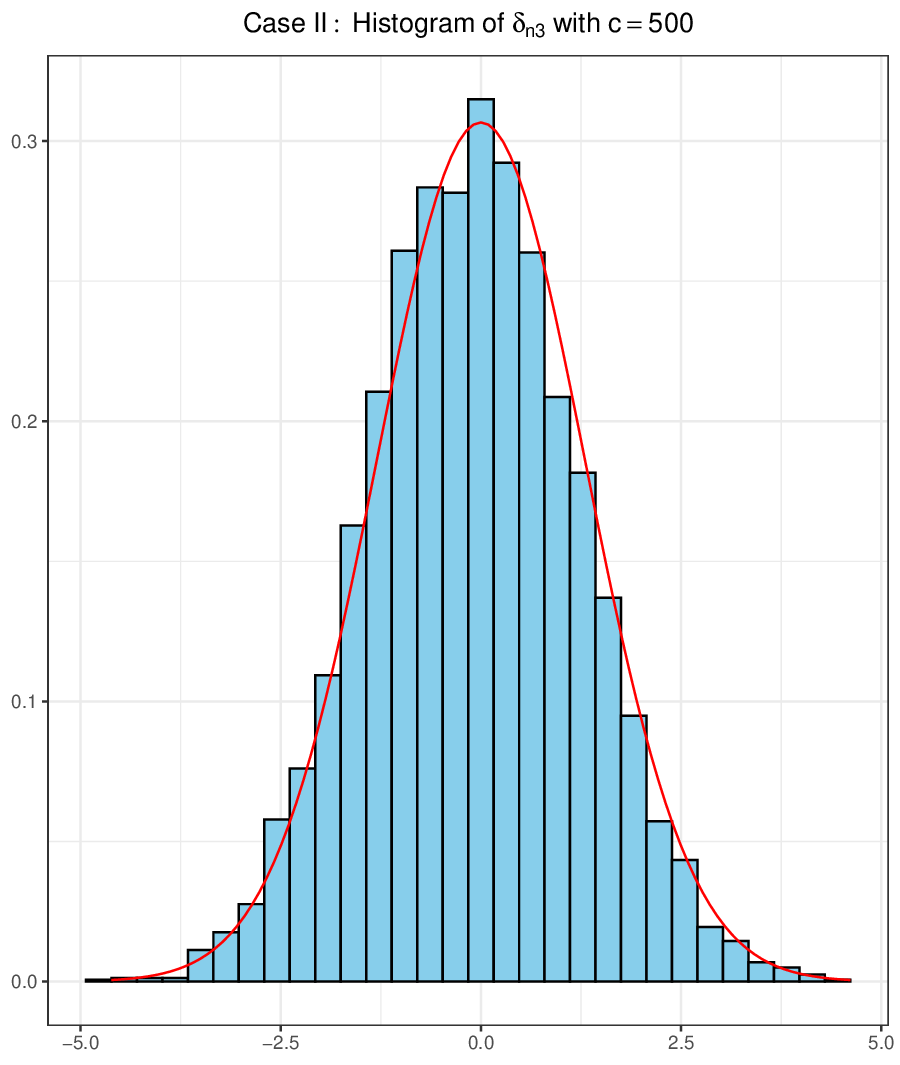}
    	\end{subfigure}\\
    	\rotatebox{0}{\tiny{Case II for $\hat\A_n$}}
    	\begin{subfigure}{.1\textwidth}
    		\centering
    		\includegraphics[width=\linewidth]{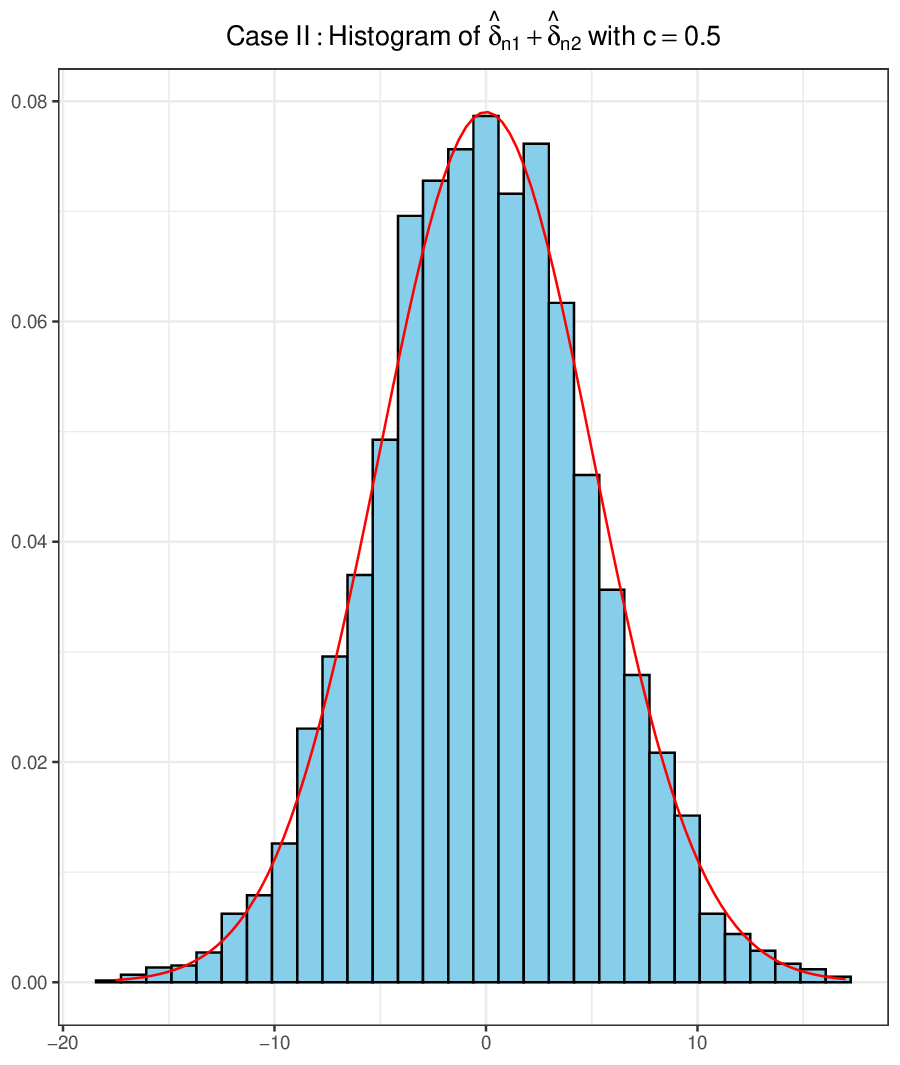}
    	\end{subfigure}
    	\begin{subfigure}{.1\textwidth}
    		\centering
    		\includegraphics[width=\linewidth]{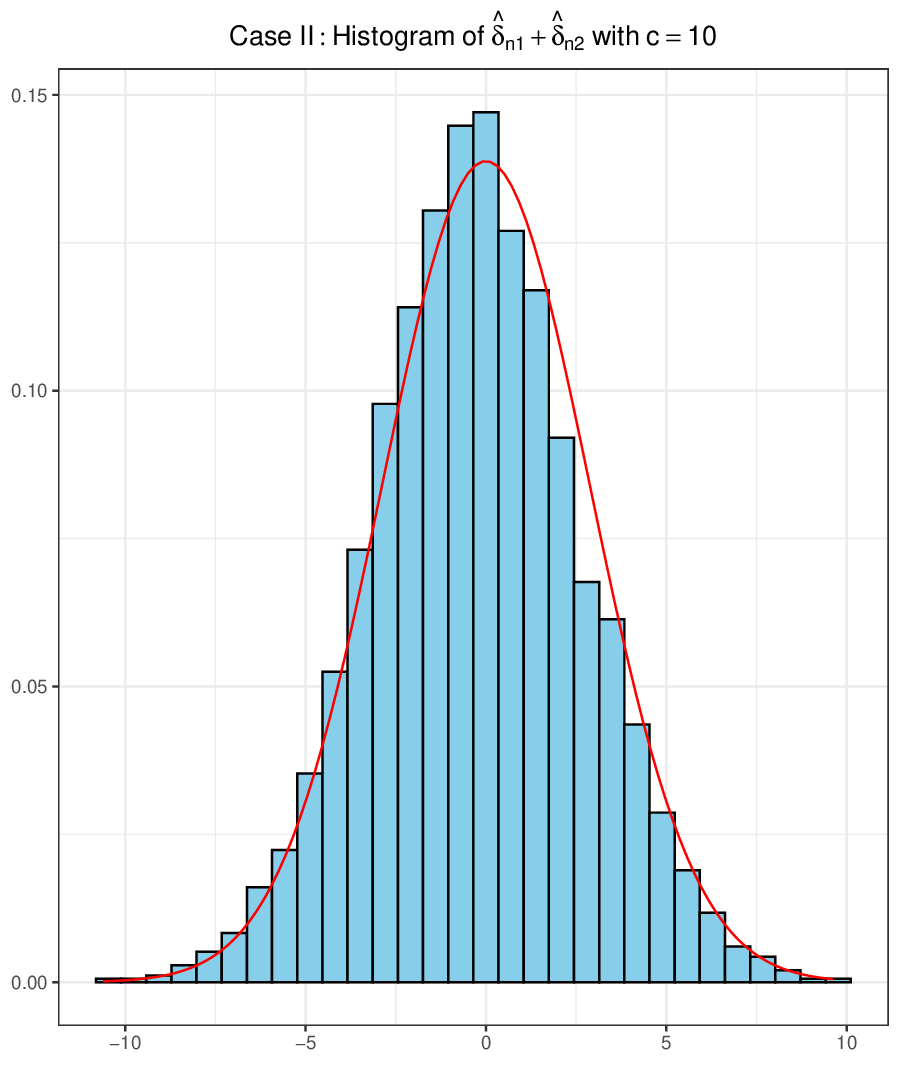}
    	\end{subfigure}
    	\begin{subfigure}{.1\textwidth}
    		\centering
    		\includegraphics[width=\linewidth]{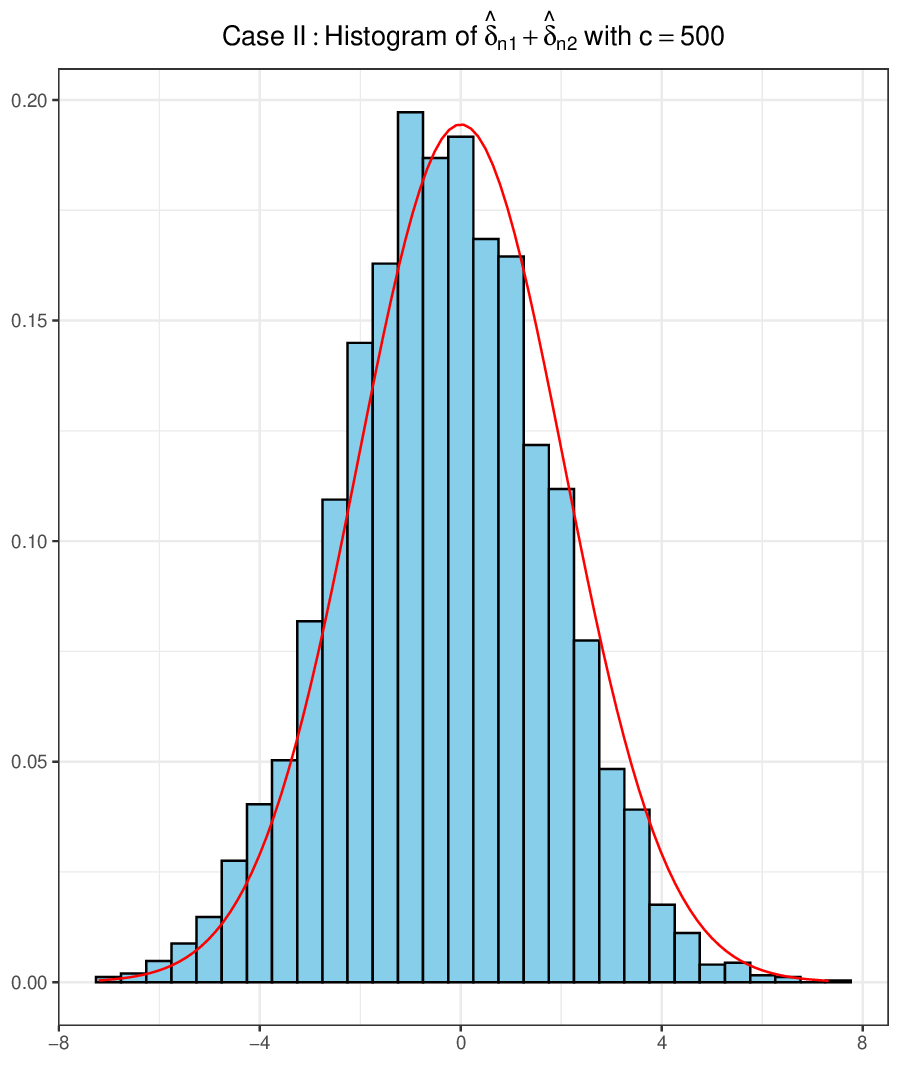}
    	\end{subfigure}\qquad
    	\begin{subfigure}{.1\textwidth}
    		\centering
    		\includegraphics[width=\linewidth]{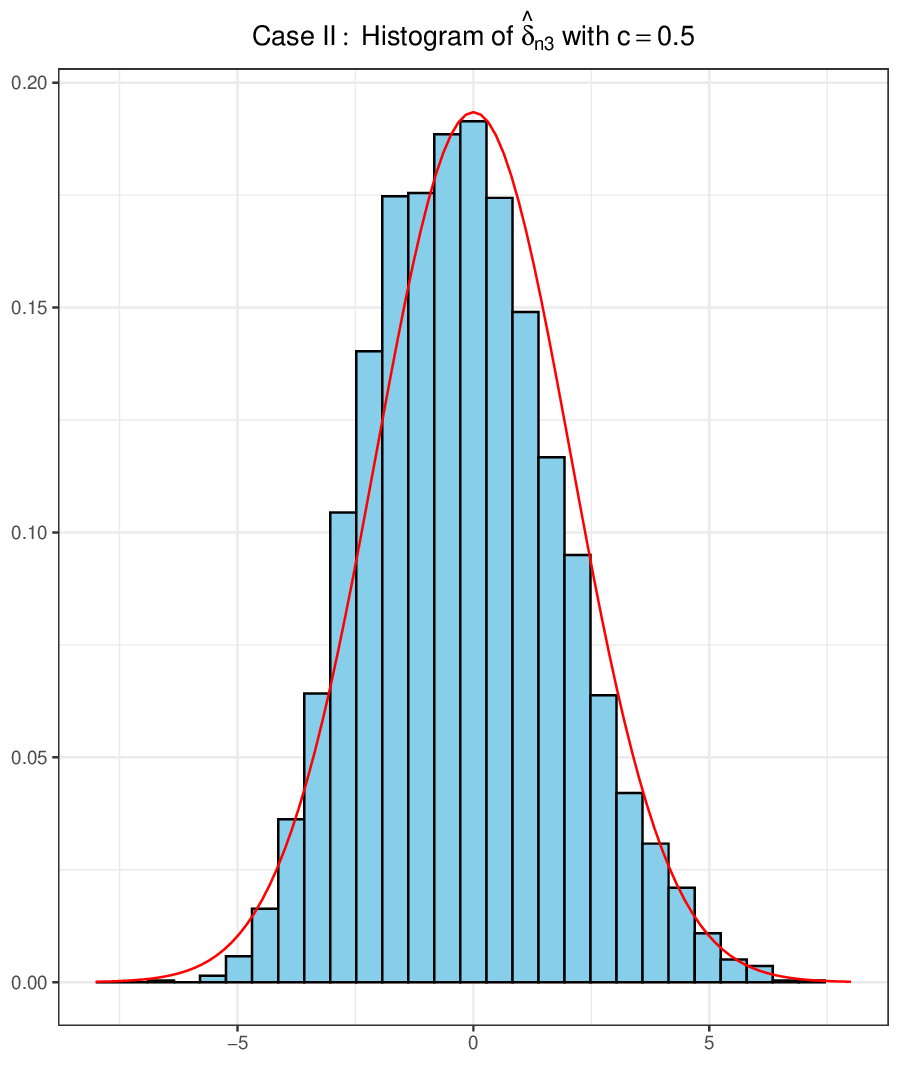}
    	\end{subfigure}
    	\begin{subfigure}{.1\textwidth}
    		\centering
    		\includegraphics[width=\linewidth]{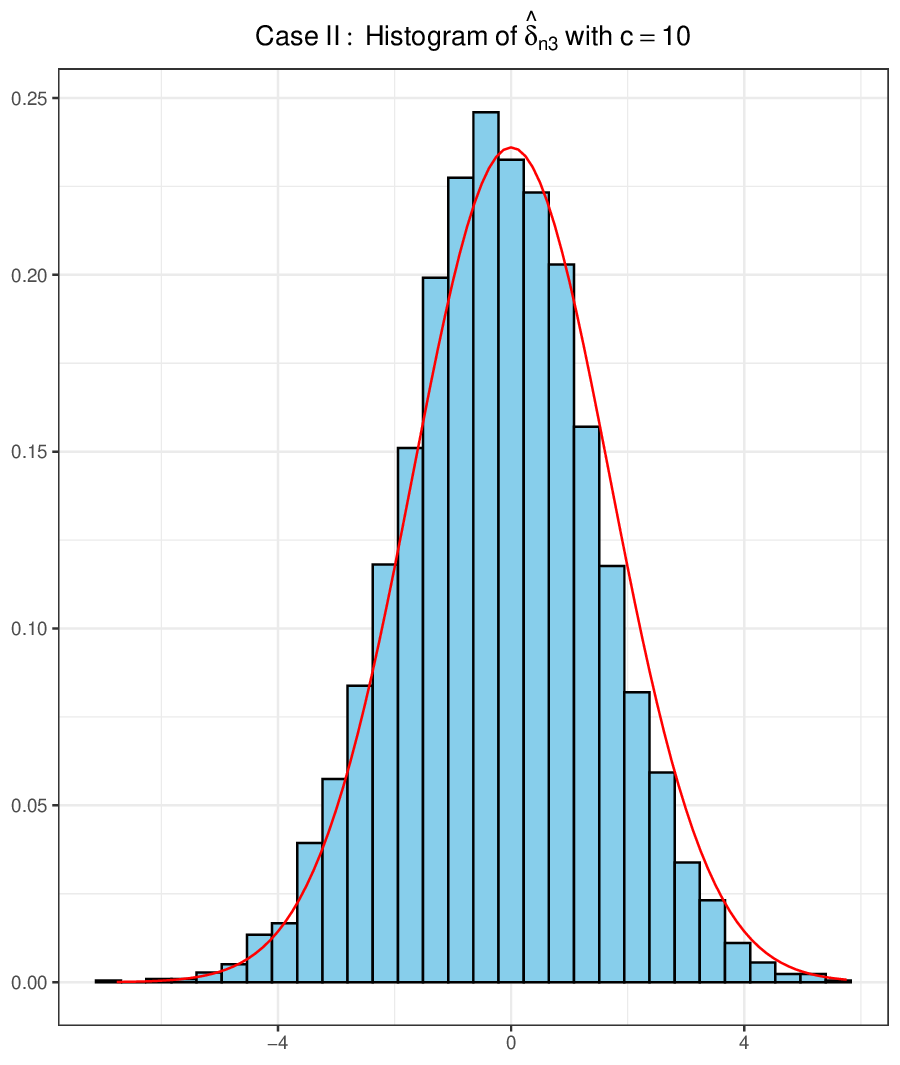}
    	\end{subfigure}
    	\begin{subfigure}{.1\textwidth}
    		\centering
    		\includegraphics[width=\linewidth]{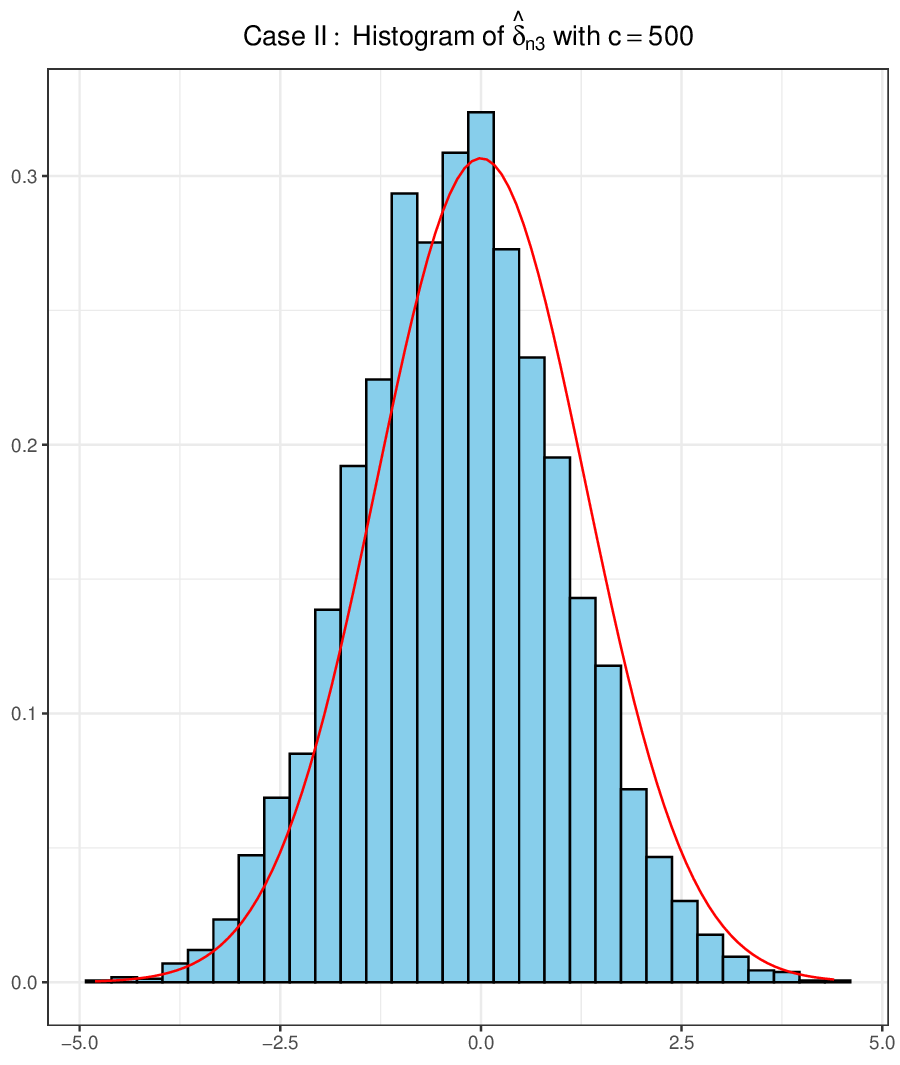}
    	\end{subfigure}\\
    	\qquad\qquad\qquad
    	\begin{subfigure}{.1\textwidth}
    		\centering
    		\caption*{\tiny{$c=0.5$\\$\delta_{n1}+\delta_{n2}$\\($\hat\delta_{n1}+\hat\delta_{n2}$)}}
    	\end{subfigure}
    	\begin{subfigure}{.1\textwidth}
    		\centering
    		\caption*{\tiny{$c=10$\\$\delta_{n1}+\delta_{n2}$\\($\hat\delta_{n1}+\hat\delta_{n2}$)}}
    	\end{subfigure}
    	\begin{subfigure}{.1\textwidth}
    		\centering
    		\caption*{\tiny{$c=500$\\$\delta_{n1}+\delta_{n2}$\\($\hat\delta_{n1}+\hat\delta_{n2}$)}}
    	\end{subfigure}\qquad
    	\begin{subfigure}{.1\textwidth}
    		\centering
    		\caption*{\tiny{$c=0.5$\\$\delta_{n3}$\\($\hat\delta_{n3}$)}}
    	\end{subfigure}
    	\begin{subfigure}{.1\textwidth}
    		\centering
    		\caption*{\tiny{$c=10$\\$\delta_{n3}$\\($\hat\delta_{n3}$)}}
    	\end{subfigure}
    	\begin{subfigure}{.1\textwidth}
    		\centering
    		\caption*{\tiny{$c=500$\\$\delta_{n3}$\\($\hat\delta_{n3}$)}}
    	\end{subfigure}
    	\caption{
    		Histograms for \(\delta_{n1} + \delta_{n2}\), \(\delta_{n3}\), \(\hat{\delta}_{n1} + \hat{\delta}_{n2}\), and \(\hat{\delta}_{n3}\) under Cases I and II from 5000 independent replications. The first row shows the histograms of \(\delta_{n1} + \delta_{n2}\) and \(\delta_{n3}\) under Case I, fitted by their Gaussian limits (red solid curves), with \(c \in\{0.5, 10, 500\}\). The second row is for \(\hat{\delta}_{n1} + \hat{\delta}_{n2}\) and \(\hat{\delta}_{n3}\) under Case I. Case II are shown in the third and fourth rows.
    	}
    	\label{fig-clt}
    \end{figure}

   Detailed limiting distributions of \(\delta_{n1} + \delta_{n2}\) and \(\delta_{n3}\) are listed in Table~\ref{table1}. Empirical histograms based on 5000 independent replications are shown in Figures \ref{fig-clt}, matching the theoretical results in Table \ref{table1}.

    \begin{table}[h!]
    	\caption{Limiting distributions of $\delta_{n1}+\delta_{n2}$ and $\delta_{n3}$ under Cases I and II with $c=0.5,\ 10$ and $500$.}\label{table1}
    	\centering
    	\begin{tabular}{cccc}
    		\hline
    		&c=0.5&\quad c=10&\quad c=500\\
    		\hline
    		\multirow{2}*{Case I}&$\delta_{n1}+\delta_{n2}\xrightarrow{d}N\left(0,  31.5876\right)$&$\delta_{n1}+\delta_{n2}\xrightarrow{d}N\left(0,  8.6294\right)$&
    		$\delta_{n1}+\delta_{n2}\xrightarrow{d}N\left(0,  4.2157\right)$\\
    		&$\delta_{n3}\xrightarrow{d} N\left(0,4.6717\right)$&
    		$\delta_{n3}\xrightarrow{d} N\left(0, 2.9699\right)$&
    		$\delta_{n3}\xrightarrow{d} N\left(0, 1.6951\right)$\\
    		\hline	
    		\multirow{2}*{Case II}&$\delta_{n1}+\delta_{n2}\xrightarrow{d}N\left(0,  25.4846\right)$&$\delta_{n1}+\delta_{n2}\xrightarrow{d}N\left(0,  8.2612\right)$&
    		$\delta_{n1}+\delta_{n2}\xrightarrow{d}N\left(0,  4.2081\right)$\\
    		&$\delta_{n3}\xrightarrow{d} N\left(0,4.2526\right)$&
    		$\delta_{n3}\xrightarrow{d} N\left(0,  2.8572\right)$&
    		$\delta_{n3}\xrightarrow{d} N\left(0,  1.6924\right)$\\
    		\hline	
    	\end{tabular}	
    \end{table}

    \section{Proofs}\label{sec:5}
This section presents the proofs of Theorems \ref{th-limit-cinf}, \ref{th-c-cinf}, \ref{th-limit} and  \ref{th-c}, along with a sketch of the proofs of Theorems \ref{th-as-cinf}-\ref{th-clt-cinf} and \ref{th-as}-\ref{th-clt}.
Notations are listed below and will be used throughout this section.
	\begin{align*}
	\delta_{n \ell}&=\sqrt{n}\left(\lambda_{k\ell}^{\mathbf{A}_n}-\lambda_{n k}\right), \ell \in \mathcal{L}_k,\quad
	\mathcal{L}_k=[1:m_k],\quad 
	\lambda_{n k}=\phi_{c_n, H_p}\left(\alpha_{n k}\right),\\ \nonumber
	\hat{\mathbf{U}}_n&=\left(\sqrt{k_{n 1}} \bar{\mathbf{x}}_1, \ldots, \sqrt{k_{n \tau}} \bar{\mathbf{x}}_\tau\right)_{p \times \tau},	\tilde{\lambda}_{k\ell}^{\mathbf{A}_n}=
	c_na_p+\sqrt{c_nb_p}\lambda_{k\ell}^{\mathbf{A}_n},\quad
	\tilde{\lambda}_{nk}=
	c_na_p+\sqrt{c_bb_p}\lambda_{nk},
	\\\nonumber
	\mathbf{R}_n&=\sqrt{n}\left[
	\left(\sqrt{\frac{p}{n b_p}} a_p+\lambda_{nk}\right) \hat{\mathbf{U}}_n^{\top}\left(\mathbf{B}_n-\tilde{\lambda}_{nk} \mathbf{I}\right)^{-1} \hat{\mathbf{U}}_n+\frac{1}{s_0(\lambda_{nk})}+\lambda_{nk}+\sqrt{\frac{p}{n b_p}} a_p\right.\\\nonumber
	&\qquad\qquad\left.+\frac{1}{s_0(\lambda_{nk})}\bU_n^{\top}\left(\bSigma_0+\sqrt{\frac{pb_p}{n }}\frac{1}{s_0(\lambda_{nk})}\bI\right)^{-1}\bU_n
	\right].
	\end{align*}
We denote by \(\mathfrak{M}\) some constants that appear in inequalities and may take different values at different appearances. The orders \(o(\cdot)\) and \(O(\cdot)\) for vectors are in terms of the Euclidean norm, and for matrices, they are in terms of the spectral norm.

	\subsection{Proofs of Theorems \ref{th-limit-cinf} and \ref{th-limit}}\label{step-as}
  
For the case where \(p/n\rightarrow c \in (0, \infty)\), results in Theorem~\ref{th-limit} can be derived from  ``exact separation" for the eigenvalues of \(\mathbf{S}_n\) in Chapter 3 of \cite{yiminghigh} and Theorems 2-3 in \cite{liu2023asymptotic}. We thus focus on the case where \(p/n\rightarrow \infty\).

Notice that Lemma \ref{th-limit1-cinf} and Lemma \ref{th-as-cinf} are established in the almost sure sense. We can thus consider a sample realization, denoted as \(\{\bx_i: i \geq 1\}\), such that the convergence in these two lemmas holds.  
In addition, from 
	$\left(\boldsymbol{\mu}_i-\boldsymbol{\mu}_j\right)^{\top}\left(\boldsymbol{\mu}_i-\boldsymbol{\mu}_j\right) \asymp \sqrt{c_n}$ and $\bar\s_i^{\top}\left(\boldsymbol{\mu}_i-\boldsymbol{\mu}_j\right)/\sqrt{c_n}\to0$, we have $\lambda_{1}^{\A_n}\leq \mathfrak{M}$. 
Then, Theorem \ref{th-limit-cinf} can be proved by demonstrating the following two claims for this particular realization.
	\begin{itemize}
		\item[]{\bf Claim 1.}
If 
\begin{align}\label{inf}
\min_{\ell \in \mathcal{L}_k}\liminf_{n\to\infty}\lambda_{k\ell}^{\mathbf{A}_n}>2,
\end{align} 
then \(\{\lambda_{k\ell}^{\mathbf{A}_n}, \ell \in \mathcal{L}_k\}\) converge to a common limit which is larger than 2. In addition, it holds that $\alpha_k> 1$.
		 \item[]{\bf Claim 2.}
If \(\alpha_k \leq 1\), then \(\{\lambda_{k\ell}^{\mathbf{A}_n}, \ell \in \mathcal{L}_k\}\) converge to $\mathfrak{b}=2$, the right edge point of the support of the semicircle law.
		 \end{itemize}
		 
		 {\bf Proof for Claim 1.} From the condition \eqref{inf}, \(\{\lambda_{k\ell}^{\mathbf{A}_n}: \ell \in \mathcal{L}_k\}\) are eigenvalues of \(\mathbf{A}_n\), while \(\{\tilde{\lambda}_{k\ell}^{\mathbf{A}_n}:  \ell \in \mathcal{L}_k\}\) are not eigenvalues of \(\mathbf{B}_n\) for all large \(p\) and \(n\). Therefore, we obtain
		\begin{align}\label{tau-eq}
		  0&=\left|\mathbf{A}_n-\lambda_{k\ell}^{\mathbf{A}_n} \mathbf{I}_n\right|=\left|\sqrt{\frac{p}{n b_p}} \frac{1}{p} \boldsymbol{\Phi}\mathbf{X}^{\top}\mathbf{X} \mathbf{\Phi}-\left(\sqrt{\frac{p}{n b_p}} a_p+\lambda_{k\ell}^{\mathbf{A}_n}\right) \mathbf{I}_n\right|\nonumber\\
		  &=\left|\mathbf{B}_n+\hat{\mathbf{U}}_n \mathbf{N}_n \hat{\mathbf{U}}_n^{\top}-\tilde{\lambda}_{k\ell}^{\mathbf{A}_n} \mathbf{I}_p\right| 
		  =\left|\mathbf{B}_n-\tilde{\lambda}_{k\ell}^{\mathbf{A}_n} \mathbf{I}_p\right|\left|\mathbf{I}_p+\left(\mathbf{B}_n-\tilde{\lambda}_{k\ell}^{\mathbf{A}_n} \mathbf{I}_p\right)^{-1} \hat{\mathbf{U}}_n \mathbf{N}_n \hat{\mathbf{U}}_n^{\top}\right| \nonumber\\
		  &=\left|\mathbf{I}_\tau+\mathbf{N}_n \hat{\mathbf{U}}_n^{\top}\left(\mathbf{B}_n-\tilde{\lambda}_{k\ell}^{\mathbf{A}_n} \mathbf{I}_p\right)^{-1} \hat{\mathbf{U}}_n \mathbf{N}_n\right|.
		  \end{align}
For any fixed \(k\) and \(\ell\), let \(\{\lambda_{k\ell}^{\A_{n_i}}\}\) be a subsequence of \(\{\lambda_{k\ell}^{\A_n}\}\) that converges to a limit, say \(\lambda_k\). From Lemma \ref{th-limit1} and Lemma \ref{th-as} for \(p/n\rightarrow\infty\), the matrix in \eqref{tau-eq} converges entry-wise such that 
\begin{align}\label{la-a}
0=\left|\frac{1}{s(\lambda_k)}\mathbf{I}_\tau+\frac{1}{\sqrt{c_nb_p}}\bN_n\bU_n^{\top}\bU_n\bN_n+o(1)\right|.
\end{align}
Then, taking the limit of \eqref{la-a} as \(n \to \infty\) and using the identity \( z=-1/s(z)-s(z)\), we obtain
$$
\alpha_k = -\frac{1}{s(\lambda_k)}=\frac{2}{\lambda_k-\sqrt{\lambda_k^2-4}} \quad \text{and}\quad
\lambda_k=\alpha_k+\frac1{\alpha_k}>2.
$$
In addition, since \(\lambda_k > 2\), we have \(\alpha_k\) strictly greater than 1. The first claim is thus verified.
	
	{\bf Proof for Claim 2.} Recall the matrix \(\mathbf{Z}\) defined in Assumption \ref{as-moment-cinf} and let 
	\begin{align}\label{par-A0}
	\mathbf{A}_0=\frac{1}{\sqrt{n p b_p}}\left(\mathbf{Z}^{\top} \boldsymbol{\Sigma}_0 \mathbf{Z}-p a_p \mathbf{I}_n\right).
	\end{align}
For the projection matrix \(\mathbf{\Phi} = \mathbf{I}_n - \mathbf{1}_n \mathbf{1}_n^{\top}/n\), by Poincar\'{e} separation theorem (Corollary 4.3.37 in \cite{horn2012matrix})
	we have $$\lambda_{2\tau}^{\A_0}\leq\lambda_{2\tau-1}^{\mathbf{\Phi}\A_0\mathbf{\Phi}}\leq \lambda_{1}^{\mathbf{\Phi}\A_0\mathbf{\Phi}}\leq \lambda_{1}^{\A_0}.$$
	Moreover,
 $\operatorname{rank}(\A_{n}-\mathbf{\Phi}\A_0\mathbf{\Phi})\leq \tau-1.
	$ Thus, by using Weyl's interlacing inequalities, we have
\begin{align}\label{eig-inq-pa}
	\lambda_{2\tau}^{\A_0}\leq\lambda_{2\tau-1}^{\mathbf{\Phi}\A_0\mathbf{\Phi}}\leq\lambda_{\tau}^{\A_n}\leq \lambda_{1}^{\mathbf{\Phi}\A_0\mathbf{\Phi}}\leq \lambda_{1}^{\A_0}.
	\end{align}
Moreover, from Theorem 3.13 in \cite{jingbingyi}, both \(\lambda_{2\tau}^{\A_0}\) and \(\lambda_{1}^{\A_0}\) converge to 2, almost surely, as \(n\to\infty\), \(p/n \rightarrow \infty\).  
Therefore, without loss of generality, we can assume that \(\lambda_{\tau}^{\A_n}\in [\lambda_{2\tau}^{\A_0}, \lambda_{1}^{\A_0}] \to 2\) for the particular realization \(\{\bx_i: i \geq 1\}\). This implies that for any convergent subsequence \(\{\lambda_{k\ell}^{\A_{n_i}}\}\) of \(\{\lambda_{k\ell}^{\A_n}\}\), its limit is either 2 or some constant greater than 2. However, according to \textbf{Claim 1}, if this limit is greater than 2, we have \(\alpha_k > 1\), which leads to a contradiction. Hence, the limit must be 2, which verifies the second claim. 

The proofs of Theorems \ref{th-limit-cinf} and \ref{th-limit} are thus complete.

 \subsection{Proofs of Theorems \ref{th-c-cinf} and \ref{th-c}}
The proof of Theorem \ref{th-c-cinf} is embedded in the proof of Theorem \ref{th-c}. Next, we show the proof of Theorem \ref{th-c}.
For distant spiked eigenvalues \(\{\lambda_{k\ell}^{\A_n}, \ell \in \mathcal{L}_k\}\), 
our strategy is to investigate the limit of the \(\tau\)-order determinant in \eqref{tau-eq} based on the asymptotic properties of random sesquilinear forms introduced in Theorems \ref{th-as} and \ref{th-clt}. The proof can be accomplished through four steps.
	\begin{itemize}
		\item[]{\bf Step 1.} By simplifying equation \eqref{tau-eq}, we link the normalized eigenvalues \(\left\{\delta_{n \ell}, \ell \in \mathcal{L}_k\right\}\) with a random matrix \(\bR_n\); see Lemma \ref{le-tau-equation}.
		\begin{lemma}\label{le-tau-equation}
			Under the assumptions of  Theorem \ref{th-c}, the determinant equation in \eqref{tau-eq} can be approximated as
			\begin{align}\label{con:huajianjieguo}
			0&=\left|\left(1-\frac{\delta_{n \ell}}{\sqrt{n}}\frac{s_0\left(\lambda_{nk}\right)+s_0^{\prime}\left(\lambda_{nk}\right)\left(
				\lambda_{nk}+\sqrt{\frac{p}{nb_p}}a_p\right)}
			{s_0\left(\lambda_{nk}\right)\left(
				\lambda_{nk}+\sqrt{\frac{p}{nb_p}}a_p\right)}\right)
			\right.\nonumber\\
			&\phantom{=\;\;}\left.\quad
			\times
			\left(\bI_{\tau}+\bN_n\bU_n^{\top}\left(\bSigma_0-\sqrt{\frac{pb_p}{n}}\alpha_{nk} \bI\right)^{-1} \bU_n \bN_n\right)
			- s_0\left(\lambda_{nk}\right)\frac{1}{\sqrt{n}}\bN_n\bbR_n\bN_n
			\right.\nonumber\\
			&\phantom{=\;\;}\left.\quad+\frac{\delta_{n\ell}}{\sqrt{n}} \frac{s_0^{\prime}\left(\lambda_{nk}\right)}{s_0^{2}\left(\lambda_{nk}\right)}\sqrt{\frac{pb_p}{n}} \bN_n\bU_n^{\top}\left(\bSigma_0-\sqrt{\frac{pb_p}{n}}\alpha_{nk} \bI\right)^{-2} \bU_n \bN_n+o_p\left(\frac{1}{\sqrt{n}}\right)\right|.
			\end{align}
		\end{lemma}
		\item[]{\bf Step 2.} We derive the weak limit \(\bR\) of the matrix \(\bR_n\), as presented in the following lemma.
		\begin{lemma}\label{le-Rn}
			Under the assumptions of Theorem \ref{th-c}, the \(\tau \times \tau\) random matrix \(\bbR_n\) converges weakly to a symmetric zero-mean Gaussian matrix \(\bbR = \left(R_{ij}\right)\). The covariances of the entries  are
			\begin{align*}
			\operatorname{Cov}(R_{ij},R_{lt})&=\frac{\alpha_k^2}{\phi_{c,H}^{\prime}(\alpha_k)}	\operatorname{Cov}(W_{ij},W_{lt})\\
			&\quad-\frac{v_3}{s^2\left(\lambda_k\right)}\left(\sqrt{k_i} h_{j l t}+\sqrt{k_j} h_{i l t}+\sqrt{k_l} h_{i j t}+\sqrt{k_t} h_{i j l}\right)
			\end{align*}
 for $i,j,l,t\in\{1,\ldots,\tau\}$,
where \(\operatorname{Cov}(W_{ij}, W_{lt})\) is given in Theorem \ref{th-c} and \(h_{ijl}= h_{ijl}\left(\alpha_k\right)\) is defined in Assumption \ref{as-limits}.
		\end{lemma}
		
		\item[]{\bf Step 3.} By the Skorokhod strong representation theorem, the convergence \(\bbR_{n} \rightarrow \bbR\) and \eqref{con:huajianjieguo} take place almost surely on an appropriate probability space. Thus, we can take the limit of the RHS of \eqref{con:huajianjieguo}, which yields the limit of \(\{\delta_{n\ell}, \ell \in \mathcal{L}_k\}\); see Lemma \ref{le-mk-equation}.
		\begin{lemma}\label{le-mk-equation}
Under the assumptions of Theorem \ref{th-c}, each \(\delta_{n\ell}\) in \(\{\delta_{n\ell}, \ell \in \mathcal{L}_k\}\) converges to a limit \(\delta\), which solves the determinant equation
			\begin{equation}\label{eq-mk}
			0=\left|\delta \bI-\frac{s^{3}\left(\lambda_{k}\right)}{s^{\prime}\left(\lambda_{k}\right)}\bQ_k^{\top} \bN \bbR \bN \bQ_k \mathbf{G}\right|.
			\end{equation}
		\end{lemma}
		
	\item[]{\bf Step 4.}	We replace $\bbR$ with
	$\alpha_k\mathbf{W}/\sqrt{\phi_{c, H}^{\prime}\left(\alpha_k\right)}$ in equation \eqref{eq-mk}; see Lemma \ref{le-v3}.
	\begin{lemma}\label{le-v3}
  The $\tau\times\tau$ random matrix $\sqrt{\phi_{c, H}^{\prime}\left(\alpha_k\right)}\bN\bbR\bN/\alpha_k \stackrel{d}{=} \bN\mathbf{W}\bN$.
	\end{lemma}

Therefore, by the strong representation theorem, we obtain the convergence of the random vector \(\{\delta_{n\ell}, \ell \in \mathcal{L}_{k}\}\).
This, together with the identities
	\begin{align}\label{relations}
	-\frac{1}{s\left(\lambda_{k}\right)}=\alpha_{k},\quad \lambda_{k}=\phi_{c,H}(\alpha_{k}),\quad s^{\prime}\left(\lambda_{k}\right) =\frac{1}{\alpha_{k}^2\phi_{c,H}^{\prime}(\alpha_{k})},
	\end{align}
	gives the conclusion of Theorem \ref{th-c}. 
	\end{itemize}
The proofs of Lemmas \ref{le-tau-equation}-\ref{le-v3} are presented in the Supplementary Material, while the proofs of Lemmas \ref{le-tau-equation}-\ref{le-Rn} rely on Theorems \ref{th-as}-\ref{th-clt}.

\subsection{A sketch of  the proofs of Theorems  \ref{th-as-cinf}-\ref{th-clt-cinf} and \ref{th-as}-\ref{th-clt}}
The proofs of Theorems  \ref{th-as-cinf}-\ref{th-clt-cinf} are contained in the proofs of Theorems \ref{th-as}-\ref{th-clt}. 
In this section, we outline the main steps for proving Theorems \ref{th-as}-\ref{th-clt}. Detailed proofs are presented in the Supplementary Material.
	
\begin{itemize}
	\item[]{\bf Step 1.} 
		By the Woodbury matrix identity, one has
	$$
	\frac{\tilde{z}_n}{\sqrt{c_n b_p}}\left(\mathbf{B}_{n}-\tilde{z}_n \mathbf{I}\right)^{-1}=
	\bH-\bH\mathbf{M}_{\bar{\mathbf{s}}}\left[\frac{\tilde{z}_n}{\sqrt{c_n b_p}}\mathbf{K}_n^{-1}+\mathbf{M}_{\bar{\mathbf{s}}}^{
		\top}\bH\mathbf{M}_{\bar{\mathbf{s}}}\right]^{-1}\bH\mathbf{M}_{\bar{\mathbf{s}}}^{
		\top},
	$$
	where 	$$\mathbf{B}_{n 0}=\frac{1}{n} \sum_{i=1}^\tau \sum_{j=1}^{n_i} \mathbf{s}_{i j} \mathbf{s}_{i j}^{\top},\quad
	\bH=\frac{\tilde{z}_n}{\sqrt{c_n b_p}}\left(\mathbf{B}_{n0}-\tilde{z}_n \mathbf{I}\right)^{-1}.
	$$
Thus, to prove Theorem \ref{th-as}, it is sufficient to prove the following lemma.
	\begin{lemma}\label{par-as}
		Under the assumptions of Theorem \ref{th-as}, for any $i,j\in[1:\tau]$, we have
	\begin{align}
	&\frac{\sqrt{c_n}}{\tilde z_n}\bar\s_i^{\top}\bH\bmu_{j}\xrightarrow{\text { a.s. }} 0,\quad
	\bmu_{i}^{\top}\bH\bmu_{j}+\bmu_{i}^{\top}\left(\sqrt{c_n b_p} \mathbf{I}+s_0(z) \boldsymbol{\Sigma}_0\right)^{-1}\bmu_{j}\xrightarrow{\text { a.s. }} 0,\nonumber\\
	&\frac{b_p}{\tilde z_n}\bar{\mathbf{s}}_i^{\top} \bH \bar{\mathbf{s}}_j -\frac{1}{k_i}\left[a_ps(z)-\sqrt{b_p/c_n}s^2(z)\right]I(i=j)\xrightarrow{\text { a.s. }} 0.
	\nonumber
	\end{align}
	\end{lemma}
	\item[]{\bf Step 2.} We simplify  \(\bL_n\) in Theorem \ref{th-clt} to the form shown in the following lemma.
	\begin{lemma}\label{par-simp}
		Under the assumptions of Theorem \ref{th-clt}, we have
		\begin{align*}
	\bL_{n}=\sqrt{n} \left(\begin{array}{cc}
\bL_{11}& \bL_{12} \\ 
\bL_{21}& \bL_{22}
	\end{array}\right) +o_p(1),
	\end{align*}
	where
	\begin{align*}
\bL_{11}&=	\frac{\sqrt{c_nb_p}}{\tilde{z}_ns^2(z)}
	\left[
	\mathbf{M}_{\bar{\mathbf{s}}}^{\top} \left(\bB_{n0}-\tilde{z}_n\bI\right)^{-1} \mathbf{M}_{\bar{\mathbf{s}}}-\left(1+\frac{\tilde{z}_ns_0(z)}{\sqrt{c_nb_p}}\right)\bK_n^{-1}\right] ,\\
\bL_{12}&=-\left[s(z)\right]^{-1}\bM_{\bar\s}^{\top}\left(\bB_{n0}-\tilde{z}_n\bI\right)^{-1}\bM_{\bmu},\quad
\bL_{21}=	-\left[s(z)\right]^{-1}\bM_{\bmu}^{\top}\left(\bB_{n0}-\tilde{z}_n\bI\right)^{-1}\bM_{\bar\s},\\
\bL_{22}&=	\frac{\tilde{z}_n}{\sqrt{c_nb_p}}
	\bM_{\bmu}^{\top}\left(\bB_{n0}-\tilde{z}_n\bI\right)^{-1}\bM_{\bmu}+\bM_{\bmu}^{\top}\left(\sqrt{c_nb_p}\bI+s_0(z)\bSigma_0\right)^{-1}\bM_{\bmu}.
	\end{align*}
\end{lemma}
   \item[]{\bf Step 3.} Using the Cram\'er-Wold device,  Theorem \ref{th-clt} can be derived by demonstrating the following lemma.
   
   \begin{lemma}\label{le-cra}
     		Under the assumptions of Theorem \ref{th-clt}, for any real number $\{a_{qi j}: q\in[1:3], ~ i, j\in[1:\tau]\}$, we have
     \begin{align*}
     \sqrt{n}\sum_{1 \leq i,j \leq \tau} 
     &\left[a_{1ij}
     \frac{\sqrt{c_n b_p}}{\tilde{z}_n s^2(z)}
     \left\{\bar\s_i^\top\left(\mathbf{B}_{n 0}-\tilde{z}_n \mathbf{I}\right)^{-1}\bar\s_j-
     \frac{I(i=j)}{k_{ni}}\left(1+\frac{\tilde{z}_n s_0(z)}{\sqrt{c_n b_p}}\right)\right\}\right.\\ 
     &\phantom{=\;\;}\left.\
     +a_{2ij}\left\{\frac{\tilde{z}_n}{\sqrt{c_n b_p}}\bmu_i^{\top}\left(\mathbf{B}_{n 0}-\tilde{z}_n \mathbf{I}\right)^{-1}\bmu_j+ {\mathbf{\bmu}}_{i}^{\top} \left(\sqrt{c_n b_p} \mathbf{I}+s_0(z) \boldsymbol{\Sigma}_0\right)^{-1} {\mathbf{\bmu}}_{j}\right\}\right.\\
      &\phantom{=\;\;}\left.\
      -a_{3ij}[s(z)]^{-1}\bar\s_i^\top\left(\mathbf{B}_{n 0}-\tilde{z}_n \mathbf{I}\right)^{-1}\bmu_j\right]
     \xrightarrow{D} N\left(0, \sigma^{2}\right),
     \end{align*}
     where
     \begin{align*}
     \sigma^{2}=\sum_{1\leq i,j,l,t\leq \tau} a_{1i j}a_{1lt}\sigma_{1ij1lt}^2+a_{2i j}a_{2lt}\sigma_{2ij2lt}^2+a_{3i j}a_{3lt}\sigma_{3ij3lt}^2+2a_{2i j}a_{3lt}\sigma_{2ij3lt}^2,
     \end{align*}
     with
     $$ 
     \begin{aligned}
     \sigma_{1ij1lt}^2 &=\begin{cases}2 \frac{s^{\prime}(z)-s^2(z)}{s^4(z)} \frac{1}{k_i^2} & \text { if } i=j=l=t, \\ \frac{s^{\prime}(z)-s^2(z)}{s^4(z)} \frac{1}{k_i k_j} & \text { if } i=l \neq j=t,\\ 0 & \text { o.w. ;}\end{cases}\\
     \sigma_{2ij2lt}^2 &=\frac{s^{\prime}(z)}{s^4(z)}\left(\zeta_{j l} \zeta_{t i}+\zeta_{jt} \zeta_{li}\right)+\frac{v_4-3}{s^2(z)}g_{ijlt};\\
     \sigma_{3ij3lt}^2 &=\frac{\zeta_{jt}}{k_i}\frac{s^{\prime}(z)}{s^4(z)}I(i=l)
     \quad\text{and}\quad
     \sigma_{2ij3lt}^2=-\frac{v_3}{s^2(z)}f_{ijt}.
     \end{aligned}
     $$
   \end{lemma}
\end{itemize}

Lemma \ref{le-cra} is proved by using the Martingale CLT (Theorem 35.12 of \cite{billingsley2008probability}). We highlight some key points in the proof, especially for the multi-sample case in the ultrahigh dimensional context. 
Detailed proofs are presented in the Supplementary Material.
	\begin{itemize}
		\item 
Central tasks in the proof involve analyzing martingale differences of the following quantities:
\[
\left\{\bar{\mathbf{s}}_i^\top\left(\mathbf{B}_{n 0} - \tilde{z}_n \mathbf{I}\right)^{-1}\bar{\mathbf{s}}_j, ~\bar{\mathbf{s}}_i^\top\left(\mathbf{B}_{n 0} - \tilde{z}_n \mathbf{I}\right)^{-1}\bmu_j: i, j \in[1:\tau]\right\}.
\]
Notice that the resolvent matrix \((\mathbf{B}_{n 0} - \tilde{z}_n \mathbf{I})^{-1}\) incorporates all the \(n\) vectors \(\{\mathbf{s}_{ij}\}\), while each of the sample mean vectors \(\{\bar{\mathbf{s}}_i\}\) only involves a part of them. This would make the martingale decomposition more complicated and redundant compared to the single population case.
To address this, we employ a unified form for all \(\bar{\mathbf{s}}_i\)'s. Specifically, write 
\[
\left(\mathbf{s}_1, \ldots, \mathbf{s}_n\right) = \left(\mathbf{s}_{11}, \ldots, \mathbf{s}_{1 n_1}, \ldots, \mathbf{s}_{\tau 1}, \ldots, \mathbf{s}_{\tau n_\tau}\right).
\]
Then, \(\{\bar{\mathbf{s}}_i, i = 1, \ldots, \tau\}\) can be represented as weighted averages of all the vectors, i.e.,
\begin{align}\label{pa-weight}
\bar{\mathbf{s}}_i = \frac{1}{n_i} \sum_{j=1}^{n_i} \mathbf{s}_{ij} = \frac{1}{n} \sum_{j=1}^{n} \mathbf{e}_j^\top \mathbf{w}_i \mathbf{s}_j = \frac{1}{n} \sum_{j=1}^{n} w_{ij} \mathbf{s}_j,
\end{align}
where \(\mathbf{w}_i = (w_{ij}) = (0, \ldots, 0, \underbrace{k_{ni}^{-1}, \ldots, k_{ni}^{-1}}_{n_i}, 0, \ldots, 0)^\top\).
This unified form facilitates the expressions of martingale differences and plays an important role in calculating the limiting covariance. Similar techniques are also used in the proof of Theorem \ref{th-T}.

		\item 
When $p/n\rightarrow\infty$, since $\tilde{z}_n=O(c_n)$, we need to split  $\left(\mathbf{B}_{n 0}-\tilde{z}_n \mathbf{I}\right)^{-1}$ as:
	\begin{align}\label{spilt-Bn}
	\left(\mathbf{B}_{n 0}-\tilde{z}_n \mathbf{I}\right)^{-1}=\frac{1}{\tilde{z}_n}\left[\frac{1}{\sqrt{n p b_p}} \boldsymbol{\Sigma}_0^{\frac{1}{2}} \mathbf{Z}\left(\mathbf{A}_{ 0 }-z\mathbf{I}\right)^{-1} \mathbf{Z}^{\top} \boldsymbol{\Sigma}_0^{\frac{1}{2}}-\mathbf{I}_p\right],
	\end{align}
	where $\A_0$ is defined in \eqref{par-A0}.
	With the help of this identity, we can obtain some new moments of quadratic forms in the ultrahigh dimensional context. For instance, in this case,  we have
		\begin{align*}
		\E\left|\bar\s_i^{\top}\bar\s_i\right|^2=O(p^2/n^2),\quad 
		\E\left|c_n\bar\s_i^{\top}\left(\mathbf{B}_{n 0}-\tilde{z}_n \mathbf{I}\right)^{-1}\bSigma_0
		\left(\mathbf{B}_{n 0}-\tilde{z}_n \mathbf{I}\right)^{-1}
		\bar\s_i\right|^2=O(p^2/n^2),\\
		\left|\bmu_i^{\top}\bmu_i\right|^2=O(p/n),\quad 
		\E\left|c_n\bmu_i^{\top}\left(\mathbf{B}_{n 0}-\tilde{z}_n \mathbf{I}\right)^{-1}\bSigma_0
		\left(\mathbf{B}_{n 0}-\tilde{z}_n \mathbf{I}\right)^{-1}
		\bmu_i\right|^2=O(n/p).
		\end{align*}
		These moments also demonstrate different effects of $\bar\s_{i}$ and $\bmu_i$ and we need to handle them carefully when $p/n\rightarrow\infty$. While they are all $O(1)$ when $p/n\rightarrow c<\infty$. Actually,  $\|c_n\left(\mathbf{B}_{n 0}-\tilde{z}_n \mathbf{I}\right)^{-1}\bSigma_0
		\left(\mathbf{B}_{n 0}-\tilde{z}_n \mathbf{I}\right)^{-1}\|$ is bounded. While if we use this result directly,  we can just only obtain $\E\left|c_n \boldsymbol{\mu}_i^{\top}\left(\mathbf{B}_{n 0}-\tilde{z}_n \mathbf{I}\right)^{-1} \boldsymbol{\Sigma}_0\left(\mathbf{B}_{n 0}-\tilde{z}_n \mathbf{I}\right)^{-1} \boldsymbol{\mu}_i\right|^2=O(p/n)$ which is much larger than $O(n/p) $ obtained by using identity \eqref{spilt-Bn}.

		\item 
		The bound for moments of quadratic forms is crucial in our proof. To address this, a new lemma is derived as follows.
		\begin{lemma}\label{le-2.7}
		Let \(\mathbf{Y} = \left(Y_1, \ldots, Y_p\right)^T\), where \(Y_i\)'s are i.i.d. real random variables with mean 0 and variance 1. Let \(\mathbf{M}\) be a deterministic complex matrix. Then for any \(k \geq 2\), we have
			$$
			\E\left|\Y^T \mathbf{M} \Y-\operatorname{tr} \mathbf{M}\right|^k \leq \mathfrak{M}_k(\E|Y_1|^k)^2\left( \operatorname{tr} \mathbf{M} \mathbf{M}^*\right)^{k / 2}+\mathfrak{M}_k v_{2k}\sum_{h=1}^{p}\left|\e_h^{\top}\mathbf{M}\e_h\right|^k,
			$$
			where $\mathbf{M}^*$ denotes the complex conjugate transpose of $\mathbf{M}$ and $v_{2k}=\E(Y_1^{2k})$.
		\end{lemma}
	Compared with the conventional Lemma 2.7 in \cite{BS98}, which gives
			$$
		\E\left|\Y^T \mathbf{M} \Y-\operatorname{tr} \mathbf{M}\right|^k \leq \mathfrak{M}_k\left(v_4 \operatorname{tr} \mathbf{M} \mathbf{M}^*\right)^{k / 2}+\mathfrak{M}_k v_{2k} \operatorname{tr}\left(\mathbf{M} \mathbf{M}^*\right)^{k / 2},
		$$
		Lemma \ref{le-2.7}  provides a more refined bound when $p/n\to \infty$. For example, with the help of Lemma \ref{le-2.7} and by truncating  the underlying random variables \((z_{ijq})\) at \(\mathfrak{M}\epsilon_n(np)^{1/4}\), where \(\{\epsilon_n\}\) is a sequence of positive numbers decreasing to zero at a slow rate, we  have $\E|\s_1^{\top}\s_2\s_2^{\top}\s_1-\s_2^{\top}\bSigma_{0}\s_2|^4=O(p^4+\epsilon_n^8n^2p^3)$. However, if we use Lemma 2.7 in \cite{BS98} directly, the bound becomes $O(\epsilon_n^4np^5)$ which is significantly larger when $p/n\to \infty$.

	\end{itemize}

	
	\begin{supplement}
		\stitle{Supplementary Material of ``On spiked eigenvalues of a renormalized sample covariance matrix from multi-population"}
		\sdescription{This supplementary document contains the proofs of 
		Lemmas \ref{th-limit1-cinf}, \ref{th-limit1},\ref{le-tau-equation}-\ref{le-2.7} and Theorems \ref{le-hatA-cinf},\ref{th-num}, \ref{th-T}, \ref{le-hatA}.
}
	\end{supplement}

	\renewcommand\refname{REFERENCES}
	\bibliographystyle{imsart-nameyear.bst}
	\bibliography{MM}

\begin{thebibliography}{27}

\bibitem[\protect\citeauthoryear{Anderson}{2003}]{Anderson}
\begin{bbook}[author]
\bauthor{\bsnm{Anderson},~\bfnm{T.~W.}\binits{T.~W.}}
(\byear{2003}).
\btitle{An introduction to multivariate statistical analysis},
\bedition{third} ed.
\bseries{Wiley Series in Probability and Statistics}.
\bpublisher{Wiley-Interscience [John Wiley \& Sons], Hoboken, NJ}.
\bmrnumber{1990662}
\end{bbook}
\endbibitem

\bibitem[\protect\citeauthoryear{Bai, Miao and Pan}{2007}]{bai2007asymptotics}
\begin{barticle}[author]
\bauthor{\bsnm{Bai},~\bfnm{Zhidong.}\binits{Z.}},
  \bauthor{\bsnm{Miao},~\bfnm{Baiqi.}\binits{B.}} \AND
  \bauthor{\bsnm{Pan},~\bfnm{Guangming}\binits{G.}}
(\byear{2007}).
\btitle{On asymptotics of eigenvectors of large sample covariance matrix}.
\bjournal{Ann. Probab.}
\bvolume{35}
\bpages{1532--1572}.
\bdoi{10.1214/009117906000001079}
\bmrnumber{2330979}
\end{barticle}
\endbibitem

\bibitem[\protect\citeauthoryear{Bai and Silverstein}{1998}]{BS98}
\begin{barticle}[author]
\bauthor{\bsnm{Bai},~\bfnm{Zhidong.}\binits{Z.}} \AND
  \bauthor{\bsnm{Silverstein},~\bfnm{Jack~W.}\binits{J.~W.}}
(\byear{1998}).
\btitle{No eigenvalues outside the support of the limiting spectral
  distribution of large-dimensional sample covariance matrices}.
\bjournal{Ann. Probab.}
\bvolume{26}
\bpages{316--345}.
\bdoi{10.1214/aop/1022855421}
\bmrnumber{1617051}
\end{barticle}
\endbibitem

\bibitem[\protect\citeauthoryear{Bai and Silverstein}{2004}]{BS04}
\begin{barticle}[author]
\bauthor{\bsnm{Bai},~\bfnm{Zhidong.}\binits{Z.}} \AND
  \bauthor{\bsnm{Silverstein},~\bfnm{Jack~W.}\binits{J.~W.}}
(\byear{2004}).
\btitle{C{LT} for linear spectral statistics of large-dimensional sample
  covariance matrices}.
\bjournal{Ann. Probab.}
\bvolume{32}
\bpages{553--605}.
\bdoi{10.1214/aop/1078415845}
\bmrnumber{2040792}
\end{barticle}
\endbibitem

\bibitem[\protect\citeauthoryear{Bai and Silverstein}{2010}]{Bai2010a}
\begin{bbook}[author]
\bauthor{\bsnm{Bai},~\bfnm{Zhidong}\binits{Z.}} \AND
  \bauthor{\bsnm{Silverstein},~\bfnm{Jack~W.}\binits{J.~W.}}
(\byear{2010}).
\btitle{Spectral analysis of large dimensional random matrices},
\bedition{second} ed.
\bseries{Springer Series in Statistics}.
\bpublisher{Springer, New York}.
\bdoi{10.1007/978-1-4419-0661-8}
\bmrnumber{2567175}
\end{bbook}
\endbibitem

\bibitem[\protect\citeauthoryear{Bai and Yin}{1988}]{baiyin}
\begin{barticle}[author]
\bauthor{\bsnm{Bai},~\bfnm{Zhidong.}\binits{Z.}} \AND
  \bauthor{\bsnm{Yin},~\bfnm{Yongquan.}\binits{Y.}}
(\byear{1988}).
\btitle{Convergence to the semicircle law}.
\bjournal{Ann. Probab.}
\bvolume{16}
\bpages{863--875}.
\bmrnumber{929083}
\end{barticle}
\endbibitem

\bibitem[\protect\citeauthoryear{Bao}{2015}]{MR3416046}
\begin{barticle}[author]
\bauthor{\bsnm{Bao},~\bfnm{Z.}\binits{Z.}}
(\byear{2015}).
\btitle{On asymptotic expansion and central limit theorem of linear eigenvalue
  statistics for sample covariance matrices when {$N/M\to 0$}}.
\bjournal{Theory Probab. Appl.}
\bvolume{59}
\bpages{185--207}.
\bdoi{10.1137/S0040585X97T987089}
\bmrnumber{3416046}
\end{barticle}
\endbibitem

\bibitem[\protect\citeauthoryear{Billingsley}{1995}]{billingsley2008probability}
\begin{bbook}[author]
\bauthor{\bsnm{Billingsley},~\bfnm{Patrick}\binits{P.}}
(\byear{1995}).
\btitle{Probability and measure},
\bedition{third} ed.
\bseries{Wiley Series in Probability and Mathematical Statistics}.
\bpublisher{John Wiley \& Sons, Inc., New York}
\bnote{A Wiley-Interscience Publication}.
\bmrnumber{1324786}
\end{bbook}
\endbibitem

\bibitem[\protect\citeauthoryear{Chen and Pan}{2012}]{bbchen}
\begin{barticle}[author]
\bauthor{\bsnm{Chen},~\bfnm{Binbin.}\binits{B.}} \AND
  \bauthor{\bsnm{Pan},~\bfnm{Guangming.}\binits{G.}}
(\byear{2012}).
\btitle{Convergence of the largest eigenvalue of normalized sample covariance
  matrices when {$p$} and {$n$} both tend to infinity with their ratio
  converging to zero}.
\bjournal{Bernoulli}
\bvolume{18}
\bpages{1405--1420}.
\bdoi{10.3150/11-BEJ381}
\bmrnumber{2995802}
\end{barticle}
\endbibitem

\bibitem[\protect\citeauthoryear{Chen and Pan}{2015}]{cbbclt}
\begin{barticle}[author]
\bauthor{\bsnm{Chen},~\bfnm{Binbin}\binits{B.}} \AND
  \bauthor{\bsnm{Pan},~\bfnm{Guangming}\binits{G.}}
(\byear{2015}).
\btitle{C{LT} for linear spectral statistics of normalized sample covariance
  matrices with the dimension much larger than the sample size}.
\bjournal{Bernoulli}
\bvolume{21}
\bpages{1089--1133}.
\bdoi{10.3150/14-BEJ599}
\bmrnumber{3338658}
\end{barticle}
\endbibitem

\bibitem[\protect\citeauthoryear{Chen et~al.}{2011}]{2011A}
\begin{barticle}[author]
\bauthor{\bsnm{Chen},~\bfnm{Lin~S.}\binits{L.~S.}},
  \bauthor{\bsnm{Paul},~\bfnm{Debashis}\binits{D.}},
  \bauthor{\bsnm{Prentice},~\bfnm{Ross~L.}\binits{R.~L.}} \AND
  \bauthor{\bsnm{Wang},~\bfnm{Pei}\binits{P.}}
(\byear{2011}).
\btitle{A regularized {H}otelling's {$T^2$} test for pathway analysis in
  proteomic studies}.
\bjournal{J. Amer. Statist. Assoc.}
\bvolume{106}
\bpages{1345--1360}.
\bdoi{10.1198/jasa.2011.ap10599}
\bmrnumber{2896840}
\end{barticle}
\endbibitem

\bibitem[\protect\citeauthoryear{Horn and Johnson}{2012}]{horn2012matrix}
\begin{bbook}[author]
\bauthor{\bsnm{Horn},~\bfnm{Roger~A}\binits{R.~A.}} \AND
  \bauthor{\bsnm{Johnson},~\bfnm{Charles~R}\binits{C.~R.}}
(\byear{2012}).
\btitle{Matrix analysis}.
\bpublisher{Cambridge university press}.
\end{bbook}
\endbibitem

\bibitem[\protect\citeauthoryear{Jing et~al.}{2024}]{jingbingyi}
\begin{barticle}[author]
\bauthor{\bsnm{Jing},~\bfnm{Bingyi}\binits{B.}},
  \bauthor{\bsnm{Li},~\bfnm{Weiming}\binits{W.}},
  \bauthor{\bsnm{Xie},~\bfnm{Jiahui}\binits{J.}},
  \bauthor{\bsnm{Zhang},~\bfnm{Yangchun}\binits{Y.}} \AND
  \bauthor{\bsnm{Zhou},~\bfnm{Wang}\binits{W.}}
(\byear{2024}).
\btitle{On Convergence Rates of Spiked Eigenvalue Estimates: A General Study of
  Global and Local Laws in Sample Covariance Matrices}.
\bjournal{manuscript}.
\end{barticle}
\endbibitem

\bibitem[\protect\citeauthoryear{Li and Zhu}{2023}]{GMM}
\begin{barticle}[author]
\bauthor{\bsnm{Li},~\bfnm{Weiming}\binits{W.}} \AND
  \bauthor{\bsnm{Zhu},~\bfnm{Junpeng}\binits{J.}}
(\byear{2023}).
\btitle{C{LT} for spiked eigenvalues of a sample covariance matrix from
  high-dimensional {G}aussian mean mixtures}.
\bjournal{J. Multivariate Anal.}
\bvolume{193}
\bpages{Paper No. 105127, 22}.
\bdoi{10.1016/j.jmva.2022.105127}
\bmrnumber{4513868}
\end{barticle}
\endbibitem

\bibitem[\protect\citeauthoryear{Li et~al.}{2020}]{lihaoran}
\begin{barticle}[author]
\bauthor{\bsnm{Li},~\bfnm{Haoran}\binits{H.}},
  \bauthor{\bsnm{Aue},~\bfnm{Alexander}\binits{A.}},
  \bauthor{\bsnm{Paul},~\bfnm{Debashis}\binits{D.}},
  \bauthor{\bsnm{Peng},~\bfnm{Jie}\binits{J.}} \AND
  \bauthor{\bsnm{Wang},~\bfnm{Pei}\binits{P.}}
(\byear{2020}).
\btitle{An adaptable generalization of {H}otelling's {$T^2$} test in high
  dimension}.
\bjournal{Ann. Statist.}
\bvolume{48}
\bpages{1815--1847}.
\bdoi{10.1214/19-AOS1869}
\bmrnumber{4124345}
\end{barticle}
\endbibitem

\bibitem[\protect\citeauthoryear{Lin et~al.}{2024}]{lin2024asymptotic}
\begin{barticle}[author]
\bauthor{\bsnm{Lin},~\bfnm{Zeqin}\binits{Z.}},
  \bauthor{\bsnm{Pan},~\bfnm{Guangming}\binits{G.}},
  \bauthor{\bsnm{Zhao},~\bfnm{Peng}\binits{P.}} \AND
  \bauthor{\bsnm{Zhou},~\bfnm{Jia}\binits{J.}}
(\byear{2024}).
\btitle{Asymptotic distribution of spiked eigenvalues in the large
  signal-plus-noise models}.
\bjournal{arXiv preprint arXiv:2401.11672}.
\end{barticle}
\endbibitem

\bibitem[\protect\citeauthoryear{Liu}{2020}]{yiminghigh}
\begin{bphdthesis}[author]
\bauthor{\bsnm{Liu},~\bfnm{Yiming.}\binits{Y.}}
(\byear{2020}).
\btitle{High dimensional clustering for mixture models},
\btype{PhD thesis},
\bpublisher{School of Physical and Mathematical Science, Nanyang Technological
  University}
\bnote{Available at:{\color{blue} \url{https://hdl.handle.net/10356/142941}}}.
\bdoi{10.32657/10356/142941}
\end{bphdthesis}
\endbibitem

\bibitem[\protect\citeauthoryear{Liu et~al.}{2023}]{liu2023asymptotic}
\begin{barticle}[author]
\bauthor{\bsnm{Liu},~\bfnm{Xiaoyu}\binits{X.}},
  \bauthor{\bsnm{Liu},~\bfnm{Yiming}\binits{Y.}},
  \bauthor{\bsnm{Pan},~\bfnm{Guangming}\binits{G.}},
  \bauthor{\bsnm{Zhang},~\bfnm{Lingyue}\binits{L.}} \AND
  \bauthor{\bsnm{Zhang},~\bfnm{Zhixiang}\binits{Z.}}
(\byear{2023}).
\btitle{Asymptotic properties of spiked eigenvalues and eigenvectors of
  signal-plus-noise matrices with their applications}.
\bjournal{arXiv preprint arXiv:2310.13939}.
\end{barticle}
\endbibitem

\bibitem[\protect\citeauthoryear{Mar\v{c}enko and Pastur}{1967}]{MP67}
\begin{barticle}[author]
\bauthor{\bsnm{Mar\v{c}enko},~\bfnm{V.~A.}\binits{V.~A.}} \AND
  \bauthor{\bsnm{Pastur},~\bfnm{L.~A.}\binits{L.~A.}}
(\byear{1967}).
\btitle{Distribution of eigenvalues in certain sets of random matrices}.
\bjournal{Mat. Sb. (N.S.)}
\bvolume{72}
\bpages{507--536}.
\end{barticle}
\endbibitem

\bibitem[\protect\citeauthoryear{Pan}{2014}]{MR3189088}
\begin{barticle}[author]
\bauthor{\bsnm{Pan},~\bfnm{Guangming.}\binits{G.}}
(\byear{2014}).
\btitle{Comparison between two types of large sample covariance matrices}.
\bjournal{Ann. Inst. Henri Poincar\'{e} Probab. Stat.}
\bvolume{50}
\bpages{655--677}.
\bdoi{10.1214/12-AIHP506}
\bmrnumber{3189088}
\end{barticle}
\endbibitem

\bibitem[\protect\citeauthoryear{Pan and Zhou}{2008}]{PZ08}
\begin{barticle}[author]
\bauthor{\bsnm{Pan},~\bfnm{Guangming.}\binits{G.}} \AND
  \bauthor{\bsnm{Zhou},~\bfnm{Wang}\binits{W.}}
(\byear{2008}).
\btitle{Central limit theorem for signal-to-interference ratio of reduced rank
  linear receiver}.
\bjournal{Ann. Appl. Probab.}
\bvolume{18}
\bpages{1232--1270}.
\bdoi{10.1214/07-AAP477}
\bmrnumber{2418244}
\end{barticle}
\endbibitem

\bibitem[\protect\citeauthoryear{Pan and Zhou}{2011}]{pan-T}
\begin{barticle}[author]
\bauthor{\bsnm{Pan},~\bfnm{Guangming.}\binits{G.}} \AND
  \bauthor{\bsnm{Zhou},~\bfnm{Wang.}\binits{W.}}
(\byear{2011}).
\btitle{Central limit theorem for {H}otelling's {$T^2$} statistic under large
  dimension}.
\bjournal{Ann. Appl. Probab.}
\bvolume{21}
\bpages{1860--1910}.
\bdoi{10.1214/10-AAP742}
\bmrnumber{2884053}
\end{barticle}
\endbibitem

\bibitem[\protect\citeauthoryear{Qiu, Li and Yao}{2023}]{qiujiaxin}
\begin{barticle}[author]
\bauthor{\bsnm{Qiu},~\bfnm{Jiaxin}\binits{J.}},
  \bauthor{\bsnm{Li},~\bfnm{Zeng}\binits{Z.}} \AND
  \bauthor{\bsnm{Yao},~\bfnm{Jianfeng}\binits{J.}}
(\byear{2023}).
\btitle{Asymptotic normality for eigenvalue statistics of a general sample
  covariance matrix when {$p/n \to\infty $} and applications}.
\bjournal{Ann. Statist.}
\bvolume{51}
\bpages{1427--1451}.
\bdoi{10.1214/23-aos2300}
\bmrnumber{4630955}
\end{barticle}
\endbibitem

\bibitem[\protect\citeauthoryear{Wang and Paul}{2014}]{wanglili}
\begin{barticle}[author]
\bauthor{\bsnm{Wang},~\bfnm{Lili}\binits{L.}} \AND
  \bauthor{\bsnm{Paul},~\bfnm{Debashis}\binits{D.}}
(\byear{2014}).
\btitle{Limiting spectral distribution of renormalized separable sample
  covariance matrices when {$p/n\to 0$}}.
\bjournal{J. Multivariate Anal.}
\bvolume{126}
\bpages{25--52}.
\bdoi{10.1016/j.jmva.2013.12.015}
\bmrnumber{3173080}
\end{barticle}
\endbibitem

\bibitem[\protect\citeauthoryear{Xu, Liu and Yao}{2023}]{MR4606323}
\begin{barticle}[author]
\bauthor{\bsnm{Xu},~\bfnm{Yuyang}\binits{Y.}},
  \bauthor{\bsnm{Liu},~\bfnm{Zhonghua}\binits{Z.}} \AND
  \bauthor{\bsnm{Yao},~\bfnm{Jianfeng}\binits{J.}}
(\byear{2023}).
\btitle{An eigenvalue ratio approach to inferring population structure from
  whole genome sequencing data}.
\bjournal{Biometrics}
\bvolume{79}
\bpages{891--902}.
\bmrnumber{4606323}
\end{barticle}
\endbibitem

\bibitem[\protect\citeauthoryear{Yao, Zheng and Bai}{2015}]{yao2015sample}
\begin{bbook}[author]
\bauthor{\bsnm{Yao},~\bfnm{Jianfeng}\binits{J.}},
  \bauthor{\bsnm{Zheng},~\bfnm{Shurong}\binits{S.}} \AND
  \bauthor{\bsnm{Bai},~\bfnm{Zhidong}\binits{Z.}}
(\byear{2015}).
\btitle{Large sample covariance matrices and high-dimensional data analysis}.
\bseries{Cambridge Series in Statistical and Probabilistic Mathematics}
\bvolume{39}.
\bpublisher{Cambridge University Press, New York}.
\bdoi{10.1017/CBO9781107588080}
\bmrnumber{3468554}
\end{bbook}
\endbibitem

\bibitem[\protect\citeauthoryear{Zheng, Bai and
  Yao}{2015}]{zheng2015substitution}
\begin{barticle}[author]
\bauthor{\bsnm{Zheng},~\bfnm{Shurong}\binits{S.}},
  \bauthor{\bsnm{Bai},~\bfnm{Zhidong}\binits{Z.}} \AND
  \bauthor{\bsnm{Yao},~\bfnm{Jianfeng}\binits{J.}}
(\byear{2015}).
\btitle{Substitution principle for {CLT} of linear spectral statistics of
  high-dimensional sample covariance matrices with applications to hypothesis
  testing}.
\bjournal{Ann. Statist.}
\bvolume{43}
\bpages{546--591}.
\bdoi{10.1214/14-AOS1292}
\bmrnumber{3316190}
\end{barticle}
\endbibitem

\end{thebibliography}

\end{document}